\documentclass[10pt,leqno]{gtpart}
        
\usepackage{hyperref}
\usepackage{calc}
\usepackage{enumerate}
\usepackage[arrow,curve,matrix,tips,2cell]{xy}
  \SelectTips{eu}{10} \UseTips
  \UseAllTwocells

\title{Equivariant covers for hyperbolic groups}
       
\author{Arthur Bartels}
\givenname{Arthur}
\surname{Bartels}
\address{Department of Mathematics\\
         Imperial College London\\
         London SW7 2AZ, UK}
\email{a.bartelsa@imperial.ac.uk}
\urladdr{http://www.ma.ic.ac.uk/~abartels/}

\author{Wolfgang L\"uck}
\givenname{Wolfgang}
\surname{L\"uck}
\address{Westf\"alische Wilhelms-Universit\"at M\"unster\\
               Mathematisches Institut\\
               Einsteinstr.~62,
               D-48149 M\"unster, Germany}
\email{lueck@math.uni-muenster.de}
\urladdr{http://www.math.uni-muenster.de/u/lueck}

\author{Holger Reich}
\givenname{Holger}
\surname{Reich}
\address{Heinrich-Heine-Universit\"at D\"usseldorf\\
         Mathematisches Institut\\
         Universit{\"a}tsstr.~1,
         D-40225 D\"usseldorf, Germany}
\email{reich@math.uni-duesseldorf.de}
\urladdr{http://reh.math.uni-duesseldorf.de/\%7Ereich/}


\volumenumber{}
\issuenumber{}
\publicationyear{}
\papernumber{}
\startpage{}
\endpage{}
\doi{}
\MR{}
\Zbl{}
\received{}
\revised{}
\accepted{}
\published{}
\publishedonline{}
\proposed{}
\seconded{}
\corresponding{}
\editor{}
\version{}

\DeclareMathAlphabet{\matheurm}{U}{eur}{m}{n}

\newcommand{\per}{\text{per}}

\DeclareMathOperator{\diam}{diam}
\DeclareMathOperator{\id}{id}

\DeclareMathOperator{\Isom}{Isom}

\newcommand{\Fin}{{\mathcal{F}\text{in}}}

\newcommand{\VCyc}{{\mathcal{VC}\text{yc}}}

  \newcommand{\IN}{\mathbb{N}}

  \newcommand{\IR}{\mathbb{R}}

  \newcommand{\calc}{\mathcal{C}}
  \newcommand{\cald}{\mathcal{D}}
  \newcommand{\cale}{\mathcal{E}}
  \newcommand{\calf}{\mathcal{F}}

  \newcommand{\call}{\mathcal{L}}
  
  \newcommand{\caln}{\mathcal{N}}

  \newcommand{\calu}{\mathcal{U}}
  \newcommand{\calv}{\mathcal{V}}
  \newcommand{\calw}{\mathcal{W}}
  \newcommand{\calx}{\mathcal{X}}

\newcommand{\x}{{\times}}

\newcommand{\e}{\epsilon}

\newcommand{\FS}{{\mathit{FS}}}

\newcommand{\dd}{{\partial}}
\newcommand{\sjX}{{\raise0.12ex\hbox{$\scriptscriptstyle\smallsmile$}
                    \mspace{-12mu}\raise0.3ex\hbox{$*$}\bar{X}}}
\newcommand{\osjX}{{\raise0.2ex\hbox{$*$}\bar{X}}}
\newcommand{\sjnobarXx}{{\raise-0.3ex\hbox{$\scriptscriptstyle{x_0}$}
                        \mspace{-15mu}{\raise0.22ex\hbox
                        {$\circ$}\mspace{-12mu}\raise0.22ex\hbox{$*$}X} }}
\newcommand{\sjnobarX}{{\raise0.2ex\hbox{$\circ\mspace{-9mu}*$}\mspace{-1.5mu}X}}
\newcommand{\good}{{\mathit{good}}}
\newcommand{\bad}{{\mathit{bad}}}

\newcommand{\stern}{{\mathit{star}}}

\newenvironment{numberlist}
  {\begin{list}{}%
   {%
    \setlength{\leftmargin}{\labelwidth+\labelsep}%
   }%
  }%
  {\end{list}}

\newcommand{\peri}[2]{\operatorname{per}^{#1}_{#2}}

\theoremstyle{plain}
\newtheorem{theorem}{Theorem}[section]

\newtheorem{lemma}[theorem]{Lemma}

\newtheorem{proposition}[theorem]{Proposition}

\newtheorem{convention}[theorem]{Convention}

\theoremstyle{definition}
\newtheorem{definition}[theorem]{Definition}

\newtheorem{remark}[theorem]{Remark}

\theoremstyle{remark}

\renewcommand{\theenumi}{(\roman{enumi})}

\makeatletter\let\c@equation=\c@theorem\makeatother
\renewcommand{\theequation}{\thetheorem}

\begin{document}

\begin{abstract}
We prove an equivariant version of the
fact that word-hyperbolic groups have finite asymptotic
dimension.
This is important in connection with our
forthcoming proof of the Farrell-Jones
conjecture for $K_*(RG)$ for every
word-hyperbolic group $G$ and every coefficient ring $R$.
\end{abstract}

\maketitle


\section{Introduction}
\label{sec:intro}

The asymptotic dimension of a metric space $X$
was introduced by Gromov in \cite[p.29]{Gromov(1993)}.
It can be defined as the smallest number $N$ such that for
every $\alpha > 0$ there exists an open cover $\calu$ of
$X$ with the following properties:
\begin{itemize}
\item $\dim \calu \leq N$;
\item The Lebesgue number of $\calu$ is at least $\alpha$,
      i.e., for every $x \in X$ there is $U \in \calu$
      such that $x^\alpha \subseteq U$, where $x^\alpha$
      is the open ball of radius $\alpha$ around $x$;
\item The members of $\calu$ have uniformly bounded diameters.
\end{itemize}
Recall that a cover $\calu$ is of dimension $\leq N$
if every $x \in X$ is contained in no more then
$N+1$ members of $\calu$.
The asymptotic dimension of a finitely generated group
is its asymptotic dimension as a metric space with
respect to any word metric.
An important result of Yu
\cite{Yu(1998a)} asserts
that the Novikov conjecture holds for groups of
finite asymptotic dimension.
This can be viewed as an injectivity result for the
assembly map in $L$-theory (after inverting $2$).
Further injectivity results for assembly maps for groups
with finite asymptotic dimension can be found in
\cite{Bartels-finite-asymp-dim},
\cite{Carlsson-Goldfarb-finite-asymptotic}
and \cite{Bartels-Rosenthal-asymptotic}.
On the other hand no surjectivity statement of
assembly maps is known for all groups of
finite asymptotic dimension and this is very much
related to the absence of any equivariance condition
for the cover $\calu$ as above.

\begin{definition}
\label{def:calf-cover}
Let $G$ be a group and $Z$ be a $G$-space.
Let $\calf$ be a collection of subgroups of $G$.
An open cover $\calu$ of $Z$ is called an
{\em $\calf$-cover} if the following two conditions
are satisfied.
\begin{enumerate}
\item For $g \in G$ and $U \in \calu$ we have
      either $g(U) = U$ or $g(U) \cap U = \emptyset$;
\item For $g \in G$ and $U \in \calu$ we have $g(U) \in \calu$;
\item For $U \in \calu$ the subgroup
      $G_U := \{ g \in G \; | \; g(U) = U \}$
      is a member of $\calf$.
\end{enumerate}
\end{definition}

Let $G$ be a word-hyperbolic group.
Fix a set of generators $S$.
Let $d_G$ be the word metric on $G$ with respect to $S$.
Let $X$ be a hyperbolic complex with an isometric
$G$-action in the sense of Mineyev
\cite{Mineyev-flows-and-joins}, see Section~
\ref{subsec:Hyperbolic complexes, double difference and Gromov product}.
Let $\dd X$ be the Gromov boundary of $X$.
(This boundary can be described as a quotient of the set of geodesic rays
in $X$, where two such rays are identified if they are asymptotic,
\cite[III.H.3]{Bridson-Haefliger-Buch}.)
Let $\overline{X} := X \cup \dd X$ be the compactification of
$X$, \cite[III.H.3]{Bridson-Haefliger-Buch}).
Let $\VCyc$ denote the collection of virtually cyclic
subgroups of $G$, that is of subgroups that have a
cyclic subgroup of finite index.
The following is our main result and should be thought of
as an equivariant version of the (much easier) fact that
hyperbolic groups have finite asymptotic dimension,
\cite[p.~31]{Gromov(1993)}, \cite{Roe-asym-dim-for-hyperbolic}.

\begin{theorem}
\label{thm:equivariant-cover-for-hyperbolic-group}
Let $G$ be word-hyperbolic and let $X$ be a hyperbolic complex.
Suppose that there is a simplicial proper cocompact $G$-action
on $X$.
Equip $G \times \overline{X}$ with the diagonal $G$-action.
Then there exists a natural number
$N = N(G,\overline{X})$ depending only on $G$ and
$\overline{X}$ such that the following holds:
For every $\alpha > 0$ there exists an open
$\VCyc$-cover $\calu$ of $G \times \overline{X}$ satisfying
\begin{enumerate}
\item $\dim(\calu) \le N$;
\item For $g_0 \in G$ and $c \in \overline{X}$ there exists
      $U \in \calu$ such that
      $g_0^{\alpha} \times \{ c \} \subseteq U$,
      where $g_0^{\alpha}$ is the open ball with center
      $g_0$ and radius $\alpha$ with respect to the
      word metric $d_G$;
\item $G\backslash \calu$ is finite.
\end{enumerate}
\end{theorem}

This result plays an important role in our proof
of the Farrell-Jones conjecture for $K_*(RG)$ for
every word-hyperbolic group $G$ and every coefficient ring $R$,
\cite{Bartels-Lueck-Reich-hyperbolic}.

The conclusion of
Theorem~\ref{thm:equivariant-cover-for-hyperbolic-group}
is formally similar to
the definition of finite asymptotic dimension discussed above.
The price we have to pay for the equivariance of the cover
$\calu$ is the space $\overline{X}$.
For the application it will be very important that
$\overline{X}$ is compact.
(If we replace $\overline{X}$ by a finite dimensional
$G$-$CW$-complex all whose isotropy groups lie in $\VCyc$,
then the conclusion follows easily from the fact that $G$
has finite asymptotic dimension.)
The members of $\calu$ are only large in the $G$-coordinate;
in the $\overline{X}$-coordinate they may be very small.
Similar covers have been used in a slightly different situation
where $\overline{X}$ is replaced by a probability space
with a measure preserving action of $G$, compare
\cite[p.300]{Gromov(1999a)},
\cite{Sauer-amenable-covers}.
It would be interesting to know if there is a version of
Theorem~\ref{thm:equivariant-cover-for-hyperbolic-group}
in this situation.

It seems reasonable to hope that the class of groups $G$ for
which there is a compact $G$-space $\overline{X}$ such that
the conclusion of
Theorem~\ref{thm:equivariant-cover-for-hyperbolic-group}
holds is bigger than the class of hyperbolic groups.

The proof of
Theorem~\ref{thm:equivariant-cover-for-hyperbolic-group}
is quite involved and uses a generalization of techniques
used and developed by Farrell-Jones in
\cite{Farrell-Jones-dynamics-I}.
Firstly, we study flows on metric spaces and prove
the existence of long and thin covers,
see Theorem~\ref{thm:long-thin-cover}.
This generalizes the long and thin cells from
\cite[Proposition~7.2]{Farrell-Jones-dynamics-I}.
Secondly, we use a variant $\FS(X)$ of
Mineyev's half open symmetric join $\sjX$
\cite{Mineyev-flows-and-joins}.
This space is a substitute for the sphere bundle
of a negatively curved manifold and is equipped with
a flow $\phi_\tau$
(corresponding to the geodesic flow on the sphere bundle).
In Theorem~\ref{thm:final-flow-estimate}
we improve upon Mineyev's flow estimate
\cite[Theorem~57 on p.~468]{Mineyev-flows-and-joins}.
The required cover is then produced by pulling back
a long and thin cover of $\FS(X)$ by the composition
of the flow $\phi_\tau$ for large $\tau$ with an
embedding $G \x \overline{X} \to \FS(X)$.
A more detailed discussion follows in
Sections~\ref{subsec:intro-long-thin}
and \ref{subsec:intro-flow}.


\subsection{Long thin covers}
\label{subsec:intro-long-thin}

The existence of long thin covers will be proven in
the following situation.

\begin{convention}
\label{conv:for-long-and-thin-covers}
Let
\begin{itemize}
\item $G$ be a discrete group;
\item $X$ be a metrizable topological space
      with a proper cocompact $G$-action on $X$;
\item $\Phi \colon X \times \IR \to X$ be a flow.
\end{itemize}
Assume that the following conditions
are satisfied:
\begin{itemize}
\item $\Phi$ is $G$-equivariant;
\item The number of closed orbits, which are not
      stationary and whose period is $ \le C$,
      of the flow induced on $G\backslash X$ is finite
      for every $C > 0$;
\item $X-X^{\IR}$ is locally connected (notation explained below);
\item If we put
      \begin{eqnarray*}
      k_G & := & \sup\{|H| \mid H \subseteq G
          \text{ subgroup with finite order } |H|\};
      \\
      d_X & := & \dim(X-X^{\IR}),
      \end{eqnarray*}
      then $k_G < \infty$ and $d_X < \infty$.
\end{itemize}
\end{convention}

Recall that a $G$-action is \emph{proper} if for
every $x \in X$ there exists an
open neighborhood $U$ such that the set
$\{g \in G \mid U \cap gU \not= \emptyset\}$ is finite.
Recall that $X$ is \emph{locally connected} if for
each  $x \in X$ and each open neighborhood $U$ of $x$
we can find a connected open neighborhood $U'$ of $x$
with $U' \subseteq U$.
A $G$-space $X$ is called \emph{cocompact}
if $G\backslash X$ is compact.
The \emph{dimension of a collection of subsets
$\{U_i \mid i \in I\}$} is $\le d$,
if every point is contained in at most $d+1$ members
of the $U_i$.
The \emph{covering dimension of a space $X$} is $\le d$
if every open covering has an open refinement
whose dimension is less or equal to $d$.
One may replace the covering dimension $d_X$
of $X - X^\IR$ appearing above
by the supremum of the covering dimensions of compact
subsets of $X - X^{\IR}$.
Recall that an \emph{equivariant flow}
$\Phi \colon \IR \times X \to X$ is a continuous
$\IR$-action, such that $\Phi_{\tau}(gx) = g\Phi_{\tau}(x)$
holds for all $g \in G$, $\tau \in \IR$ and $x \in X$.
We denote by $X^{\IR}$ the $\IR$-fixed point set,
i.e., the set of points $x \in X$ for which
$\Phi_{\tau}(x) = x$  for all $\tau \in \IR$.
The period  of a closed orbit of $\Phi$ which is not
stationary is the smallest number $\tau > 0$ such that
$\Phi_{\tau}(x) = x$ holds for all $x$ in this orbit.

The following is our main result in this situation.

\begin{theorem}
\label{thm:long-thin-cover}
There exists  a natural number
$N$ depending only on $k_G$, $d_X$ and the action of $G$ on
an arbitrary small neighborhood of $X^\IR$
such that for every $\alpha > 0$ there is an
$\VCyc$-cover $\calu$ of $X$ with the following two
properties:
\begin{enumerate}
\item $\dim \calu \leq N$;
\item For every $x \in X$ there exists $U \in \calu$ such
      that
      \[
      \Phi_{[-\alpha,\alpha]} (x) :=
      \{ \Phi_\tau(x) \; | \; \tau \in [-\alpha,\alpha] \}
                                   \subseteq U;
      \]
\item $G\backslash \calu$ is finite.
\end{enumerate}
\end{theorem}

The main difference between
Theorem~\ref{thm:long-thin-cover} and
\cite[Proposition~7.2]{Farrell-Jones-dynamics-I}
is that we deal with
metric spaces rather than manifolds.
This requires a different type of general position argument
(compare Section~\ref{sec:general-position-metric-spaces})
and forces us to work with open covers rather than
cell structures.
While cell structures of a manifold are automatically
finite dimensional, in our situation more care is needed to
establish the bound on the dimension of $\calu$ and our
bound is much larger then the dimension of the metric space $X$.
Finally, we deal with an honest proper action and do not require
a torsion free subgroup of finite index, as is used in
\cite{FJ-iso-conj}.

The proof of
Theorem~\ref{thm:long-thin-cover}
will be given in
Section~\ref{sec:construction-long-VCyc-cover-X}
and depends on Sections~\ref{sec:boxes},
\ref{sec:general-position-metric-spaces} and
\ref{sec:covering-G->gamma}.


\subsection{The flow space}
\label{subsec:intro-flow}

Let $G$ be a hyperbolic group.
Fix a set of generators $S$.
Let $d_G$ be the word metric on $G$ with respect to $S$.
Let $X$ be a hyperbolic complex
and $\overline{X} = X \cup \dd X$ be
its compactification as before.
Assume that $G$ acts isometrically on $X$.
In Section~\ref{sec:Mineyevs_flow_space}
we introduce the metric space  $(\FS(X),d_{\FS})$.
This space is equipped with an isometric
$G$-action and a $G$-equivariant flow $\phi_\tau$.

Our main flow estimate is the following.

\begin{theorem}
\label{thm:final-flow-estimate}
There exists a continuous $G$-equivariant
(with respect to the diagonal $G$-action on the source)
map $j \colon G \x \overline{X} \to \FS(X)$
such that for every $\alpha > 0$ there exists
a number $\beta = \beta(\alpha)$
such that the following holds:

If $g, h \in G$ with $d_G(g,h) \leq \alpha$
and $c \in \overline{X}$ then there is $\tau_0 \in \IR$
with $|\tau_0| \leq \beta$ such that for all $\tau \in \IR$
\[
d_{\FS} ( \phi_\tau j(g,c),
    \phi_{\tau + \tau_0} j(h,c) ) \leq f_\alpha(\tau).
\]
Here $f_\alpha \colon \IR \to [0,\infty)$ is a function
that depends only on $\alpha$ and has
the property that
$\lim_{\tau \to \infty} f_\alpha(\tau) = 0$.
\end{theorem}

An important ingredient of the proof of this result is
Theorem~\ref{the:exponential_estimate_for_d_{overline{FS}_x_0}}
which is an improvement of Mineyev's
\cite[Theorem~57 on p.~468]{Mineyev-flows-and-joins}.
The main differences are  that we consider points not
necessary on the same horosphere, and that we consider
the action of the flow $\phi_\tau$ and not translation
by length.
In addition, Mineyev's estimate is in terms of a
pseudo-metric, not in terms of the metric $d_\FS$.

In order to apply Theorem~\ref{thm:long-thin-cover}
to $\FS(X)$ we need further properties of the flow space
and $G$.

\begin{proposition}
\label{prop:further-properties-of-FS}
$ $
\begin{enumerate}
\item \label{prop:further-properties-of-FS:finite-subgroups}
      The order of finite subgroups in $G$ is bounded.
\item \label{prop:further-properties-of-FS:connected+dimension}
      $\FS(X) - \FS(X)^\IR$ is locally connected
      and has finite covering dimension.
\item \label{prop:further-properties-of-FS:proper+cocompact}
      If the action of $G$ on $X$ is cocompact and proper,
      then action of $G$ on $\FS(X)$ is also cocompact
      and proper.
\item \label{prop:further-properties-of-FS:closed-orbits}
      If the action of $G$ on $X$ is cocompact and proper,
      then the number of closed orbits, which are not
      stationary and whose period is $ \le C$,
      of the flow induced on $G \backslash \FS(X)$ is finite
      for every $C > 0$.
\end{enumerate}
\end{proposition}

The proof of
Theorem~\ref{thm:final-flow-estimate} will be given in
Section~\ref{sec:The_flow_estimates_for_the_map_iota}
and depends only on Sections~\ref{sec:Mineyevs_flow_space} and
\ref{sec:flow_estimates}.
The proof of Proposition~\ref{prop:further-properties-of-FS}
will be given in
Section~\ref{sec:Properness_of_the_induced_action}
and depends only on
Sections~\ref{sec:boxes} and~\ref{sec:Mineyevs_flow_space}. 


\subsection{Construction of the cover}

Using the results from Sections~\ref{subsec:intro-long-thin}
and \ref{subsec:intro-flow}
we can now give the proof of
Theorem~\ref{thm:equivariant-cover-for-hyperbolic-group}.
During this proof we will use the following notation:
if $A$ is a subset of a metric space $Z$ and $\delta > 0$,
then $A^\delta$ denotes the set of all points $z \in Z$ for
which $d(z,A) < \delta$, compare
Definition~\ref{def:thickenings-in-metric-space}.

\begin{proof}
Consider any $\alpha > 0$.
Let $\beta = \beta(\alpha)$ be the number
appearing in Theorem~\ref{thm:final-flow-estimate}.
It follows from Proposition~\ref{prop:further-properties-of-FS}
that Theorem~\ref{thm:long-thin-cover} can be applied
to $\FS(X)$.
Thus there is a number $N$ (independent of $\alpha$)
such that there exists an $\VCyc$-cover
$\calv$ of $\FS(X)$ of dimension no more than $N$
with the following property:
For every $\xi \in \FS(X)$ there exists $V_{\xi} \in \calv$
such that
\[
\phi_{[-2\beta,2\beta]}(\xi)
   = \{\phi_{\tau}(\xi) \mid \tau \in [-2\beta,2\beta]\}
                              \subseteq V_{\xi}.
\]
Since $\phi_{[-2\beta,2\beta]}(\xi)$ is compact,
$V_{\xi}$ is open and
$\phi_{[-2\beta,2\beta]}(\xi) \subseteq V_{\xi}$,
we can find
$\delta_{\xi} > 0$ (depending on $\xi$ and $\beta$, $V_\xi$) such that
\[
\left(\phi_{[-2\beta,2\beta]}(\xi)\right)^{\delta_{\xi}}
      \subseteq V_{\xi}.
\]
Because $G$ acts by isometries,
we can arrange that $\delta_{\xi} = \delta_{g \xi}$
holds for all $g \in G$.
In particular we get
$g \cdot \left(\phi_{[-2\beta,2\beta]}
         (\xi) \right)^{\delta_{\xi}} =
\left(\phi_{[-2\beta,2\beta]}(g \xi) \right)
                        ^{\delta_{g \xi}}$.
For $\xi \in \FS(X)$ pick $\e_\xi > 0$ such that
\[
0 < e^\beta \e_\xi < \delta_\xi / 2.
\]
Again we arrange that $\e_{g\xi} = \e_\xi$ holds for all $g \in G$.
Obviously the collection
\[
\left\{\left(\phi_{[-\beta,\beta]}(\xi)\right)
  ^{\e_{\xi}} \mid \xi \in \FS(X)\right\}
\]
is an open covering of $\FS(X)$.
Since $G$ acts cocompactly, we can find finitely many points
$\xi_i$ for $i = 0,1,2, \ldots, I$ for some positive
natural number $I$ such that the $G$-cofinite collection
\[
\left\{\left(\phi_{[-\beta,\beta]}
            (g \xi_i) \right)
            ^{\e_{g \xi_i}} \mid g \in G,
                            i \in \{0,1,2 \ldots, I\}\right\}
\]
is an open covering of $\FS(X)$.
Consider $\xi \in \FS(X)$.
Then we can find $i = i(\xi) \in \{0,1,2 \ldots, I\}$
and $g = g(\xi) \in G$ such that
$\xi \in \left(\phi_{[-\beta,\beta]}
               (g  \xi_i) \right)^{\e_{g \xi_i}}$.
In particular, there is $\tau \in [-\beta,\beta]$ such that
$d_\FS(\xi, \phi_{\tau}(g \xi_i)) < \e_{g \xi_i}$.
Let
\[
 \delta :=
 \min\{\delta_{\xi_i} / 2 \mid i = 0,1,2 \ldots, I\}.
\]
Consider
$\zeta \in \left(\phi_{[-\beta,\beta]}(\xi)\right)^{\delta}$.
Choose $\sigma \in [-\beta,\beta]$ satisfying
$d_{\FS}(\zeta,\phi_{\sigma}(\xi)) < \delta$.
In the following estimate we will use
Lemma~\ref{lem:_exponential_estimate_on_flow}.
(In this lemma the more careful notation
$d_{\FS,x_0}$ is used for $d_\FS$.)
\begin{eqnarray*}
d_{\FS}(\zeta, \phi_{\sigma + \tau}(g \xi_i))
& \le &
d_{\FS}(\zeta,\phi_{\sigma}(\xi)) +
      d_{\FS}(\phi_{\sigma}(\xi), \phi_{\sigma + \tau}(g \xi_i))
\\
& < &
\delta + e^{|\sigma|} \cdot
     d_{\FS}(\xi,\phi_{\tau}(g \xi_i))
\\
& < & \delta +  e^{\beta} \cdot \e_{g \xi_i}
\\
& < &  \delta_{g \xi_i}
\end{eqnarray*}
Since $\sigma + \tau \in [-2\beta,2\beta]$, this implies
\[
\left(\phi_{[-\beta,\beta]}(\xi) \right)^{\delta}
~ \subseteq ~
\left(\phi_{[-2\beta,2\beta]}(g \xi_i) \right)
                         ^{ \delta_{g \xi_i}}
~ \subseteq ~ V_{g \xi_i}.
\]
Thus we have found $\delta > 0$ such that for every
$\xi \in \FS(X)$ there exists $V_{\xi} \in \calv$
such that
\begin{eqnarray}
\left(\phi_{[-\beta,\beta]}(\xi)\right)^{\delta}
               & \subseteq & V_{\xi}.
\label{eq:beta^delta-in-V}
\end{eqnarray}

We will construct the desired open covering
$\calu$ of $G \times \overline{X}$ by pulling back
$\calv$ with the composition
\[
G \times \overline{X} \xrightarrow{j}
\FS(X) \xrightarrow{\phi_{\tau}} \FS(X)
\]
for an appropriate real number $\tau$, where
$j$ is the map from Theorem~\ref{thm:final-flow-estimate}.
Obviously $\calu$ has for every choice of $\tau$
all the desired properties
except for the property that there exists
$U_{(g_0,c)} \in \calu$ such that
$g_0^{\alpha} \times \{ c \} \subseteq U_{(g_0,c)}$
for every $c \in \overline{X}$
and every $g_0 \in G$.

We conclude from
Theorem~\ref{thm:final-flow-estimate}
for $\tau \in \IR$ and the function $f_\alpha$
appearing in Theorem~\ref{thm:final-flow-estimate}
\[
\phi_{\tau} \circ j(g ,c) ~ \in
\left(\phi_{[-\beta,\beta]}
        (\phi_{\tau} \circ j (g_0 ,c)) \right)^{f_\alpha(\tau)}
\]
for all  $c \in \overline{X}$ and all $g \in G$
with $d_G(g_0,g) < \alpha$.
By Theorem~\ref{thm:final-flow-estimate} there is $\tau$
such that $f_\alpha(\tau) < \delta$.
For such a choice of $\tau$ we conclude
from \eqref{eq:beta^delta-in-V} that
\[
\phi_{\tau} \circ j(g ,c) ~ \in
\left(\phi_{[-\beta,\beta]}
        (\phi_{\tau} \circ j (g_0 ,c)) \right)^{\delta}
\subset V_{\phi_\tau \circ j(g_0,c)}
\]
for all $c \in \overline{X}$ and all
$g \in G$ with $d_G(g,g_0) < \alpha$.
This finishes the proof of
Theorem~\ref{thm:equivariant-cover-for-hyperbolic-group}.
\end{proof}


\subsection*{Acknowledgements}

We thank Igor Mineyev for useful conversations on
his paper and helpful comments on a preliminary version
of Sections~\ref{sec:Mineyevs_flow_space} to
\ref{sec:Properness_of_the_induced_action}.
We thank Tom Farrell who a long time ago
explained the proof of
\cite[Proposition~7.2]{Farrell-Jones-dynamics-I}
to us. Moreover we are grateful to the referee who read through
the manuscript very carefully and made a lot of helpful comments.


\section{Boxes}
\label{sec:boxes}

\begin{convention} \label{conv:properties-of-X-and-Phi-for-boxes}
Throughout this section we consider
\begin{itemize}
\item a discrete group $G$;
\item a metrizable topological space $X$;
\item a proper cocompact $G$-action on $X$;
\item a $G$-equivariant flow $\Phi \colon X \times \IR \to X$
      such that $X - X^\IR$ is locally connected.
\end{itemize}
\end{convention}


\subsection{Basics about boxes}
\label{subsec:Basics_about_boxes}

In this subsection we introduce and study the notion of a box.

\begin{definition}
\label{def:calf-subset-hyp}
Let $B$ be a subset of a $G$-space. Define a subgroup of $G$ by
$$G_B ~ := ~ \{g \in G \mid gB = B\},$$
where $gB := \{g b \mid b \in B\}$.

A subset $B$ of a topological $G$-space is
called an \emph{$\calf$-subset}  for a collection $\calf$
of subgroups of $G$, if
$G_B$ belongs to $\calf$ and for all $g \in G$
we have the implication
$g B \cap B \not= \emptyset ~\Rightarrow B = gB$.
\end{definition}

Notice that $gB = B$ does \emph{not} imply that
$g  b = b$ holds for all $b \in B$.
We denote by $\Fin$ the collection of finite subgroups.

\begin{definition} \label{def:box_wolfgang}
A \emph{box} $B$ is a subset $B \subseteq X$ with the following properties:
\begin{enumerate}
\item $B$ is a compact $\Fin$-subset;

\item There exists a real number $l = l_{B} > 0$,
called the \emph{length} of the box $B$, with the property that for every $x \in B$
there exists real numbers $a_-(x) \le 0 \le a_+(x)$ and $\epsilon(x) > 0$ satisfying
\begin{eqnarray*}
l & = & a_+(x) - a_-(x);
\\
\Phi_{\tau}(x) & \in & B ~ \text{ for } \tau \in [a_-(x),a_+(x)];
\\
\Phi_{\tau}(x) & \not\in & B ~ \text{  for }\tau \in (a_-(x) - \epsilon(x), a_-(x)) \cup
(a_+(x),a_+(x) + \epsilon(x)).
\end{eqnarray*}
\end{enumerate}
\end{definition}

\begin{definition}
\label{def:datas_associated to a box}
Let $B \subseteq X$ be a box.
Then the following data are associated to it:
\begin{itemize}
\item The \emph{length} $l_B > 0$;
\item Let $G_B \subseteq G$ be the finite subgroup
      $\{g \in G \mid gB = B\}$;
\item We denote by $B^{\circ}$  the (topological) interior and
      by $\partial B$ the (topological) boundary of
      $B \subseteq X$;
\item Let $S_B \subseteq B$  be the set of points
      $\{x \in B \mid a_-(x) + a_+(x) = 0\}$.
      We call $S_B$ the \emph{central slice} of $B$;
\item Let $\partial_{\pm}B$ be the set of points
      $\{ x \in B \mid a_{\pm} (x) = 0\} =
       \{ \phi_{a_\pm(x)} (x) \mid x \in S_B \}$.
      We call $\partial_-B$ the \emph{bottom} and
      $\partial_+B$ the \emph{top} of $B$.
      Define the \emph{open bottom} and \emph{open top}
      $\partial_{\pm} B^{\circ} :=
      \{ \phi_{a_\pm(x)}(x) \mid x \in S_B \cap B^\circ  \}$;
\item Let $\pi_B \colon B \to S_B$ be the retraction onto the central
      slice which sends $x$ to
      $\Phi_{\frac{a_+(x) + a_-(x)}{2}}(x)$.
\end{itemize}
\end{definition}

\begin{remark}\label{rem:data_of_a_box}
A box  does not intersect $X^{\IR}$ but may intersect a closed
orbit.
A box  does never contain a closed orbit.
It may happen that a non-closed orbit
meets the central slice
infinitely many times, but whenever it meets the
central slice it has
to leave the box before it comes back to the central slice.
We do not require that
the central slice is connected.
We have for $x \in B$, $\tau \in [a_-(x),a_+(x)]$ that $\Phi_{\tau}(x) \in B$ and
\begin{eqnarray*}
a_-(\Phi_{\tau}(x)) & = & a_-(x) - \tau;
\\
a_+(\Phi_{\tau}(x)) & = & a_+(x) - \tau.
\end{eqnarray*}
\end{remark}

\begin{lemma} \label{lem:basics_about_boxes}
Let $B \subseteq X$ be a box of length $l = l_B$. Then
\begin{enumerate}
\item \label{lem:basics_about_boxes:Fin}
      We get for $g \in G_B$ and $x \in X$
      \begin{eqnarray*}
      a_-(gx) & = & a_-(x);
      \\
      a_+(gx) & = & a_+(x);
      \end{eqnarray*}
      The bottom $\partial_-B$, the open bottom
      $\partial_-B^{\circ}$,
      the top $\partial_+B$, the open top $\partial_+D^{\circ}$,
      the central slice $S_B$ and the interior $B^{\circ}$ are
      $\Fin$-subsets of $G$ and satisfy
      (unless they are empty)
      $$G_B = G_{B^{\circ}} = G_{\partial_-B} =
              G_{\partial_-B^{\circ}} = G_{\partial_+B} =
              G_{\partial_+B^{\circ}} = G_{S_B};$$
\item \label{lem:basics_about_boxes:a_pm_continuous}
      The maps
      $$a_{\pm} \colon B \to \IR, \quad x \mapsto a_{\pm}(x)$$
      are continuous;
\item \label{lem:basics_about_boxes:shape_of_B_in_terms_of_central_slice}
      The  maps
      $$\mu \colon S_B \times [-l/2,l/2] \xrightarrow{\cong} B, \quad (x,\tau) ~ \mapsto \Phi_{\tau}(x)$$
      and
      \begin{eqnarray*}
      && \mu^{-1} \colon B \xrightarrow{\cong} S_B \times [-l/2,l/2], \quad
      x \mapsto \left(\Phi_{\frac{a_-(x) + a_+(x)}{2}}(x), l/2 - a_+(x)\right)
      \end{eqnarray*}
      are to one another inverse $G_B$-homeomorphisms, where $G_B = G_{S_B}$ acts on
      $S_B \times [-l/2,l/2]$ by $g \cdot (s,t) = (gs,t)$.

      We have
      $$\begin{array}{lcl}
      B^{\circ}
      & = &
      \mu\left(\left(S_B \cap B^{\circ}\right) \times (-l/2,l/2)\right);
      \\
      \partial B
      & = &
      \mu\left(\left(S_B \cap \partial B\right) \times [-l/2,l/2] \cup
      \left(S_B \times \{-l/2,l/2\}\right)\right);
      \\
      \partial_{\pm} B
      & = &
      \mu\left(S_B \times \{\pm l/2\}\right);
      \\
      \partial_{\pm} B^{\circ}
      & = &
      \mu\left((S_B \cap B^{\circ}) \times \{\pm l/2\}\right).
      \end{array}
      $$
\item \label{lem:basics_about_boxes:S_BcapB^circ_locally_connected}
      The space $S_B \cap B^{\circ}$ is locally connected;
\item \label{lem:basics_about_boxes:uniform_epsilon(x)}
      There exists $\epsilon_B > 0$ depending only on $B$ such that the
      numbers $\epsilon(x)$ appearing in
      Definition~\ref{def:box_wolfgang} can be chosen so that
      $\epsilon(x) \ge \epsilon_B$ holds for all $x \in B$.
\end{enumerate}
\end{lemma}

\begin{proof}
\ref{lem:basics_about_boxes:Fin}
For $x \in B$ and $g \in G_B$ we have
$\Phi_{\tau}(x) \in B \Leftrightarrow g \cdot \Phi_{\tau}(x) = \Phi_{\tau}(gx) \in B$.
This implies $a_{\pm}(gx) = a_{\pm}(x)$ for $x \in B$ and $g \in G_B$.
We conclude from the definition of the bottom
$\partial_-B$, the open bottom $\partial_-B^{\circ}$, the top $\partial_+B$,
the open top $\partial_+B^{\circ}$, the central slice $S_B$ and the interior $B^{\circ}$
that these sets are $G_B$ invariant and contained in the $\Fin$-subset $B$.
Hence they are themselves $\Fin$-subsets of $X$ and satisfy
$G_B = G_{B^{\circ}} = G_{\partial_-B} = G_{\partial_+B} = G_{S_B}$
if non-empty.
\\[1mm]
\ref{lem:basics_about_boxes:a_pm_continuous}
Consider $x \in B$ and $\epsilon > 0$ with $\epsilon < \epsilon(x)$,
where $\epsilon(x)$ is the number
appearing in Definition~\ref{def:box_wolfgang}.
The points $\Phi_{a_{\pm}(x) \pm \epsilon}(x)$ lie outside
$B$. Since $B$ is compact and $X$ is a Hausdorff space, we can find an open neighborhood
$V_{\pm}$ of $\Phi_{a_{\pm}(x) \pm \epsilon}(x)$ such that $V_{\pm}$
does not meet $B$. Put $U_{\pm} = (\Phi_{a_{\pm}(x) \pm \epsilon})^{-1}(V_{\pm})$.
Then $x \in U_{\pm}$ and $\Phi_{a_{\pm}(x) \pm \epsilon}(u)$ does not lie in $B$ for $u \in U_{\pm}$.
This implies $a_+(u) < a_+(x) + \epsilon$ for $u \in U_+ \cap B $ and $a_-(x) - \epsilon < a_-(u)$ for
$u \in U_- \cap B$. Put $U = U_- \cap U_+\cap B$. Then $U \subseteq B$ is an open neighborhood
of $x$ in $B$ such that $a_-(x) - \epsilon < a_-(u)$ and $a_+(u) < a_+(x) + \epsilon$ holds for $u \in U$.
Since $a_+(u) - a_-(u) = l$ for all $u \in U$, we conclude
$a_{\pm}(u) \in (a_{\pm}(x) - \epsilon,a_{\pm}(x) + \epsilon)$ for all $u \in U$. Hence
$a_{\pm}$ is continuous.
\\[1mm]
\ref{lem:basics_about_boxes:shape_of_B_in_terms_of_central_slice}
The maps $\mu$ and $\mu^{-1}$ are continuous since
$\Phi$ and by assertion~\ref{lem:basics_about_boxes:a_pm_continuous} the maps
$a_+$ and $a_-$ are continuous. One easily checks
that they are inverse to one another.

Since the flow is compatible with the $G$-action and $G_B = G_{S_B}$, the map $\mu$ is
$G_{S_B} = G_B$-equivariant.

Next we prove
\begin{eqnarray}
\mu\left(\left(S_B \cap B^{\circ}\right) \times (-l/2,l/2)\right) & \subseteq  B^{\circ};
\label{lem:basics_about_boxes:shape_of_B_in_terms_of_central_slice(1)}
\\
\mu\left( \left(S_B \cap \partial B\right) \times [-l/2,l/2] \cup
\left(S_B \times \{-l/2,l/2\}\right)\right) & \subseteq  \partial B.
\label{lem:basics_about_boxes:shape_of_B_in_terms_of_central_slice(2)}
\end{eqnarray}
Consider $(x,\tau)  \in (S_B \cap B^{\circ}) \times (-l/2,l/2)$.
Since $a_-$ and $a_+$ are continuous by
assertion~\ref{lem:basics_about_boxes:a_pm_continuous}
and $a_-(x) = -l/2$ and $a_+(x) = l/2$, we can find an open neighborhood
$U \subseteq B^{\circ}$ of $x$ such
that $\tau \in (a_-(u), a_+(u))$
holds for all $u \in U$. Hence $\Phi_{\tau}(U)$
is contained in $B$. Since $\Phi_{\tau}(U)$ is an open subset of
$X$, we have $\Phi_{\tau}(U) \subseteq B^{\circ}$.
Since $\mu(x,\tau) = \Phi_{\tau}(x)$ lies in $\Phi_{\tau}(U)$,
the inclusion~\eqref{lem:basics_about_boxes:shape_of_B_in_terms_of_central_slice(1)} is proven.

Consider $x \in S_B$.
Let $U \subseteq X$ be an open neighborhood of $\Phi_{l/2}(x)$.
Since $\IR \to X,~ \tau \mapsto \Phi_{\tau}(x)$ is a continuous map,
there is an $\epsilon$ with
$0 < \epsilon < l/2$ such that $\Phi_{\tau}(x) \in U$ holds for
$\tau \in (l/2 - \epsilon,l/2 + \epsilon)$.
Since $\{\Phi_{\tau}(x) \mid \tau \in (l/2 - \epsilon,l/2)\}$ is contained in $B$
and $\{\Phi_{\tau}(x) \mid \tau \in (l/2,l/2 + \epsilon)\}$ is contained in $X-B$,
the open neighborhood $U$ of $\Phi_{l/2}(x)$ intersects both $B$ and $X-B$.
This shows $\Phi_{l/2}(x) = \mu(x,l/2) \in \partial B$.
Analogously one proves $\Phi_{-l/2}(x) = \mu(x,-l/2) \in \partial B$.

Consider $x \in S_B \cap \partial B$ and
$\tau \in (-l/2,l/2)$. We want to show $\mu(x,\tau) \in \partial B$.
Suppose the converse. Since $\mu(x,\tau) = \Phi_{\tau}(x)$ belongs to $B$, there must be an
open neighborhood $U$ of $\Phi_{\tau}(x)$ such that $U \subseteq B$.
Since the functions $a_-$ and $a_+$ are continuous
by assertion~\ref{lem:basics_about_boxes:a_pm_continuous} and $a_-(x) = -l/2 < -\tau < a_+(x) = l/2$,
we can arrange by making $U$ smaller that
$-\tau \in (a_-(u),a_+(u))$ holds for all $u \in U$. Hence $\Phi_{-\tau}(U)$ is an open subset of $X$
which is contained in $B$ and contains $x$. This contradicts $x \in \partial B$.
This finishes
the proof of \eqref{lem:basics_about_boxes:shape_of_B_in_terms_of_central_slice(2)}.
Now assertion~\ref{lem:basics_about_boxes:shape_of_B_in_terms_of_central_slice} follows
from~\eqref{lem:basics_about_boxes:shape_of_B_in_terms_of_central_slice(1)} and
\eqref{lem:basics_about_boxes:shape_of_B_in_terms_of_central_slice(2)}.
\\
\ref{lem:basics_about_boxes:S_BcapB^circ_locally_connected}
Since $B^{\circ}$ is an open subset of the locally connected space
$X-X^{\IR}$, it is itself locally connected. Because of
assertion~\ref{lem:basics_about_boxes:shape_of_B_in_terms_of_central_slice}
the space $S_B \cap B^{\circ} \times (-l/2,l/2)$ is locally connected.
Since the projection
$S_B \cap B^{\circ} \times  (-l/2,l/2) \to S_B \cap B^{\circ}$
is an open continuous map and the image of a connected set
under a continuous maps is again connected, $S_B \cap B^{\circ}$ is
locally connected.
\\[1mm]
\ref{lem:basics_about_boxes:uniform_epsilon(x)}
Suppose that such $\epsilon_B$ does not exists.
Then we can find a sequence
$(x_n)_{n \ge 0}$ of elements in $B$ and a sequence
$(\tau_n)_{n \ge 0} > 0$ of positive real numbers
with $\lim_{n \to \infty} \tau_n = 0$ such that
\begin{itemize}
\item[(a)] $\Phi_{a_-(x_n) - \tau_n}(x_n) \in B$ and
           $\Phi_{a_-(x_n) - \tau}(x_n) \not\in B$ for
           $\tau \in (0,\tau_n)$
\end{itemize}
or
\begin{itemize}
\item[(b)] $\Phi_{a_+(x_n) + \tau_n}(x_n) \in B$ and
           $\Phi_{a_+(x_n) + \tau}(x_n) \not\in B$
           for $\tau \in (0,\tau_n)$
\end{itemize}
holds for $n \ge 0$.
By passing to a subsequence we can arrange that
$x_n$ converges to some point $x \in B$ and (a) holds for all
$n \ge 0$ or (b) holds for all $n \ge 0$. We only treat the case (a),
where $\Phi_{a_-(x_n) - \tau_n}(x_n) \in B$ and
$\Phi_{a_-(x_n) - \tau}(x_n) \not\in B$ for $\tau \in (0,\tau_n)$
holds for all $n \ge 0$, the proof in the other case
(b) is analogous.
Put $y_n = \Phi_{a_-(x_n) - \tau_n}(x_n)$. Then $y_n\in \partial_+B$ for all $n \ge 0$
since $y_n \in B$ and $\Phi_{\tau}(y_n) = \Phi_{a_-(x_n) - \tau_n + \tau}(x_n) \not\in B$ holds
for $\tau \in (0,\tau_n)$.
We conclude $\lim_{n \to \infty} a_-(x_n) = a_-(x)$ from assertion~\ref{lem:basics_about_boxes:a_pm_continuous}.
Hence $\lim_{n \to \infty} y_n = \Phi_{a_-(x)}(x)$. Since $y_n \in \partial_+B$ for $n \ge 0$,
we have $\lim_{n \to \infty} y_n \in \partial_+ B$. This contradicts $\Phi_{a_-(x)}(x) \in \partial_-B$
since $l_B > 0$.
\end{proof}

We mention that in general $S_B$ itself is not
locally connected.

\begin{remark} \label{rem:compactness_and_continuity}
It is a little bit surprising that the function $a_{\pm}$ is
continuous as stated in
Lemma~\ref{lem:basics_about_boxes}~\ref{lem:basics_about_boxes:a_pm_continuous}
since there seem to be no link between different flow lines
entering the box. The point here is that we require the box to be
compact. If $G$ is trivial and we consider the flow $\phi_{\tau}(x,y)
= (x + \tau,y)$ on $\IR^2$, the subset of $\IR^2$ given by
$$B ~ := ~ \{(x,y) \mid x,y \in [-1,1], y \not= 0\} \cup \{(x,0) \mid x \in [0,2]\}$$
satisfies all the requirements of a box of length $2$ except for
compactness and the functions $a_{\pm}$ are not continuous at
$(0,0)$.
\end{remark}

\begin{definition} \label{def:restriction of boxes}
Consider a  box $B$ of length $l_B$.
Let $V \subseteq S_B$ be a non-empty closed $\Fin$-subset
of the $G_{S_B}$-space $S_B$
and $a,b$ real numbers with
$-l_B/2 \le a < b \le l_B/2$. Define a new box of length $b-a$ by
$$B(V;a,b) ~ := ~ \Phi_{[a,b]}(V):= \{\Phi_{\tau}(v) \mid v \in V, \tau \in [a,b]\}.$$
If $a = -v/2$ and $b = v/2$ for some $v \in [0,w]$ we abbreviate
$$B(V;v) := B(V;-v/2,v/2).$$
If $a = -l_B/2$ and $b = l_B/2$, we abbreviate
$$B(V) := B(V;-l_B/2,l_B/2),$$
and call $B(V)$ the \emph{restriction} of $B$ to $V$.
\end{definition}

We have to show that $B(V;a,b)$ is again a box. Since $V$ is
a closed subset of the compact set $S_A$,
it is compact. Hence $V \times [a,b]$ is compact.
We conclude that $B(V,a,b)$ as the image of a compact set under the
continuous map $\Phi \colon X \times \IR \to X$ is compact.
{}From Lemma~\ref{lem:basics_about_boxes} we get $G_B = G_{S_B}$
and the $G_B$-equivariant homeomorphism
\[
\mu_B \colon S_B \times [-l_B/2,l_B/2] \xrightarrow{\cong} B,
\quad (x,\tau) ~ \mapsto \Phi_{\tau}(x).
\]
The subset $V$ of the $G_B$-space $S_B$ is a $\Fin$-subset.
Hence
$B(V;a,b) = \mu_B(V \times [a,b])$ is a $\Fin$-subset of the
$G_B$-space $B$. Since $B$ is a $\Fin$-subset of the $G$-space $X$,
$B(V;a,b)$ is a $\Fin$-subset of the $G$-space $X$. Consider
$x \in B(V;a,b)$. We can write it as $x = \Phi_{\tau}(v)$ for $v \in S_B$ and $\tau \in [a,b]$.
Put $a_-(x) = a - \tau$ and $a_+(x) = b- \tau$.
Let $\epsilon(x)$ for $x \in B$ be the number
appearing in the Definition~\ref{def:box_wolfgang} of a box for $B$.
Now one easily checks that the collections $a_{\pm}(x)$ and $\epsilon(x)$
have the desired properties appearing in Definition~\ref{def:box_wolfgang} for $B(V;a,b)$.
This shows that $B(V;a,b)$ is a box. We have
\begin{eqnarray*}
\partial_-B(V;a,b) & = & \Phi_{a}(V);
\\
\partial_+B(V;a,b) & = & \Phi_{b}(V);
\\
\partial_-B(V;a,b)^{\circ} & = & \Phi_{a}(V \cap B^{\circ});
\\
\partial_+B(V;a,b) & = & \Phi_{b}(V \cap B^{\circ});
\\
S_{B(V;a,b)} & = & \Phi_{\frac{a+b}{2}}(V).
\end{eqnarray*}
In particular $B(V;v)$ is a box of length $v$ with central slice $V$ and
$B(V)$ is  box of length $l_{B}$ with central slice $V$.


\subsection{Constructing boxes}
\label{subsec:Constructing_boxes}

\begin{lemma}
\label{lem:existence-of-boxes}
For every $x \in X - X^\IR$ there exists a
non-equivariant box whose interior contains $x$.
\end{lemma}

The following proof is a variation of an argument
used in \cite[Theorem~1.2.7]{Palais-slice-non-compact-groups}.

\begin{proof}[Proof of Lemma~\ref{lem:existence-of-boxes}]
Because the $G$-action on the metrizable space $X$ is
proper and cocompact $X$ is locally compact, compare
\cite[Theorem 1.38 on p.27]{Lueck(1989)}.
Let $\alpha > 0$ with $\Phi_\alpha (x) \neq x$.
Let $W_\pm$ be a closed neighborhood of $\Phi_{\pm \alpha}(x)$
that does not contain $x$.
By continuity there exist compact neighborhoods $U' \subset U$
of $x$ and $\e > 0$ such that
$\Phi_{[\pm \alpha - \e,\pm \alpha + \e]} U \subset W_{\pm}$ and
$\Phi_{[ - \e, \e]} U' \subset U$ and $U$ is disjoint from $W_- \cup W_+$.
Let $f \colon X \to [0,\infty)$ be a continuous function with
$f(y) = 1$ for $y \in U$ and $f(y) = 0$ for $y \in W_- \cup W_+$.
Define $\psi \colon U \to \IR$ by
\[
\psi(y) = \ln \left(
     \int_{-\alpha}^\alpha f(\Phi_\tau(y)) e^{-\tau} d\tau
               \right).
\]
(The logarithm makes sense because the integrant is
non-negative, continuous, and positive for $\tau = 0$.)
Let $y \in U'$ .
If $s \in [\pm \alpha - \e,\pm \alpha + \e]$
then $\Phi_{s}(y) \in W_\pm$ and therefore $f(\Phi_{s}(y)) = 0$.
Using this we compute for $\delta \in [-\e,\e]$
\begin{eqnarray*}
\psi( \Phi_\delta (y) ) & = &
\ln \left(
\int_{-\alpha}^{\alpha} f(\Phi_{\tau+\delta}(y)) e^{-\tau} d\tau
\right) \\
& = &
\ln \left(
\int_{-\alpha+\delta}^{\alpha+\delta} f(\Phi_\tau(y)) e^{-\tau + \delta} d\tau
\right) \\
& = &
\ln \left( e^\delta
\int_{-\alpha+\delta}^{\alpha+\delta} f(\Phi_\tau(y)) e^{-\tau} d\tau
\right) \\
& = &
\delta + \ln \left(
\int_{-\alpha}^{\alpha} f(\Phi_\tau(y)) e^{-\tau} d\tau
\right) \\
& = &
\delta + \psi(y).
\end{eqnarray*}
Define $S := U' \cap \psi^{-1}(\psi(x))$.
Then
$B := \{ \Phi_\tau (s) \; | \; s \in S, \tau \in [-\e/2, \e/2]\}$
is a box whose interior contains $x$.
\end{proof}

\begin{definition} \label{def:deviation_of_subsets}
Let $C$ be a box of length $l_C$. Let
$$\mu_C \colon S_C \times [-l_C/2,l_C/2] \to C,
  \quad (x,\tau) \mapsto \Phi_{\tau}(x)$$
be the homeomorphism appearing in
Lemma~\ref{lem:basics_about_boxes}
~\ref{lem:basics_about_boxes:shape_of_B_in_terms_of_central_slice}.

Consider a subset $S \subseteq C$.
It is called \emph{transversal to the flow with respect to $C$}
if $\mu_C^{-1}(S) \cap \{x\} \times [-l_C/2,l_C/2]$
consists of at most one point for every $x \in S_C$.
\end{definition}

\begin{lemma}\label{lem:achieving_tranversality_to_the_flow}
Let $C$ be box of length $l_C$.
Let $B$ be a box with $B \subseteq C$.
Then we can find for every $x \in S_B$ a closed
neighborhood $U \subseteq S_B$ of $x$ satisfying
\begin{enumerate}
\item $U$ is a $G_x$-invariant $\Fin$-subset of the
      $G_B$-space $S_B$;
\item $U$ is transversal to the flow with respect to $C$.
\end{enumerate}
\end{lemma}

\begin{proof}
Let $\tau_{S_B} \colon S_B \to [-l_C/2,l_C/2]$ be the continuous function
given by the restriction to $S_B$ of the
composite of the projection
$S_C \times  [-l_C/2,l_C/2]  \to  [-l_C/2,l_C/2]$ and $\mu_C^{-1}$.
Let $U_1 \subseteq S_B$ be a closed neighborhood of $x \in S_B$
such that $|\tau_{S_B}(u) - \tau_{S_B}(x) | < l_B / 2$ holds for $u \in U_1$.
Choose a closed neighborhood $U_2 \subseteq S_B$ of $x$
such that $gU_2 \cap U_2 \not= \emptyset \Rightarrow g \in G_x$ holds
for $g \in G_B$. Put $U = \bigcap_{g \in G_x} g(U_1 \cap U_2)$.
This is a closed neighborhood of $x$ in $S_B$
which is $G_x$-invariant, a $\Fin$-subset of the $G_B$-space $S_B$
and satisfies  $|\tau_{S_B}(u) - \tau_{S_B}(x)| \le l_B / 2$ for $u \in U$.

It remains to show  that $U$ is transversal to the flow
with respect to $C$. Suppose the converse. So we can find
$u_0,u_1 \in U$, $v \in S_C$ and $\tau_0, \tau_1$ such that
$u_0 = \Phi_{\tau_0}(v)$, $u_1 = \Phi_{\tau_1}(v)$ and $\tau_0 \not= \tau_1$.
Note that $\tau_0 = \tau_{S_B}(u_0)$ and $\tau_1 = \tau_{S_B}(u_1)$.
Since $B$ is a box of length $l_B$,
we can find for $i = 0,1$ real numbers $\epsilon_i > 0$ such
that $\Phi_{[-l_B/2,l_B/2]}(u_i) \subseteq B$
and $\Phi_{\tau}(u_i) \not\in B$ for
$\tau \in (-l_B/2 - \epsilon_i,-l_B/2) \cup (l_B/2,l_B/2 + \epsilon_i)$.
This implies $|\tau_1 - \tau_0| > l_B$.
But by the definition of $U$,
$| \tau_1 - \tau_0 | = | \tau_{S_B}(u_1) - \tau_{S_B}(u_0) |
  \leq | \tau_{S_B}(u_1) - \tau_{S_B} (x) | +
         | \tau_{S_B} (x) - \tau_{S_B}(u_0) | \leq l_B$.
This is the required contradiction.
\end{proof}

Next we show for $x \in X$ that the existence of non-equivariant box containing $x$ in its interior
already implies the existence of an equivariant box containing $x$ in its interior.
The basic idea of proof is an averaging process in the time direction of the flow applied to
the central slice of a non-equivariant box.

\begin{lemma} \label{lem:from_boxes_to equivariant boxes}
Suppose for the point $x \in X$ that there is a non-equivariant box
whose interior contains $x$.

Then there exists a box $B$ in the sense
of Definition~\ref{def:box_wolfgang} satisfying
\begin{enumerate}
\item $G_B = G_{S_B} = G_x$;
\item $x \in S_B \cap B^{\circ}$;
\item $S_B$ is connected.
\end{enumerate}
\end{lemma}

\begin{proof}
Let $C$ be a non-equivariant box, i.e. a box in the sense
of  Definition~\ref{def:box_wolfgang} in the case, where the group $G$
is trivial, such that $x \in C^{\circ}$. Let $l = l_C$ be the length of $C$.
Since the $G$-action on $X$ is proper by assumption and $X$ is metrizable,
we can find a closed neighborhood $U$ of $x$ such that $U$ is a
$\Fin$-subset of $X$ with $G_U = G_x$.
We can assume without loss of generality that $x \in S_C \cap C^{\circ}$
and $C \subseteq U$ holds and for every
$\tau \in [-l,-l/2] \cup [l/2,l]$ and
$s \in S_C$ we have $\Phi_{\tau}(s) \notin C^{\circ}$,
otherwise replace $C$ by an appropriate restriction.
Let $S_C$ be the central slice of $C$. Let
$$\mu \colon S_C \times [-l/2,l/2]\to C,
      \quad  (s,\tau) \mapsto \Phi_{\tau}(s)$$
be the homeomorphism of
Lemma~\ref{lem:basics_about_boxes}~\ref{lem:basics_about_boxes:shape_of_B_in_terms_of_central_slice}.
Since $S_C$ is compact, $C^\circ \subseteq X$
is open, $G_x$ is finite and $X$ is metrizable, we can find a compact  neighborhood
$S_0 \subseteq S_C \cap C^{\circ}$ of $x \in S_C \cap C^{\circ}$ such
that $g S_0 \subseteq C^{\circ}$
holds for all $g \in G_x$. Define
$$S_1 ~ = ~ \bigcap_{g \in G_x}
        S_0 \cap \pi_{C}(\mu^{-1}(gS_0)),$$
where $\pi_C \colon S_C \times [-l/2,l/2] \to S_C$ is the projection.
Then $S_1 \subseteq S_0$ is a compact neighborhood of $x$ in $S_0$.
By construction there exists for every $g \in G_x$ and $s \in S_1$ precisely one
element $\tau_g(s) \in (-l/2,l/2)$ such that $\Phi_{\tau_g(s)}(s) \in gS_0$.
The function $\tau_g(s)$ is continuous in $s$ and has image in $(-l/2 + \delta,l/2 - \delta)$ for some
$\delta$ with $0 < \delta <l/2$, since it is the restriction to $S_1$ of the continuous function with a compact source
$$
\tau_g \colon S_0  \to (-l/2,l/2),
\quad s \mapsto - \pi_C
\circ \mu^{-1}(g^{-1}s).
$$
Define the continuous function
$$\tau \colon S_1 \to (-l/2, l/2),  \quad s ~ \mapsto ~ \frac{1}{|G_x|} \cdot
\sum_{g \in G_x} \tau_g(s).$$
Put
$$S_2 = \left\{\Phi_{\tau(s)}(s) \mid s \in S_1\right\}.$$
Next we show that $S_2 \subseteq C$ is $G_x$-invariant. Consider $g'
\in G_x$ and $u \in  S_2$. Write $u = \Phi_{\tau(s)}(s)$ for appropriate $s\in  S_1$.
Let $s' \in S_0$ be the element uniquely determined by
$\Phi_{\tau_{(g')^{-1}}(s)}(s) = (g')^{-1}s'$.
Then we get for $g \in G_x$
$$\Phi_{\tau_g(s) - \tau_{(g')^{-1}}(s)}(s') ~ = ~ \Phi_{\tau_g(s)} \circ \Phi_{-\tau_{(g')^{-1}}(s)}(s')
~ = ~  \Phi_{\tau_g(s)}(g's) = g'\Phi_{\tau_g(s)}(s) \in g'gS_0.$$
Since $g'gS_0 \subseteq C^{\circ}$ and $\tau_g(s)$ and
$\tau_{(g')^{-1}}(s)$ belong to $(-l/2,l/2)$
and hence $\tau_g(s) - \tau_{(g')^{-1}}(s) \in (-l,l)$
we conclude $\tau_g(s) - \tau_{(g')^{-1}}(s) \in (-l/2,l/2)$.
Hence $s' \in S_1$ and
$\tau_{g'g}(s') = \tau_g(s) - \tau_{(g')^{-1}}(s)$.
Since this
implies $\tau(s') = \tau(s) - \tau_{(g')^{-1}}(s)$, we conclude
\begin{multline*}
g' \cdot u ~ = ~ g' \cdot \Phi_{\tau(s)}(s) ~ = ~ \Phi_{\tau(s)}(g's) ~ = ~ \Phi_{\tau(s') + \tau_{(g')^{-1}}(s)}(g's)
\\
~ = ~  \Phi_{\tau(s')}\left(g' \cdot \Phi_{\tau_{(g')^{-1}}(s)}(s)\right)
~ = ~  \Phi_{\tau(s')}\left(g' \cdot (g')^{-1} s'\right)
~ = ~ \Phi_{\tau(s')}s' ~ \in S_2.
\end{multline*}
Since $S_1 \subseteq S_C \cap C^{\circ}$ is compact,
$S_2$ is a compact $G_x$-invariant subset of
$C^{\circ}$ with $x \in S_2$.
Let $S_3$ be the component of $S_2$ which contains $x$.
Then $S_3$ is a connected closed subset of $S_2$.
Since $gS_3  \cap S_3$ contains $x$ for $g \in G_x$,
the subset $S_3$ is $G_x$-invariant.
Thus $S_3 \subseteq C^{\circ}$ is a compact
connected $G_x$-invariant subspace containing $x$.

We can find $\delta$ with $0 <\delta < l/2$ such that
$B \subseteq C^{\circ}$ holds for
$B  := \Phi_{[-\delta/2,\delta/2]}(S_3)$.

Next we show that $B$ is a box of length $\delta$.
Since $S_3$ is $G_x$-invariant and the flow $\Phi$ commutes with the $G$-action
the subset $B \subseteq X$ is $G_x$-invariant.
Recall that $B \subseteq C$ holds and $C$ is a  $\Fin$-subset
of $X$ with $G_C = G_x$. Hence  $B$ is a compact $\Fin$-subset of $X$.
Consider $y \in B$.
There is precisely one element $s \in S_3$ and
$\tau \in [-\delta/2,\delta/2]$ satisfying
$y = \Phi_{\tau}(s)$ since $S_C$ and hence $S_3$ is transversal to
the flow with respect to $C$.
Put $a_-(y) = - \delta/2 - \tau$, $a_+(y) = \delta/2 - \tau$,
$\epsilon_-(y) = \epsilon_+(y) = l/2$. Then
\begin{eqnarray*}
\delta & = & a_+(y) - a_-(y);
\\
\Phi_{\tau}(y) & \in & B ~ \text{ for } \tau \in [a_-(y),a_+(y)];
\\
\Phi_{\tau}(y) & \not\in & B ~ \text{  for }
   \tau \in (a_-(y) - \epsilon(y), a_-(y)) \cup (a_+(y),a_+(y) + \epsilon(y)).
\end{eqnarray*}
Hence $B$ is a box with connected central slice $S_B = S_3$.
We have $x \in S_3$.
The projection $\pi_{C} \colon C \to
S_{C}$ induces a homeomorphism $S_2 \to S_1$ and maps the
component $S_3$ of $S_2$ to a component $S_1'$ of  $S_1$. Since $S_1$ is an
open neighborhood of $x$ in the space $S_{C} \cap C^{\circ}$ which is locally connected by
Lemma~\ref{lem:basics_about_boxes}~\ref{lem:basics_about_boxes:S_BcapB^circ_locally_connected},
the component $S_1'$ is a
neighborhood of $x$ in the space $S_{C} \cap C^{\circ}$. Since $\tau$ is continuous
we conclude from Lemma~\ref{lem:basics_about_boxes}~\ref{lem:basics_about_boxes:shape_of_B_in_terms_of_central_slice}
that $x$ lies in the interior of $B$.
\end{proof}

\begin{definition}  \label{def:G-period}
For $x \in X$ define its \emph{$G$-period}
$$
\peri{G}{\Phi}(x) ~ = ~ \inf\{ \tau  \mid \tau > 0,
        \; \exists \; g \in G \text{ with }
        \Phi_{\tau}(x) = gx\} \quad \in [0,\infty],
$$
where the infimum over the empty set is defined to be
$\infty$. If $L \subseteq X$ is an orbit of the flow $\Phi$, define its
\emph{$G$-period} by
$$\peri{G}{\Phi}(L) ~ := ~ \peri{G}{\Phi}(x)$$
for any choice of $x \in X$ with $L = \Phi_{\IR}(x)$.

For $r \ge 0$ put
\begin{eqnarray*}
X_{>r} & := & \{x \in X \mid \peri{G}{\Phi}(x) > r\};
\\
X_{\le r} & := & \{x \in X \mid \peri{G}{\Phi}(x) \le r\}.
\end{eqnarray*}
\end{definition}

Consider $x \in X$.
Then the $G$-period $\peri{G}{\Phi}(x)$ vanishes if and
only if $x \in X^{\IR}$.
We have $\peri{G}{\Phi}(x) = \infty$ if and only if the
orbit through $x$ is not periodic
and
$g\Phi_{\IR}(x) \cap \Phi_{\IR}(x) = \emptyset$ holds for all
$g \neq 1$,
or, equivalently, the orbit through $Gx$ in the quotient space $G\backslash X$ with respect to the induced flow
is not periodic.  If $0 < \peri{G}{\Phi}(x) < \infty$, then the
properness of the $G$-action implies the existence of
$g \in G$ such that $\Phi_{\peri{G}{\Phi}(x)}(x) = gx$ and  $\peri{G}{\Phi}(x)$ is the period of the
periodic orbit through $Gx$ in the quotient space $G\backslash X$ with respect to the induced flow.

Next we show for a point $x$ that we can find  an equivariant box around a given
compact part of the flow line, where the compact part is as long
as the $G$-orbit length allows. The idea of proof is to take an
equivariant box which contains $x$ in its interior, making its
central slice very small by restriction and then prolonging the
box  along the flow line though $x$.

\begin{lemma} \label{lem:long_equivariant_boxes_around_x}
Suppose for the point $x \in X$ that there is a
non-equivariant box whose interior contains $x$.
Consider a real number $l$ with $0 < l < \peri{G}{\Phi}(x)$.

Then we can find a box $C$
which satisfies
\begin{itemize}
\item $l_C = l$;
\item $G_C = G_x$;
\item $x \in S_C \cap C^{\circ}$;
\item $S_C$ is connected.
\end{itemize}
\end{lemma}

\begin{proof}
{}From Lemma~\ref{lem:from_boxes_to equivariant boxes}
we conclude the existence of a box $B$ in the sense
of Definition~\ref{def:box_wolfgang}
which satisfies $G_B = G_x$, $S_B$ is connected
and $x \in S_B \cap B^{\circ}$.
Let $l_B$ be the length of $B$.
{}From Lemma~\ref{lem:basics_about_boxes}
~\ref{lem:basics_about_boxes:uniform_epsilon(x)}
we obtain a number $\epsilon_B> 0$ such that for
every $y \in S_B$ and
$\tau \in (-l_B/2 - \epsilon_B,-l_B/2)
           \cup (l_B/2,l_B/2 + \epsilon_B)$
the element $\Phi_{\tau}(y)$ does not belong to $B$.
We can arrange by restricting $B$ and diminishing
$\epsilon_B$ that $l_B< l$  and
$l + \epsilon_B < \peri{G}{\Phi}(x)$ holds.

Next we show
$$
\Phi_{[-l/2,l/2]}(x) \cap g \Phi_{[-l/2,l/2]}(x)
             \not= \emptyset \Rightarrow g \in G_x.
$$
Namely, consider
$y \in \Phi_{[-l/2,l/2]}(x) \cap g \Phi_{[-l/2,l/2]}(x)$.
Then $y = \Phi_{\tau}(x) = g\Phi_{\sigma}(x) $
for appropriate $\tau,\sigma \in [-l/2,l/2]$.
This implies $\Phi_{\tau - \sigma}(x) = gx$ and
$|\tau - \sigma| \le l < \peri{G}{\Phi}(x)$.
We conclude $\tau - \sigma = 0$ and hence $g \in G_x$.

Since the $G$-action on $X$ is proper,
we can find a closed neighborhood
$V_x^1 \subseteq X$ of $x$ such that
$\Phi_{[-l/2,l/2]}(V_x^1) \cap g \Phi_{[-l/2,l/2]}(V_x^1)
                       \not= \emptyset \Rightarrow g \in G_x$
holds.
{}From $l + \epsilon_B < \peri{G}{\Phi}(x)$ we conclude that
\begin{eqnarray*}
\Phi_{[-l/2+l_B,l/2]}(x) \cap \Phi_{[-l/2- \epsilon_B,-l/2]} (x)
           & = &\emptyset;
\\
\Phi_{[-l/2,l/2-l_B]}(x) \cap \Phi_{[l/2,l/2 +\epsilon_B] } (x)
           & = &\emptyset.
\end{eqnarray*}
Since $\Phi$ is continuous and $[-l/2+l_B,l/2]$,
$[-l/2- \epsilon_B,-l/2]$,
$[-l/2,l/2-l_B]$ and $[l/2,l/2 +\epsilon_B]$ are compact,
we can find a closed neighborhood
$V^2_x \subseteq X$ of $x$ such that
\begin{eqnarray*}
\Phi_{[-l/2+l_B,l/2]}(V^2_x) \cap
   \Phi_{[-l/2- \epsilon_B,-l/2]} (V^2_x)  & = &\emptyset;
\\
\Phi_{[-l/2,l/2-l_B]}(V^2_x) \cap
   \Phi_{[l/2,l/2 +\epsilon_B] } (V^2_x)  & = &\emptyset.
\end{eqnarray*}
Put
$$
V_x = \bigcap_{g \in G_x} g  \cdot
              \left(V_x^1 \cap V_x^2\right).
$$
Then $V_x \subseteq X$ is a closed
$G_x$-invariant neighborhood of $x$ with the properties
\begin{itemize}
\item $\Phi_{[-l/2,l/2]}(V_x)$ is a
      $\Fin$-subset of the $G$-space $X$;

\item $G_{\Phi_{[-l/2,l/2]}(V_x)} = G_x$;

\item $\Phi_{[-l/2+l_B,l/2]}(V_x) \cap
        \Phi_{[-l/2- \epsilon_B,-l/2]} (V_x) ~ = ~ \emptyset$;

\item $\Phi_{[-l/2,l/2-l_B]}(V_x) \cap
        \Phi_{[l/2,l/2 +\epsilon_B] } (V_x)  ~ = ~ \emptyset$.

\end{itemize}
Let $V_x^{\circ} \subseteq X-X^{\IR}$ be the interior of $V_x$.
Let $T$ be the component of
$S_B \cap B^{\circ} \cap V_x^{\circ}$
that contains $x$.
Since $S_B \cap B^{\circ}$ is locally connected by
Lemma~\ref{lem:basics_about_boxes}
~\ref{lem:basics_about_boxes:S_BcapB^circ_locally_connected}
and $S_B \cap B^{\circ} \cap V_x^{\circ}$ is an open subset of
$S_B \cap B^{\circ}$, the component $T$ is an open subset
of $S_B \cap B^{\circ} \cap V_x^{\circ}$ and hence of $S_B$.
Let $\overline{T}$ be the closure of $T$ in $S_B$.
This is a closed connected $G_x$-invariant neighborhood of
$x \in S_B$ which is contained in $V_x$.
Since $S_B$ is a $\Fin$-subset of $X$ with $G_{S_B} = G_x$,
$\overline{T}$ is a $\Fin$-subset of $S_B$
and we can consider the restriction $B(\overline{T})$.
We can assume without loss of generality that the central slice
$S_B$ is a $G_x$-invariant connected subset of $V_x$,
otherwise replace $B$ by the restriction $B(\overline{T})$.

We define $C :=\Phi_{[-l/2,l/2]}(S_B)$.

Next we show that $C$ is a box of length $l$.
Since $C \subseteq \Phi_{[-l/2,l/2]}(V_x)$ and
$C$ is $G_x$-invariant, $C$ is a compact $\Fin$-subset
of the $G$-space $X$.
Consider $y \in C$.
We can write it as $y = \Phi_{\tau_y}(s)$ for
$\tau_y \in [-l/2,l/2]$ and $s \in S_B$.
Put $a_-(y) = -l/2 - \tau_y$ and $a_+(y) = l/2 - \tau_y$.
Obviously
$l = a_+(y) - a_-(y)$ and $\Phi_{\tau}(y) \in C$
for $\tau \in [a_-(y),a_+(y)]$.
It remains to show that
$\Phi_{\tau'}(y) \not\in  C$ holds for
$\tau' \in (a_-(y) - \epsilon_B, a_-(y))
          \cup (a_+(y),a_+(y) + \epsilon_B)$.
This is equivalent  to showing that
$\Phi_{\tau'}(s) \not\in  C$ holds for
$\tau' \in (-l/2 - \epsilon_B, -l/2) \cup (l/2,l/2 + \epsilon_B)$.
Since $s \in S_B \subseteq V_x$, we have
\begin{eqnarray*}
\Phi_{[-l/2+l_B,l/2]}(S_B) \cap \Phi_{[-l/2- \epsilon_B,-l/2]} (s)
                            & = &\emptyset;
\\
\Phi_{[-l/2,l/2-l_B]}(S_B) \cap \Phi_{[l/2,l/2 +\epsilon_B] } (s)
                            & = &\emptyset.
\end{eqnarray*}
The main property of $\epsilon_B$ is
$$
\Phi_{[-l_B/2,l_B/2]}(S_B) \cap \Phi_{(-l_B/2- \epsilon_B,-l_B/2)
       \cup (l_B/2,l_B/2 + \epsilon_B)} (s) = \emptyset.
$$
Applying $\Phi_{(l_B - l)/2}$ respectively
$\Phi_{(l - l_B)/2}$ we obtain
\begin{eqnarray*}
\Phi_{[-l/2,-l/2+l_B]}(S_B) \cap \Phi_{(-l/2- \epsilon_B,-l/2)} (s)
                    & = &\emptyset;
\\
\Phi_{[l/2-l_B,l/2]}(S_B) \cap \Phi_{(l/2,l/2 +\epsilon_B) } (s)
                    & = &\emptyset.
\end{eqnarray*}
We conclude
$$
\Phi_{[-l/2,l/2]}(S_B) \cap \Phi_{(-l/2- \epsilon_B,-l/2)
\cup (l/2,l/2 +\epsilon_B) } (s)  ~ = ~ \emptyset.
$$
Hence $C$ is a box of length $l$.
By construction $S_C = S_B$ is connected, $G_C = G_{S_C} = G_x$
and $x \in S_B \cap B^{\circ}$ and hence
$x \in S_C \cap C^{\circ}$.
\end{proof}

\begin{lemma}
\label{lem:A-B-C-cover-for-K}
Consider real numbers $a,b,c > 0$ satisfying $c > a + 2b$.
Let $K$ be a cocompact $G$-invariant
subset of $X_{> a + 2b+2c}$.

Then there exist a $G$-set $\Lambda$
and for every $\lambda \in \Lambda$ boxes
$A_\lambda \subseteq B_\lambda \subseteq C_\lambda$
such that
\begin{enumerate}
\item \label{lem:A-B-C-cover-for-K:cofinite}
      $\Lambda$ is $G$-cofinite;

\item \label{lem:A-B-C-cover-for-K:long} We have
      \begin{eqnarray*}
      l_{A_{\lambda}} & = & a;
      \\
      l_{B_{\lambda}} & = & a + 2b;
      \\
      l_{C_{\lambda}} & = & a + 2b +2c;
      \end{eqnarray*}

\item \label{lem:A-B-C-cover-for-K:S_C is connected}
      $S_{C_{\lambda}}$ is connected;

\item \label{lem:A-B-C-cover-for-K:parallel} We have
      $S_{A_\lambda} \subseteq S_{B_\lambda} \subseteq S_{C_{\lambda}}$;

\item \label{prop:A_in B_inC}
   $A_{\lambda} \subseteq B_{\lambda}^{\circ}$ and $B_{\lambda} \subseteq C_{\lambda}^{\circ}$;

\item \label{lem:A-B-C-cover-for-K:cover}
      $K \subseteq \bigcup_{\lambda \in \Lambda} A_{\lambda}^{\circ}$;

\item \label{lem:A-B-C-cover-for-K:equivariant}
      $g A_\lambda = A_{g \lambda}$,
      $g B_\lambda = B_{g \lambda}$ and
      $g C_\lambda = C_{g \lambda}$
      for $g \in G$;
\item \label{lem:A-B-C-cover-for-K:intersection-control}
      If $B_\lambda \cap B_{\lambda'} \neq \emptyset$, then $B_\lambda \subseteq C_{\lambda'}^{\circ}$
      and $S_{B_{\lambda}}$ is transversal to the flow with respect to $C_{\lambda'}$.
\end{enumerate}
\end{lemma}

\begin{proof}
Lemma~\ref{lem:long_equivariant_boxes_around_x} implies that we can
find for every $x \in X_{> a + 2b+2c}$ a box $C_x$ of length
$a + 2b +2c$ such that $x \in S_{C_x} \cap C_x^{\circ}$
and $G_{C_x} = G_x$ holds and $S_{C_x}$ is connected.
Since $G_x$ is finite, we can find a $G_x$-invariant closed
neighborhood $T_x$ of $x$ in $S_{C_x}$ such that
$T_x \subseteq  S_{C_x} \cap C_x^{\circ}$ holds.
We can arrange that $gC_x = C_{gx}$ and $gT_x = T_{gx}$
holds for $g \in G$ and $x \in X_{>a + 2b  +2c}$.
Obviously $K \subseteq  \bigcup_{x \in K} C_x(T_x;a)^{\circ}$,
where we use the notation from
Definition~\ref{def:restriction of boxes}.
Since $K$ is cocompact and $G$-invariant, we can find a cofinite
$G$-subset $I \subseteq K$  satisfying
\begin{eqnarray}
K & \subseteq & \bigcup_{x \in I} C_x(T_x;a)^{\circ}.
\label{lem:A-B-C-cover-for-K:(1)}
\end{eqnarray}

Fix $x \in I$. Consider $y \in T_x$.
Since the $G$-action is proper, $g \cdot C_z(T_z) = C_{gz}(T_{gz})$
holds for $z \in I$ and $g \in G$ and $C_{x}(\{y\};a + 2b)$
and  $C_z(T_z;a+2b)$ are compact, we can find a closed $G_y$-invariant
neighborhood
$V_y \subseteq T_x$ of $y$ such that for all $z \in I$
$$C_{x}(\{y\};a + 2b) \cap C_z(T_z;a+2b) = \emptyset
~ \Rightarrow ~
C_{x}(V_y;a + 2b) \cap C_z(T_z;a+2b) = \emptyset.$$

For $y \in T_x$ we define
$$I_y = \{z \in I \mid C_{x}(\{y\};a + 2b)
               \cap C_z^{\circ} \not= \emptyset\}.$$
Since the set $G$-set $I$ is cofinite, $C_{x}(\{y\};a + 2b)$ and
$C_z$ are compact  and the $G$-action on $X$ is proper, the
set $I_y$ is finite.
{}From $a + 2b < c$ and $l_{C_z} = a + 2b +2c$
we conclude for $z \in I_y$ that $C_{x}(\{y\};a + 2b) =
\Phi_{[-a/2 - b,a/2 + b]}(y)$ is contained in $C_z^{\circ}$.
Since $C_{x}(\{y\};a + 2b)$ is compact, we can find for $z \in I_y$ a
closed $G_y$-invariant neighborhood $U_y(z) \subseteq T_x$ of $y$ such that
$C_x(U_y(z);a + 2b) = \Phi_{[-a/2 - b,a/2 + b]}(U_y(z))$ is
contained in $C_z^{\circ}$. Because of
Lemma~\ref{lem:achieving_tranversality_to_the_flow}
applied to $C_x(U_y(z);a + 2b)
\subseteq C_z$ we can assume without loss of generality that
$U_y(z)$ is transversal to the flow with respect to $C_z$ for
every $z \in I_y$.

Put $U_y  :=  V_y \cap \bigcap_{z \in I_y} U_y(z)$.
Then $U_y \subseteq T_x$ is a $G_y$-invariant closed
neighborhood of $y$ such that
\begin{eqnarray}
\label{lem:A-B-C-cover-for-K:(2)}
& & C_x(U_y;a+2b) \subseteq C_z^{\circ}
\quad \text{ if } z \in I_y;
\\
\label{lem:A-B-C-cover-for-K:(3)} & &  U_y \text{ is
transversal to the flow with respect to } C_z  \text { for } z \in
I_y;
\\
\label{lem:A-B-C-cover-for-K:(4)}
& & C_x(U_y;a + 2b) \cap C_z(T_z;a+2b) = \emptyset
\\
\nonumber & & \hspace*{30mm}  ~ \text{ if }z \in I \text{ and }
C_x(\{y\};a + 2b) \cap C_z(T_z;a+2b) = \emptyset.
\end{eqnarray}
Choose a $G_y$-invariant closed neighborhood
$W_y \subseteq T_x$ of $y$ such that $W_y \subseteq U_y^{\circ}$.
Obviously $T_x = \bigcup_{y \in T_x} W_y^{\circ}$.
Since $T_x$ is compact, we can find
$y(x)_1, y(x)_2, \ldots ,y(x)_{n(x)}$ in $T_x$ such that
$T_x = \bigcup_{i = 1}^{n(x)} W_{y(x)_i}$.
We can arrange $W_{gy} = gW_y$, $n(x) = n(gx)$ and
$y(gx)_i = gy(x)_i$ for $g \in G$ since $C_x$ is a
$\Fin$-subset of $X$ with $G_{C_x} = G_x$.
Obviously
\begin{eqnarray}
C_x(T_x;a)^\circ & = &\bigcup_{i=1}^{n(x)} C_x(W_{y(x)_i};a)^\circ.
\label{lem:A-B-C-cover-for-K:(5)}
\end{eqnarray}
Define
$$\Lambda = \left\{y_i(x) \mid x \in I,
         i \in \{1,2, \ldots, n(x)\}\right\}.$$
This is a cofinite $G$-set.
Define for $\lambda = y_i(x)$ in $\Lambda$
\begin{eqnarray*}
C_{\lambda} & = & C_x;
\\
B_{\lambda} & = & C_x(U_{y_i(x)};a+2b);
\\
A_{\lambda} & = & C_x(W_{y_i(x)};a).
\end{eqnarray*}
It remains to check that this collection of
 boxes has the desired properties.
This is obvious for assertions
\ref{lem:A-B-C-cover-for-K:cofinite},
\ref{lem:A-B-C-cover-for-K:long},
\ref{lem:A-B-C-cover-for-K:S_C is connected},
\ref{lem:A-B-C-cover-for-K:parallel}
and~\ref{lem:A-B-C-cover-for-K:equivariant}.
Assertion~\ref{prop:A_in B_inC} follows from
$W_y \subseteq U_y^{\circ}$ and
$U_y \subseteq S_{C_x} \cap C_x^{\circ}$.
Assertion~\ref{lem:A-B-C-cover-for-K:cover} follows
from \eqref{lem:A-B-C-cover-for-K:(1)}
and~\eqref{lem:A-B-C-cover-for-K:(5)}.

Finally we prove
assertion~\ref{lem:A-B-C-cover-for-K:intersection-control}.
Suppose that $B_{\lambda} \cap B_{\lambda'} \not= \emptyset$.
Write $\lambda = y_i(x)$ and $\lambda' = y_{i'}(x')$.
Since $B_{\lambda} = C_x(U_{y_i(x)};a+2b)$ and
$B_{\lambda'} \subseteq C_{x'}(T_{x'};a+2b)$, we conclude that
$ C_x(U_{y_i(x)};a+2b) \cap C_{x'}(T_{x'};a+2b) \not= \emptyset$.
By \eqref{lem:A-B-C-cover-for-K:(4)}
we have $x' \in I_{y_i(x)}$.
We conclude from~\eqref{lem:A-B-C-cover-for-K:(2)}
$$B_{\lambda} = C_x(U_{y_i(x)};a+2b)  \subseteq C_{x'}^{\circ}
                                  = C_{\lambda'}^{\circ}.$$
The central slice $S_{B_{\lambda}} = U_{y_i(x)}$ is transversal
to the flow with respect to $C_{\lambda'}$
by \eqref{lem:A-B-C-cover-for-K:(3)}.
This finishes the proof of Lemma~\ref{lem:A-B-C-cover-for-K}.
\end{proof}


\section{General position in metric spaces}
\label{sec:general-position-metric-spaces}

\begin{definition}
\label{def:thickenings-in-metric-space}
Let $Z$ be a metric space, $A \subseteq Z$ and $\delta > 0$.
Then we define the sets
\begin{align*}
A^\delta & = \{ x \in Z \mid \exists a \in A \;
        \mbox{such that} \; d(x,a) < \delta \},
\\
A^{-\delta} & = \{ x \in A \mid d(x,Z-A) > \delta\}.
\end{align*}
For $x \in Z$, we will abbreviate $x^\delta = \{ x \}^\delta$.
\end{definition}

Notice that $A^{\delta}$ and $A^{-\delta}$ are open.
The following
Propositions~\ref{pro:general-position-for-subsets} and
~\ref{prop:cover-Z-compatitible-with-calu} are the main
results of this section which is entirely devoted to their proof.

\begin{proposition}
\label{pro:general-position-for-subsets}
Let $Z$ be a compact metrizable space of covering dimension $n$
with an action of a finite group $F$.
Let $\calu$ be a finite collection of open $\Fin$-subsets
such that $gU \in \calu$ for $g \in F$, $U \in \calu$.
Suppose that we are given for each $U \in \calu$ an open subset
$U'' \subseteq U$
satisfying $\overline{U''} \subseteq U$.
Put $m = (n+1) \cdot |F|$.

Then for each $U \in \calu$ we can find an open
subset $U' \subseteq Z$ such that
\begin{enumerate}
\item \label{thm:general:delta-smaller}
      $U'' \subseteq \overline{U''} \subseteq U'
           \subseteq \overline{U'} \subseteq U$;
\item \label{thm:general:intersection}
      If $\calu_0 \subseteq \calu$ has more than $m$
      elements, then
      \begin{equation*}
      \bigcap_{U \in \calu_0} \dd U' = \emptyset;
      \end{equation*}
\item \label{thm:general:equivariant}
      $(gU)' = g(U')$ for $g \in F$, $U \in \calu$.
\end{enumerate}
\end{proposition}

\begin{proposition}
\label{prop:cover-Z-compatitible-with-calu}
Let $Z$ be a compact metric space of covering dimension $n$
with an action of a finite group $F$ by isometries.
Let $Y \subseteq Z$ be an open locally connected $F$-invariant subset.
Let $\calu$ be a finite collection of
open subsets such that
$gU \in \calu$ for $g \in F$, $U \in \calu$
and $\overline{U} \subseteq Y$ for $U \in \calu$ holds.
Assume that there is $k$ such that
for every subset $\calu_0 \subseteq \calu$ with
more than $k$ elements we have
\begin{equation*}
\bigcap_{U \in \calu_0} \dd U = \emptyset.
\end{equation*}
Let $\delta > 0$.
Put $m = (n+1) \cdot |F|$, $l = k|F|$.

Then  there are finite
collections $\calv^j$,
$j=0,\dots,m$ of open subsets of $Z$,
such that
\begin{enumerate}
\item \label{prop:cover-Z:cover}
      $\calv = \calv^0 \cup \dots \cup \calv^m$ is an open
      cover of $Z$;
\item \label{prop:cover-Z:diam-delta}
      $\diam(V) < \delta$ for every $V \in \calv$;
\item \label{prop:cover-Z:intersect-not-contained}
      For $V \in \calv$ there are at most
      $l$ different sets $U \in \calu$
      such that $U$ and $V$ intersect, but
      $U$ does not contain $V$;
\item \label{prop:cover-Z:intersect-V}
      For fixed $j$ and $V_0 \in \calv^j$ we have $V_0 \cap V \not= 0$ for at most $2^{j+1} - 2$
      different sets
      $V \in \calv^0 \cup \dots \cup \calv^j$;
\item \label{prop:cover-Z:mutually-disjoint}
      For fixed $j$ and $V_0,V_1 \in \calv^j$ we have either $V_0 = V_1$ or $\overline{V_0} \cap \overline{V_1} = \emptyset$;
\item \label{prop:cover-Z:F-invariant}
      Each $\calv^i$ is $F$-invariant, i.e.,
      $gV \in \calv^i$ for $g \in F$, $V \in \calv^i$;
\item \label{prop:cover-Z:V_circ=V} $V = \overline{V}^{\circ}$ for $V \in \calv$.
\end{enumerate}
\end{proposition}

In order to prove these two propositions
we will compare the metric space $Z$ to
the nerve of a suitable open cover of $Z$.
The results will be first proven for simplicial
complexes and then be pulled back to $Z$ using
the map from the next remark.

\begin{remark}
\label{rem:map-to-realization}
Let $Z$ be a metric space and $\calu$ be a
locally finite open cover of $Z$.
Denote by $\caln(\calu)$ the simplicial
complex given by the nerve of $\calu$.
Then
\begin{equation}
z \mapsto \sum_{U \in \calu}
     \left( \frac{d(z,Z-U)}{\sum_{V \in \calu} d(z,Z-V)}
                 \right) [U]
\end{equation}
defines a map $f_{\calu} \colon Z \to \caln(\calu)$, where
$[U]$ denotes the vertex of $\caln(\calu)$ that corresponds to
$U \in \calu$.
If a group $G$ acts by isometries on $Z$ and
on the cover $\calu$
(i.e., $gU \in \calu$ for $g \in G$, $U \in \calu$)
then $f_\calu$ is equivariant for the induced action
on $\caln(\calu)$.
\end{remark}

\begin{lemma}
\label{lem:for-beta-there-is-alpha}
Let $f \colon Y \to Z$ be a continuous map between
metric spaces.
Assume that $Y$ is compact.
Let $U$ be an open subset of $Z$ and $\beta > 0$.
Then there exists $\alpha > 0$ such that
\begin{equation*}
(f^{-1}(U))^{-\beta} \subseteq f^{-1}(U^{-\alpha}).
\end{equation*}
\end{lemma}

\begin{proof}
We proceed by contradiction and assume that
for every $n \in \IN$ there is $x_n \in Y$
such that
\begin{align}
\label{eq:x-is-in}
x_n & \in f^{-1}(U)^{-\beta},
\\
\label{eq:x-is-not-in}
x_n & \not\in f^{-1}(U^{-\frac{1}{n}}).
\end{align}
By compactness of $Y$, we may assume that
$x_n \to x$ as $n \to \infty$.
We derive from~\eqref{eq:x-is-not-in}
that there is $z_n \in Z - U$ such that
$d(z_n, f(x_n)) \le \frac{2}{n}$.
Then $z_n \to f(x)$ as $n \to \infty$ and
therefore $f(x) \not\in U$.
On the other hand \eqref{eq:x-is-in}
implies that $f((x_n)^\beta) \subseteq U$.
Since $x \in (x_n)^\beta$ for sufficiently large
$n$ we have $f(x) \in U$.
This is the desired contradiction.
\end{proof}

In the sequel interior of a simplex means simplicial
interior, i.e.,  the simplex with all its
proper faces removed.

\begin{lemma}
\label{lem:simplices-of-subdivision}
Let $Z$ be a simplicial complex and
let $Z^{(n)}$ be the $n$-th
barycentric subdivision of $Z$, $n \geq 1$.
Let $\Delta$ be a simplex of $Z^{(n)}$. Let
$\sigma$, $\tau$ be simplices of $Z$ such that both
the interior of $\sigma$ and the interior of $\tau$
intersect $\Delta$.
Then $\sigma \subseteq \tau$ or
$\tau\subseteq \sigma$.
\end{lemma}

\begin{proof}
Let $\Delta'$ be a simplex of the first barycentric
subdivision $Z' = Z^{(1)}$
of $Z$ which contains $\Delta$.
Then both the interior of $\sigma$ and the interior
$\tau$  intersect $\Delta'$.
Thus it suffices to prove the claim for the
first barycentric subdivision.

Now a $d$-simplex $\Delta'$
of the first barycentric subdivision $Z'$ of $Z$
is given by a sequence $\sigma_0, \sigma_1, \ldots, \sigma_d$
such that each $\sigma_i$ is a simplex of $Z$
and $\sigma_{i-1}$ is a proper face of $\sigma_i$
for $i = 1,2, \ldots, d$.
This is the simplex in the barycentric subdivision
whose vertices are the barycenters
of $\sigma_0, \sigma_1, \ldots, \sigma_d$.
Then the simplices of $Z$ whose interior intersect $\Delta'$
are precisely the $\sigma_i$.
\end{proof}

Recall that the \emph{open star} of a vertex $e$
of a simplicial complex $Z$
is defined to be the set of all points $z \in Z$
whose carrier simplex has $v$ as an vertex.
Equivalently, one obtains the
star of $v$ by taking the union of all simplices
$\sigma$ which have
$v$ as vertex and then deleting those faces of these
simplices $\sigma$ which do not have $v$ as a vertex.
We will denote it by  $\stern(v)$.
The set of open stars of vertices of $Z$ is an open covering of $Z$.

Let $\sigma$ be a simplex of $Z$. Define $\stern'(\sigma)$ to be
the star in the first barycentric subdivision $Z'$ of the vertex
in $Z'$ given by $\sigma$.

\begin{proof}
[Proof of Proposition~\ref{prop:cover-Z-compatitible-with-calu}]
Since $\calu$ is finite and $\partial U$
closed for $U \in \calu$,
the assumptions on $\calu$ imply that we can find for
every $z \in Z$ an open neighborhood $W_z$ of $z$
such that for all $z \in Z$ the following holds.
\begin{itemize}
\item $\diam(W_z) < \frac{\delta}{|F|}$;

\item $W_z$ intersects the boundary of at most
      $k$ sets in $\calu$;

\item If $W_z$ intersects $\overline{U}$ for some
      $U \in \calu$, then $W_z \subseteq Y$.
\end{itemize}

Since the covering dimension of $Z$ is $n$ by assumption
and $Z$ is compact,
we can choose
a finite open  refinement $\calw$ of the open cover
$\{ W_z \; | \; z \in Z \}$ such that $\dim (\calw) \leq n$.
Let  $\calw_F$ be the collection of subsets of $Z$ which consists
of the components of subsets of the form
$\bigcup_{g \in F} gW$ for $W \in \calw$
with $W \subseteq Y$ and subsets of the form $g W$
for $g \in F$ and $W \in \calw$ with $W \not\subseteq Y$.
Since $Y$ is a locally connected open subset of $Z$,
the above components are open subsets of $Z$.
Thus the elements in the finite set $\calw_F$ are open subsets.
Hence $\calw_F$ is a finite open covering of $Z$.
Since $\dim (\calw) \le n$ every orbit of $F$ in $Z$ meets
at most $(n+1) \cdot |F|$ members of $\calw$.
We conclude
$\dim (\calw_F) \leq m = (n+1) \cdot |F|$.
Obviously $g V \in \calw_F$ holds if $V \in \calw_F$.
If $V$ is a component of
$\bigcup_{g \in F} gW$ for some $W \in \calw$, we can find
elements $g_1, g_2, \ldots, g_r$ in $F$ such that
$V$ is contained in $\bigcup_{i = 1}^r g_iW$
and $\left(\bigcup_{j=1}^i g_j W\right) \cap g_{i+1} W \not= \emptyset$
holds for
$i =1,2, \ldots ,(r-1)$.
One shows by induction over
$i = 1,2 \ldots, r$ that the diameter of $\bigcup_{j=1}^i g_j W$ is less or
equal to the sum of the diameters
of the sets $g_1 W$, $g_2 W$,\dots , $g_i \cdot W$.
Hence the diameter of any element
$V$ of $\calw_F$ is bounded by $\delta$.

Consider $U \in \calu$ and an element $V \in \calw_F$
such that $V$ intersects $U$ but is not contained in $U$.
By construction $V$ is a component
of $\bigcup_{g \in F} gW$ for some $W \in \calw$ with $W \subseteq Y$.
Since $V$ is connected we must have
$V \cap \partial U \not= \emptyset$.
Hence there exists $g \in F$ with $gW
\cap \partial U \not= \emptyset$, or, equivalently, with $W \cap
g^{-1} \partial U \not= \emptyset$. Since each set $W_z$
intersects the boundary of at most $k$ sets in $\calu$ and $g^{-1}
U \in \calu$, there are at most $l = |F| \cdot k$ elements $U \in
\calu$ satisfying $W \cap g^{-1} \partial U \not= \emptyset$
for some $g \in F$.
Hence for every $V \in \calw_F$ there are at most $l$ elements $U
\in \calu$ such that $V$ intersects $U$ but is not contained in
$U$.

Consider the map $f = f_{\calw_F} \colon Z \to \caln(\calw_F)$
from Remark~\ref{rem:map-to-realization}. Put for $j = 0,1, \ldots, m$
$$\calx^{j} ~ = ~ \{\stern'(\sigma) \mid \sigma \in \caln(\calw_F), \dim(\sigma) = j\}.$$
Then $\calx = \calx^{0} \cup \calx^{1} \cup \ldots \cup \calx^{m}$
is an open cover of $\caln(\calw_F)$.

Let $\sigma$ and $\tau$ be simplices in $Z$ with
$\stern'(\sigma) \cap \stern'(\tau) \not= \emptyset$.
Let $\Delta$ be the carrier simplex in $Z'$ of some point in
$\stern'(\sigma) \cap \stern'(\tau)$.
Then the two vertices of $Z'$ given by $\sigma$ and $\tau$ belong to
$\Delta$. Hence the barycenters of $\sigma$ and $\tau$ belong to
$\Delta$. In particular the interior of $\sigma$ and the interior of $\tau$ intersect
$\Delta$. We conclude from  Lemma~\ref{lem:simplices-of-subdivision}
that $\sigma \subseteq \tau$ or $\tau \subseteq \sigma$ holds.
This implies that the elements in
$\calx^{j}$ are pairwise  disjoint and every $V \in \calx^{j}$ intersects at most $2^{j+1} -2$ elements in
$\{\calx^{0} \cup \calx^{1} \cup \ldots \cup \calx^{j}\}$
since a $j$-simplex has $2^{j+1} -2$ proper faces.

If the simplex $\sigma$ in $\caln(\calw_F)$ is given by the subset $\{V_0,V_1, \ldots ,V_d\} \subseteq \calw_F$,
then $f^{-1}(\stern'(\sigma))$ is contained in $V_0 \cap V_1 \cap \ldots \cap V_d$,
since $\stern'(\sigma)$ is contained in the intersection of the stars in $\caln(\calw_F)$
of the vertices of $\sigma$ and for $V \in \caln(\calw_F)$ we have $f^{-1}(\stern(V)) \subseteq V$.
Hence for every $X \in \calx$ there exists $W \in \calw_F$ with
$f^{-1}(X) \subseteq W$. Put
$$\widehat{\calv}^{j} = \{ f^{-1}_{\calw_F}(X) \; | \; X \in \calx^j \}.$$
Then $\widehat{\calv}^{j}$ has the properties
\begin{itemize}
\item $\widehat{\calv} = \widehat{\calv}^0 \cup \dots \cup \widehat{\calv}^m$ is an open
      cover of $Z$ consisting of finitely many elements.
\item $\diam \widehat{V} < \delta$ for every $\widehat{V} \in \widehat{\calv}$;
\item For $\widehat{V} \in \widehat{\calv}$ there are at most
      $l$ different sets $U \in \calu$
      such that $U$ and $\widehat{V}$ intersect, but
      $U$ does not contain $\widehat{V}$;
\item For fixed $j$, every $\widehat{V}_0 \in \widehat{\calv}^j$ intersects at most
      $2^{j+1} - 2$
      different sets
      $\widehat{V} \in \widehat{\calv}^0 \cup \dots \cup \widehat{\calv}^j$;
\item For fixed $j$ and $\widehat{V}_0,\widehat{V}_1 \in
      \widehat{\calv}^j$
      we have either $\widehat{V}_0 = \widehat{V}_1$ or $\widehat{V}_0 \cap \widehat{V}_1 = \emptyset$;
\item Each $\widehat{\calv}^i$ is $F$-invariant, i.e.,
      $g\widehat{V} \in \widehat{\calv}^i$ for $g \in F$, $\widehat{V} \in \widehat{\calv}^i$.
\end{itemize}
For $\e > 0$ define
$$\calv^j = \left\{\left(\overline{(\widehat{V})^{-\e}}\right)^{\circ} \mid \widehat{V} \in
\widehat{\calv}^{j}\right\}.$$ Then $\calv^j$ has
properties~\ref{prop:cover-Z:diam-delta},
\ref{prop:cover-Z:intersect-not-contained},
\ref{prop:cover-Z:intersect-V},
\ref{prop:cover-Z:mutually-disjoint},
\ref{prop:cover-Z:F-invariant} and~\ref{prop:cover-Z:V_circ=V}
since $\left(\overline{(\widehat{V})^{-\e}}\right)^{\circ} \subseteq \widehat{V}$
for $\widehat{V} \in \calv^{j'}$ and for any open subset $T$ of a topological
space we have $\overline{T} = \overline{\overline{T}^{\circ}}$.
Since $Z$ is compact, we can
choose $\e$ so small that also property
~\ref{prop:cover-Z:cover} is satisfied. This finishes the proof of
Proposition~\ref{prop:cover-Z-compatitible-with-calu}.
\end{proof}

\begin{proposition}
\label{prop:subdivision-and-general-position}
Let $Z$ be a simplicial complex and let
$Z^{(m)}$ be the $m$-th
barycentric subdivision of $Z$, $m \geq 1$.
Let $B$ be a subcomplex of $Z$.
Let $B_0$ be the union of all simplices of $Z^{(m)}$
that are contained in the interior of $B$.
Let $\beta$ be a simplex of $Z$
and let $\sigma$ be a simplex of $Z^{(m)}$
that is contained in the intersection of
the boundary of $B_0$ with $\beta$.
Then $\dim (\sigma) < \dim (\beta)$.
\end{proposition}

We recall that the topological interior and boundary of
a subcomplex  of a simplicial complex can be described
combinatorial as follows: the interior is the union of all
open simplices (simplices with all proper faces removed)
that are contained in the subcomplex;
its boundary is the union of all simplices that are contained in
the subcomplex and are in a addition a face of simplex not contained in
the subcomplex.

\begin{proof}
[Proof of Proposition~\ref{prop:subdivision-and-general-position}]
By the definition of $B_0$ there exists a simplex $\Delta$
of $Z^{(m)}$ such that $\sigma \subseteq\Delta$, but $\Delta$
is not contained in the interior of $B$.
Thus there exists a simplex $\alpha$ of $Z$ that is
contained in the boundary of $B$ and intersects $\Delta$.
By passing to a face of $\alpha$, we can arrange
that $\alpha$ is contained in the boundary of $B$ and the interior of
$\alpha$ intersects $\Delta$. Since $\sigma$ is contained in $\beta$,
we can find a face $\beta'$ of $\beta$ such that $\sigma$ intersects
the interior of $\beta'$. Hence the interior of the simplex $\beta'
\subseteq Z$ intersects $\Delta$.
Lemma~\ref{lem:simplices-of-subdivision} implies
$\alpha \subseteq\beta'$ or $\beta' \subseteq\alpha$.
In the second case,
$\sigma \subseteq\beta' \subseteq\alpha$.
But $\sigma$ has to be disjoint from $\alpha$,
because $\alpha$ lies on the boundary of $B$, while
$\sigma$ is contained in the interior of $B$.
We conclude $\alpha \subseteq\beta'$ and hence $\alpha \subseteq \beta$.
Therefore, $\tau = \beta \cap \Delta$
is a simplex of $Z^{(m)}$, that contains the simplex $\sigma$ as a face
and intersects $\alpha$.
Since $\sigma $ and $\alpha$ are
disjoint, $\sigma$ is a proper face of $\tau$. This implies
$\dim (\sigma) < \dim (\tau) \leq \dim (\beta)$.
\end{proof}

\begin{proof}
[Proof of Proposition~\ref{pro:general-position-for-subsets}]
The strategy is  first to prove a simplicial version and then use
the map appearing in Remark~\ref{rem:map-to-realization} to
handle the general case of a metric space.

Since we assume that $Z$ is metrizable and $F$ is finite, we can choose a
metric $d_Z$ on $Z$ which is $F$-invariant. Since $Z$ and hence each $\overline{U''}$ is compact and
the collection $\calu$ is finite, we can find $\delta > 0$ such that $\overline{U''} \subseteq U^{-\delta}$ holds for
$U \in \calu$. Hence we can assume in the sequel without loss of generality that $U'' = U^{-\delta}$.

So we start with the special case where $Z$ is in addition
a simplicial complex, each $U \in \calu$ is the interior
of a subcomplex of $Z$ and $F$ acts simplicially on $Z$.
Let $\{ U_1,\dots,U_r \} \subseteq \calu$ contain
exactly one element from every $F$-orbit in the $F$-set $\calu$.
Pick $m \geq 0$ such that the simplices of the $m$-th
barycentric subdivision $Z^{(m)}$ of $Z$ have diameter $< \delta$.
For $i = 1,2, \ldots , r$ let $Z^{(m+i)}$ be the $(m+i)$-th
barycentric subdivision of $Z$ and
let $A_i$ be the union of all those simplices of $Z^{(m+i)}$ which are
contained in $U_i$. Then $A_i$ is the largest subcomplex of $Z^{(m+i)}$
that is contained in $U_i$. Since each simplex of $Z^{(m+i)}$
has diameter $< \delta$, we get $U_i^{-\delta} \subseteq A_i
\subseteq U_i$. Define $U_i'$ to be the interior of $A_i$.
For $U \in \calu$ define $U' = gU_i'$ for any choice of $g \in F$ and $i
\in \{1,2, \ldots ,r\}$ satisfying $U = gU_i$.
One easily checks that $i$ is uniquely determined by $U$ and the
choice of $g \in F$ does not matter in the definition of $U'$
since each $U$ is by assumption a $\Fin$-subset.
Obviously \ref{thm:general:delta-smaller}
and \ref{thm:general:equivariant} are satisfied
for $\calu' = \{ U' \mid U \in \calu \}$.

Next we  prove \ref{thm:general:intersection} but with $m$ replaced
by $n = \dim(Z)$.
Consider a subset $\calu_0 \subseteq \calu$ with $|\calu_0| = k$.
Notice that $ \bigcap_{U \in \calu_0} \dd U'$ is a subcomplex
of $Z^{(m+a)}$ for some $a$ since the intersection of a subcomplex of
$Z^{(n+a)}$ with a subcomplex of $Z^{(n+b)}$ is a subcomplex of
$Z^{(n+b)}$ if $a \le b$.
It suffices to show
$\dim \left( \bigcap_{U \in \calu_0} \dd U' \right) \leq
        n - k$
since this implies
\begin{equation}
\label{eq:empty-boundary-intersection-for-n}
\bigcap_{U \in \calu_0} \dd U' = \emptyset
           \quad \mbox{if $\calu_0$ contains more
                                   than $n$ elements}.
\end{equation}
Choose $i_1, i_2, \ldots , i_k$ in $\{1,2,\ldots ,r\}$
and $g_1, g_2, \ldots , g_k$ in $F$ with the property that
$\calu_0$ consists of the mutually different elements $g_1U_{i_1}$,
$g_2U_{i_2}$, $\ldots,$ $ g_kU_{i_k}$ and
$i_1 \le i_2 \le \ldots \le i_k$ holds. If for some
$j \in \{1,2,\ldots ,(r-1)\}$ we have $i_j = i_{j+1}$, then
$g_jU_{i_j} \not= g_{j+1}U_{i_{j+1}}$ implies already
$g_jU_{i_j} \cap g_{j+1}U_{i_{j+1}} = \emptyset$ and the claim is
obviously true. Hence we can assume without loss of generality
$i_1 < i_2 < \ldots < i_k$. Next we show by induction for
$j = 1,2, \ldots ,k$ that
$$\dim \left( \bigcap_{l=1}^j \dd g_{i_l}U_{i_l}'  \right) \leq n - j.$$
The induction beginning is obvious since the dimension of the boundary of a simplicial
complex is smaller than the dimension of the simplicial complex
itself. The induction step from $j-1$ to $j$ is done as follows.
By induction hypothesis the dimension of the simplicial subcomplex
$\bigcap_{l=1}^{j-1} \dd g_{i_l}U_{i_l}'$ of $Z^{(n + i_{j-1})}$ is less or
equal to $(n-j +1)$. Let $\sigma$ be a simplex of the subcomplex
$\bigcap_{l=1}^j \dd g_{i_l}U_{i_l}'$ of $Z^{(n + i_j)}$. We can find a
simplex $\beta$ in $\bigcap_{l=1}^{j-1} \dd g_{i_l}U_{i_l}'$ such that
$\sigma$ is contained in $\beta$.  Recall that by assumption
$U_j$ is the interior of a subcomplex $B_j$ of $Z^{(n+i_{j-1})}$.
Proposition~\ref{prop:subdivision-and-general-position}
applied to the case  $m = i_j - i_{j-1}$ and $B = g_jB_j$ and
$\sigma$ and $\beta$ as above implies
$$\dim(\sigma) < \dim(\beta) \le n-j+1$$
since in the notation of
Proposition~\ref{prop:subdivision-and-general-position} we have $B_0
=g_jU_{i_j}'$. Hence $\dim(\sigma) \le n -
j$. This finishes the proof of
Proposition~\ref{pro:general-position-for-subsets} in the special case
where $Z$ is a simplicial complex, each $U \in \calu$ is the interior
of a subcomplex of $Z$ and $F$ acts simplicially on $Z$.

In the general case we start with an open cover $\calv$ of $Z$ such
that each $V \in \calv$ has diameter $< \epsilon = \frac{\delta}{3}$
and $\dim(\calv) \leq n$. Let $\calv_F$ be the cover of $Z$ whose
members are the open sets of the form $g V$ with $V \in \calv$ and
$g \in F$. Then $\calv_F$ is an open cover of $Z$ whose members have
diameter $< \epsilon$ and we have $g V \in \calv_F$ for $V \in
\calv_F$ and $g \in F$. Analogously as in
Proposition~\ref{prop:cover-Z-compatitible-with-calu} one shows that
$\dim (\calv_F) \leq m$.

For $U \in \calu$ let
$\calv_U = \{ V \in \calv_F \; | \; V \subseteq U \}$
and consider $\caln(\calv_U)$ as a subcomplex of $\caln(\calv_F)$
(this is the subcomplex spanned by the vertices $[V]$
with $V \in \calv_U$) and denote by $\widehat{U}$ the interior
of $\caln(\calv_U)$.
Consider the map $f = f_{\calv_F} \colon Z \to \caln(\calv_F)$
from Remark~\ref{rem:map-to-realization}.
If $x \in f^{-1}(\caln(\calv_U))$ and $x \in V$, $V \in \calv_F$
then by the construction of $f$, $V \in \calv_U$.
Therefore
\begin{equation*}
f^{-1}(\widehat{U}) \subseteq f^{-1}(\caln(\calv_U)) \subseteq U.
\end{equation*}
Let $x \in U^{-\epsilon}$.
If $x \in V$ with $V \in \calv_F$ then $V \in \calv_U$,
because the diameter of $V$ is  $< \epsilon$.
Therefore $f(x) \in \caln(\calv_U)$.
If $f(x)$ lies on the boundary of $\caln(\calv_U)$,
then there are $V \in \calv_U$, $V' \in \calv_F - \calv_U$
such that $x \in V$ and $V \cap V' \neq \emptyset$.
In particular this implies $x \not\in U^{-2\epsilon}$.
We have thus shown that
\begin{equation}
\label{eq:U-minus-2e-subset}
U^{-2 \epsilon} \subseteq f^{-1}(\widehat{U}) \subseteq U.
\end{equation}
Equip $\caln(\calv_F)$ with a metric such that the action of $F$ is by
isometries. By Lemma~\ref{lem:for-beta-there-is-alpha} there is
$\alpha$ such that
\begin{equation}
\label{eq:minus-e-subset-minus-alpha}
(f^{-1}(\widehat{U}))^{-\epsilon} \subseteq f^{-1}( \widehat{U}^{-\alpha} )
\end{equation}
for all $U \in \calu$.
By the special case treated in the first part of this
proof for each $\widehat{U}$ there is $\widehat{U}'$ such that
$\widehat{U}^{-\alpha} \subseteq \overline{\widehat{U}^{-\alpha}}
 \subseteq \widehat{U}' \subseteq \overline{\widehat{U}'}
\subseteq \widehat{U}$,
$(g \widehat{U})' = g (\widehat{U}')$ for $g \in F$ and
\begin{equation*}
\bigcap_{U \in \calu_0} \dd \widehat{U}' = \emptyset
      \quad \mbox{if $\calu_0$ contains more than $m$ elements}.
\end{equation*}
(See \eqref{eq:empty-boundary-intersection-for-n} and
recall that $\dim (\calv_F) \leq m$.)
Finally set $U' = f^{-1}(\widehat{U}')$ for $U \in \calu$.
Since $\dd f^{-1}(U) \subseteq f^{-1}(\dd U)$
and taking preimages preserves inclusions and
intersections \ref{thm:general:intersection}
and \ref{thm:general:equivariant} are satisfied.
Moreover, by
\eqref{eq:U-minus-2e-subset} and
\eqref{eq:minus-e-subset-minus-alpha}
\begin{equation*}
U^{-\delta} = U^{-3\epsilon} \subseteq (f^{-1}(\widehat{U}))^{-\epsilon}
 \subseteq f^{-1}(\widehat{U}^{-\alpha}) \subseteq f^{-1}( \widehat{U}' )
 \subseteq f^{-1}(\widehat{U}) \subseteq U
\end{equation*}
and therefore \ref{thm:general:delta-smaller}
is also satisfied.
\end{proof}


\section{Covering $X_{>\gamma}$ by long boxes}
\label{sec:covering-G->gamma}

Throughout this Section we will assume that we
are in the situation of
Convention~\ref{conv:for-long-and-thin-covers}.
In particular $k_G$ is the maximum over the orders of
finite subgroups of $G$ and $d_X$ is the dimension
of $X - X^\IR$.
Both $k_G$ and $d_X$ are finite.
Also recall the notation $X_{>r}$ introduced in
Definition~\ref{def:G-period}.
This section is entirely devoted to the proof of the following

\newlength{\origlabelwidth}
\setlength\origlabelwidth\labelwidth

\begin{proposition}
\label{prop:cover-K-by-boxes-hyp}
\setlength\labelwidth\origlabelwidth
There exists a natural number $M = M(k_G,d_X)$
depending only on $k_G$ and $d_X$
which has the following property:

For every $\alpha, \epsilon \in \IR$ with $0 < \epsilon < \alpha$
there exists $\gamma = \gamma(\alpha,\epsilon,M) > 0$
such that for every cocompact $G$-invariant subset $K$ of
$X_{> \gamma}$, there is a collection $\cald$ of boxes satisfying
\begin{numberlist}
\item[\label{prop:cover-K-hyp:cover}]
      $K \subseteq \bigcup_{D \in \cald} \Phi_{(-\epsilon,\epsilon)} (D^{\circ})$;
\item[\label{prop:cover-K-hyp:long}]
      For every $x \in X$ which lies on the open
      bottom or open top of a box in $\cald$,
      the set
      $\Phi_{[-\alpha,-\epsilon] \cup [\epsilon,\alpha]} (x)$
      does not intersect the open bottom or the open top of
      a box in $\cald$;
\item[\label{prop:cover-K-hyp:each-box-long}]
      For any $x \in X$ there is no box $D \in \cald$ such that
      $\Phi_{[0,\alpha]} (x)$ intersects both the
      open bottom and open top of $D$;
\item[\label{prop:cover-K-hyp:dim}]
      The dimension of the collection $\{D^{\circ} \mid D \in \cald\}$
      is less or equal to $M$, i.e. the intersection of
      $(M+2)$ pairwise distinct elements is always empty;
\item[\label{prop:cover-K-hyp:equivariant}]
      For $g \in G$, $D \in \cald$ we have $g D \in \cald$;
\item[\label{prop:cover-K-hyp:cofinite}]
      There is a finite subset $\cald_0 \subseteq \cald$
      such that for
      every $D \in \cald$ there exists $g \in G$ with
      $gD \in \cald_0$;
\item[\label{prop:cover-K-hyp:Fin-subset}]
      $\Phi_{[-\alpha-\epsilon,\alpha + \epsilon]}(D)$ is a
      $\Fin$-subset of $X$ for all $D \in \cald$.
\end{numberlist}
\end{proposition}

The idea of the proof is very roughly as follows.
Conditions \eqref{prop:cover-K-hyp:long} and
\eqref{prop:cover-K-hyp:each-box-long} require the boxes
to be very long, but we have still the freedom to
make the boxes very thin.
Proposition~\ref{pro:general-position-for-subsets}
applied to the transversal directions
will be important to arrange the boxes during the construction
to be in general position.
This will allow the application of
Proposition~\ref{prop:cover-Z-compatitible-with-calu},
where property~\ref{prop:cover-Z:intersect-not-contained}
will be crucial in order control how many boxes from previous
steps of the construction interfere at each step in the
construction.


\subsection{Preliminaries and the basic induction structure}
\label{subsec:preliminaries-and-basic-induction-structure}

We begin with fixing some numbers and collection of boxes.
Define numbers
\begin{eqnarray*}
m & := & k_G \cdot (d_X + 1);
\\
M & := &  (k_G)^2 \cdot (d_X + 1) + 2^{m+1};
\\
a & := & \e / 2;
\\
b & := & 4M \cdot (\alpha + 2\epsilon)
                      + 3(\alpha + \epsilon);
\\
c & := & a + 2(b+\e) + 1 + \e;
\\
\gamma & := & a + 2 b + 2c.
\end{eqnarray*}

Notice that $m$ and $M$ depend only on $k_G$ and $d_X$ and all
the other numbers depend only on $\alpha$, $\epsilon$ and $M$.
(The reader may wonder why we picked $a$ small, we are after
all looking for long boxes.
But for our construction it is only important that $b$ and $c$
are large and in the proof of
Lemma~\ref{lem:good-things-about-D-W}~
\ref{lem:good-things-D-W:each-long}
our choice of a small $a$ will be convenient.)

Let $A_\lambda \subseteq B_\lambda \subseteq C_\lambda$ for
$\lambda \in \Lambda$ be three collection of boxes as in the
assertion of Lemma~\ref{lem:A-B-C-cover-for-K}, where we
use $a$  as defined above and replace ${b}$ by $b + \epsilon$,
and $c$ by $c - \epsilon$.
Then we replace $B_{\lambda}$ by the restriction ${B}_{\lambda}(a + 2b)$.
The resulting collections satisfy

\begin{itemize}
\item $\Lambda$ is $G$-cofinite;

\item $l_A = a$, $l_B = a + 2b$ and $l_C = a + 2b +2c$;

\item $S_{C_{\lambda}}$ is connected;

\item $S_{A_\lambda} \subseteq S_{B_\lambda}
                         \subseteq S_{C_{\lambda}}$;

\item $A_{\lambda} \subseteq B_{\lambda}^{\circ}$ and
      $B_{\lambda} \subseteq C_{\lambda}^{\circ}$;

\item $K \subseteq \bigcup_{\lambda \in \Lambda}
                                   A_{\lambda}^{\circ}$;

\item $g A_\lambda = A_{g \lambda}$,
      $g B_\lambda = B_{g \lambda}$ and
      $g C_\lambda = C_{g \lambda}$
      for $g \in G$ and $\lambda \in \Lambda$;

\item If $B_\lambda \cap B_{\lambda'} \neq \emptyset$, then
      $\phi_{[-\epsilon,\epsilon]}(B_\lambda)
                        \subseteq C_{\lambda'}^{\circ}$
      and $S_{B_{\lambda}}$ is transversal to the flow with
      respect to $C_{\lambda'}$.

\end{itemize}

Next we discuss the main strategy of the proof.
We will construct the desired collection $\cald$ inductively
over larger and larger parts of $K$.
The boxes $C_\lambda$ will be used to control
the group action, i.e., they will allow us to restrict
attention to the finite group $G_{C_\lambda}$.
The boxes $B_\lambda$ will be used to
control properties \eqref{prop:cover-K-hyp:long} and
\eqref{prop:cover-K-hyp:each-box-long}.
The boxes $A_\lambda$ will be used to control which
part of $K$ is already covered.
We will need to sharpen the induction assumptions.
We introduce a minor but useful variation
of \eqref{prop:cover-K-hyp:long} as follows.

\begin{definition} \label{def:delta-overlong}
A collection of boxes is \emph{$\delta$-overlong} for
$0 \leq \delta < \epsilon$ if for every $x \in X$ which lies
on the open bottom or open top of a box in this collection,
the set
$\Phi_{[-\alpha-\delta,-\epsilon + \delta]
         \cup [\epsilon - \delta,\alpha + \delta]} (x)$
does not intersect the open bottom or the open top of
any box in this collection.
\end{definition}

The assertion \eqref{prop:cover-K-hyp:long}
is then that $\cald$ is $0$-overlong.
Clearly, $\delta$-overlong for some
$0 \leq \delta < \epsilon$ implies $0$-overlong.

\begin{definition} \label{def:not-huge}
We say that a box $D$ is
\emph{not huge}
if for every $\lambda \in \Lambda$ such that
$D$ intersects $B_\lambda$
we have
$\phi_{[-\epsilon,\epsilon]} (D) \subset C_\lambda^{\circ}$
and
$S_D$ is transversal to the flow with respect to $C_{\lambda}$.
\end{definition}

Every box which is obtained from one of the boxes
$B_{\lambda}$ by restriction is
automatically not huge.

\begin{definition}
We say a collection $\cale$ of boxes is a
\emph{$\delta$-good box cover}
of a subset $S$ of $X$ if it has the following properties:

\begin{itemize}

\item $\cale$ is $\delta$-overlong;

\item Every box in $\cale$ is not huge;

\item Assertions \eqref{prop:cover-K-hyp:each-box-long},
                 \eqref{prop:cover-K-hyp:dim},
                 \eqref{prop:cover-K-hyp:equivariant},
                 \eqref{prop:cover-K-hyp:cofinite} and
                 \eqref{prop:cover-K-hyp:Fin-subset}
                 hold for $\cale$ in place of $\cald$;

\item $S \subseteq \bigcup_{E \in \cale}
              \Phi_{(-\epsilon ,\epsilon)} (E^{\circ})$.

\end{itemize}
\end{definition}

To prove Proposition~\ref{prop:cover-K-by-boxes-hyp} we will
construct a $0$-good box cover of $K$.
Let $N$ be the number of $G$-orbits of $\Lambda$.

Put
\begin{equation}
\label{eq:delta-r}
\delta_r ~ = ~ \frac{N-r}{N+1} \cdot \epsilon
                      \quad \text{for} \quad r= 0,1, \dots,N.
\end{equation}
Clearly,
$$\e > \frac{N}{N+1}\e = \delta_0 > \dots > \delta_r >
        \delta_{r+1} > \dots > \delta_N = 0.$$
We will show inductively for $r = 0,1,2, \ldots, N$ that for any
subset  $\Xi \subseteq \Lambda$ consisting of
$r$ $G$-orbits there exists a $\delta_r$-good box cover of
$K_\Xi := \bigcup_{\xi \in \Xi} A_\xi$.
The induction beginning $r=0$ is trivial, take $\cald = \emptyset$.
The induction step from $r$ to $r+1$
is summarized in the next lemma.

\begin{lemma}[Induction step : $r$ to $r + 1$]
\label{lem:induction_step_for_prop}
Let $\Xi \subseteq \Lambda$ consist of $r$ $G$-orbits
and assume that $\cald$ is a $\delta_r$-good box cover of
$K_\Xi = \bigcup_{\xi \in \Xi} A_\xi$.
Let $\lambda \in \Lambda - \Xi$ and $\Xi' := G \lambda \cup \Xi$.
Then there is a $\delta_{r+1}$-good box cover $\cald'$ of
$K_{\Xi'} = \bigcup_{\xi \in \Xi'} A_\xi$.
\end{lemma}

Clearly Lemma~\ref{lem:induction_step_for_prop} implies
Proposition~\ref{prop:cover-K-by-boxes-hyp}.
The proof of Lemma~\ref{lem:induction_step_for_prop}
will occupy the remainder of
Section~\ref{sec:covering-G->gamma}.


\subsection{Boxes in $C_{\lambda}^{\circ}$}

Before we explain the construction of $\cald'$, we have to
introduce some notation and to rearrange $\cald$ as follows.
In the sequel everything will take place in the interior of
the box $C_{\lambda}$.
Recall for the sequel the $G_{C_{\lambda}}$-homeomorphism
\begin{eqnarray*}
& & \mu_{C_{\lambda}} \colon
  S_{C_{\lambda}} \times [-a/2-b-c,a/2+b+c] \to
  C_{\lambda}, \quad (x,\tau) ~ \mapsto ~ \Phi_{\tau}(x)
\end{eqnarray*}
from
Lemma~\ref{lem:basics_about_boxes}
~\ref{lem:basics_about_boxes:shape_of_B_in_terms_of_central_slice}.
Let
$\pi_{C_{\lambda}} \colon C_{\lambda} \to S_{C_{\lambda}},
             x \mapsto \Phi_{\frac{a_+(x) + a_-(x)}{2}}(x)$
be the retraction onto the central slice.
Closures and interiors of subsets of
$S_{C_{\lambda}}$ are always understood with respect to
$S_{C_{\lambda}}$.

\begin{definition} \label{def:U_E}
Let $E\subseteq C_{\lambda}^{\circ}$ be a box such that $S_E$ is
transversal to the flow with respect to $C_{\lambda}$.

Define an open  subset of
$S_{C_{\lambda}} \cap C_{\lambda}^{\circ}$ by
$$U_E ~ := ~ \pi_{C_{\lambda}}(S_E \cap E^{\circ}) =
                         \pi_{C_{\lambda}}(E^{\circ}).$$

Define the
continuous $G_E$-invariant map
$$\tau_E \colon \pi_{C_{\lambda}}(S_E) \to [-l_C/2,l_C/2]$$
to be the composite of the inverse of
$\pi_{C_{\lambda}}|_{S_E} \colon S_E \xrightarrow{\cong}
                                     \pi_{C_{\lambda}}(S_E)$,
the inverse of $\mu_{C_{\lambda}}$ restricted to $S_E$ and
the projection
$S_{C_{\lambda}} \times [-l_C/2,l_C/2] \to [-l_C/2,l_C/2]$.

For a subset $T \subseteq S_{C_{\lambda}}$
define a subset of $S_E$ by
$$\sigma_E(T) ~ := \pi_{C_{\lambda}}^{-1}(T) \cap S_E.$$
\end{definition}

\begin{lemma} \label{lem:properties_of_U_E_and_sigma_E(T)}
Let $E,E' \subseteq C_{\lambda}^{\circ}$ be boxes
such that $S_E$ and $S_{E'}$ are transversal to the flow
with respect to $C_{\lambda}$.
Then
\begin{enumerate}
\item
\label{lem:properties_of_U_E_and_sigma_E(T):G_E_subsetG_{C_lambda}}
If $gE \cap E' \not= \emptyset$ for some $g \in G$,
then $g \in G_{C_{\lambda}}$.
In particular $G_E$ is a subgroup of $G_{C_{\lambda}}$;
\item
\label{lem:properties_of_U_E_and_sigma_E(T):characterization_of_tau_E}
The map $\tau_E$ is uniquely characterized by
$$\mu_{C_{\lambda}}(x) =
   (\pi_{C_{\lambda}}(x),\tau_E \circ\pi_{C_{\lambda}}(x))$$
for $x \in S_E$;
\item
\label{lem:properties_of_U_E_and_sigma_E(T):E(T)_is_defined}
If $T \subseteq S_{C_{\lambda}}$ is a closed $\Fin$-subset of the
$G_E$-space $S_{C_{\lambda}}$, then $\sigma_E(T)$ is a closed
$\Fin$-subset of the $G_E$-space $S_E$ and the restriction
$E(\sigma_E(T))$ is defined;
\item
\label{lem:properties_of_U_E_and_sigma_E(T):identifyingE(T)^circ}
Let $U \subseteq S_{C_{\lambda}}$ be an open $\Fin$-subset of the
$G_E$-space $S_{C_{\lambda}}$ with
$\overline{U} \subseteq \pi_{C_{\lambda}}(S_E)$
and $U = \overline{U}^{\circ}$.
Then
\begin{eqnarray*}
E(\sigma_E(\overline{U})) & = &
       \Phi_{[-l_E/2,l_E/2]}(\sigma_E(\overline{U}));
\\
E(\sigma_E(\overline{U}))^{\circ} & = &
       \Phi_{(-l_E/2,l_E/2)}(\sigma_E(U));
\\
U_{E(\sigma_E(\overline{U}))} & = & U.
\end{eqnarray*}
\end{enumerate}
\end{lemma}

\begin{proof}
\ref{lem:properties_of_U_E_and_sigma_E(T):G_E_subsetG_{C_lambda}}
Suppose $gE \cap E' \not= \emptyset$.
Since then
$gC_{\lambda} \cap C_{\lambda} \not= \emptyset$,
we get $g \in G_{C_{\lambda}}$.
\\[1mm]
\ref{lem:properties_of_U_E_and_sigma_E(T):characterization_of_tau_E}
This follows from the definitions.
\\[1mm]
\ref{lem:properties_of_U_E_and_sigma_E(T):E(T)_is_defined}
This is obvious since $\pi_{C_{\lambda}}$ is $G_E$-equivariant.
\\
\ref{lem:properties_of_U_E_and_sigma_E(T):identifyingE(T)^circ}
We conclude
$E(\sigma_E(\overline{U}))  =
    \Phi_{[-l_E/2,l_E/2]}(\sigma_E(\overline{U}))$
from the definition of
the restriction.
The set $\Phi_{(-l_E/2,l_E/2)}(\sigma_E(U))$ is mapped under the
homeomorphism $\mu_{C_{\lambda}}$ to the set
$\{(u,t) \mid u \in U, |t - \tau_E(u)| < l_C/2\}$.
Since $\tau_E$ is continuous
and $U$ is open in $S_{C_{\lambda}}$, this set and hence
$\Phi_{(-l_E/2,l_E/2)}(\sigma_E(U)) \subseteq C_{\lambda}^{\circ}$
are open.
This implies
\begin{eqnarray*}
\Phi_{(-l_E/2,l_E/2)}(\sigma_E(U)) & \subseteq &
                     E(\sigma_E(\overline{U}))^{\circ};
\\
\sigma_E(U) & \subseteq &
    \sigma_E(\overline{U}) \cap E(\sigma_E(\overline{U}))^{\circ}.
\end{eqnarray*}
If we apply $\pi_{C_{\lambda}}$ to the latter inclusion,
we conclude
$$U \subseteq \overline{U} \cap
        \pi_{C_{\lambda}}(E(\sigma_E(\overline{U}))^{\circ})
    \subseteq \overline{U}^{\circ} = U.$$
This implies
$\sigma_E(U) = \sigma_E(\overline{U})
       \cap E(\sigma_E(\overline{U}))^{\circ}$
and $U_{E(\sigma_E(\overline{U}))} = U$.
Lemma~\ref{lem:basics_about_boxes}
~\ref{lem:basics_about_boxes:shape_of_B_in_terms_of_central_slice}
implies
$\Phi_{(-l_E/2,l_E/2)}(\sigma_E(U)) =
            E(\sigma_E(\overline{U}))^{\circ}$.
\end{proof}


\subsection{Rearranging the data of the induction beginning}

Let $\cald$ be the collection of boxes,
from the hypothesis of
Lemma~\ref{lem:induction_step_for_prop}
(with respect to $\Xi \subset \Lambda$ and $\lambda \in \Lambda - \Xi$).

\begin{definition} \label{def:cald_lambda_and calu_lambda}
Put
\begin{eqnarray*}
\cald_{\lambda} & := & \{ D \in \cald \mid D \cap B_{\lambda}
                     \neq \emptyset \};
\\
\calu(D_\lambda) & := & \{U_D \; | \; D \in \cald_{\lambda} \}.
\end{eqnarray*}
\end{definition}

Since the $G$-action on $X$ is proper and
$B_{\lambda}$ is compact,
property~\eqref{prop:cover-K-hyp:cofinite} implies
that $\cald_\lambda$ is finite.
We will consider the $G_{C_{\lambda}}$ sets
\begin{align*}
G_{C_{\lambda}}\cdot \cald_{\lambda}
 & = \{gD \mid D \in \cald_{\lambda}, g \in G_{C_{\lambda}}\};
\\
G_{C_{\lambda}}\cdot  \calu(\cald_\lambda)
 & = \{gU_D \mid D \in \cald_{\lambda}, g \in G_{C_{\lambda}}\}
   = \{U_E \mid E \in G_{C_{\lambda}} \cdot \cald_{\lambda} \}.
\end{align*}
Since every $D \in\cald_{\lambda}$ is not huge,
the set $D$ is contained in $C_{\lambda}^{\circ}$ and
$S_D$ is transversal to the flow with respect to $C_{\lambda}$.
Both these properties also hold for every
$D \in G_{C_{\lambda}}\cdot \cald_{\lambda}$.

We use Proposition~\ref{pro:general-position-for-subsets}
to diminish the elements in $\cald$ slightly
in order to obtain a general position property for $\cald$.
At this point it is important, that we arranged the central slice
$S_{C_\lambda}$ to be connected.

We main goal of this subsection will be the proof of
the following lemma.

\begin{lemma} \label{lem:additional_properties_of_cald}
We can assume without loss of generality that $\cald$ has the
following general position properties
\begin{numberlist}
\item[\label{eq:general-position-for-cald}]
      $\bigcap_{U \in \calu_0} \dd U = \emptyset$
      whenever
      $\calu_0 \subseteq G_{C_{\lambda}}\cdot
       \calu(\cald_\lambda)$  fulfills
      $\left|\calu_0\right| > m = k_G \cdot (d_X + 1)$;
\item[\label{counting_cald_{lambda}_in_terms_of_cald_{lambda}}]
      If $D$, $D' \in G_{C_\lambda} \cdot \cald_\lambda$
      and $U_D = U_{D'}$ then $D = D'$.
\end{numberlist}
\end{lemma}

The proof of Lemma~\ref{lem:additional_properties_of_cald}
will use the following lemma, that we will prove first.

\begin{lemma}
\label{lem:can-shrink-and-still-cover}
There exists collections
$\{V_D \mid D \in G_{C_{\lambda}}\cdot \cald_{\lambda}\}$ and
$\{W_D \mid D \in G_{C_{\lambda}}\cdot \cald_{\lambda}\}$
of open subsets of $S_{C_{\lambda}}$ satisfying:
\begin{enumerate}
\item
\label{lem:can-shrink-and-still-cover:(1)}
For $D \in G_{C_{\lambda}}\cdot \cald_{\lambda}$
the sets $W_D$ $,\overline{W_D}$, $V_D$ and $\overline{V_D}$ are
$G_D$-invariant subsets of the $G_D$-space $U_D$;
\item
\label{lem:can-shrink-and-still-cover:(2)}
We have $W_D \subseteq \overline{W_D} \subseteq V_D \subseteq
\overline{V_D}
\subseteq U_D$ for $D \in G_{C_{\lambda}}\cdot \cald_{\lambda}$;
\item
\label{lem:can-shrink-and-still-cover:(2.5)}
We have $W_D = (\overline{W_D})^{\circ}$ and $V_D =
(\overline{V_D})^{\circ}$
 for $D \in G_{C_{\lambda}}\cdot \cald_{\lambda}$;
\item
\label{lem:can-shrink-and-still-cover:(3)}
$W_{gD} = gW_D$ and $V_{gD} = gV_D$ holds for
$g \in G_{C_{\lambda}}$
and $D \in G_{C_{\lambda}}\cdot \cald_{\lambda}$;
\item
\label{lem:can-shrink-and-still-cover:(4)}
$K_{\Xi} \subseteq \bigcup_
    {\substack{D \in \cald\\D \not\in G\cdot \cald_{\lambda}}}
    \Phi_{(-\epsilon,\epsilon)}  (D^{\circ})
    ~ \cup ~ \bigcup_{g \in G} \bigcup_{D \in \cald_{\lambda}}
    g \cdot \Phi_{(-\epsilon,\epsilon)}
    \left(D(\sigma_D(\overline{W_D}))^{\circ}\right)$;

\item
\label{lem:can-shrink-and-still-cover:(5)}
If $V_D = V_E$ for $D,E \in  G_{C_{\lambda}}\cdot \cald_{\lambda}$,
then $D = E$;

\item
\label{lem:can-shrink-and-still-cover:(6)}
If for $D,E \in  G_{C_{\lambda}}\cdot \cald_{\lambda}$
the intersection $V_D \cap V_E$ contains both
$W_D$ and $W_E$, then $V_D = V_E$.

\end{enumerate}
\end{lemma}

\begin{proof}
Choose a metric $d$ on $S_{C_{\lambda}}$
which is $G_{C_{\lambda}}$-invariant.
Consider $D \in G_{C_{\lambda}}\cdot \cald_{\lambda}$.
Put
$$V_D(n) = \left(\overline{U_D^{-1/n}}\right)^{\circ}.$$
Notice for the sequel that for an open subset $Y$ of a
topological space we have
$\overline{Y} = \overline{\overline{Y}^{\circ}}$ and
$Y \subseteq \overline{Y}^{\circ}$ but in general
$Y \not=\overline{Y}^{\circ}$ .
Hence $V_D(n)$ is a $G_D$-invariant open subset of $S_{C_{\lambda}}$
with $V_D(n) = \left(\overline{V_D(n)}\right)^{\circ}$.
We get for $D \in G_{C_{\lambda}}\cdot \cald_{\lambda}$
and $g \in G_{C_{\lambda}}$
that $\overline{V_D(n)} \subseteq  V_D(n+1) \subseteq U_D$,
$U_D = \bigcup_{n \ge 1} V_D(n)$ and $gV_D(n) = V_{gD}(n)$ holds.
Denote for $D \in G_{C_{\lambda}}\cdot \cald_{\lambda}$
and $n \ge 1$ the restriction by
$$D_n ~ = ~ D\left(\sigma_D(\overline{V_D(n)}),l_D - 1/n\right).$$
Lemma~\ref{lem:properties_of_U_E_and_sigma_E(T)}
~\ref{lem:properties_of_U_E_and_sigma_E(T):identifyingE(T)^circ}
and
Lemma~\ref{lem:basics_about_boxes}
~\ref{lem:basics_about_boxes:shape_of_B_in_terms_of_central_slice}
imply $D^{\circ} = \bigcup_{n \ge 1} D_n^{\circ}$ and
$D_n \subseteq D_{n+1}^{\circ}$ for $n \ge 1$.

Put
\begin{eqnarray*}
K_{\Xi}' & := &  K_{\Xi} - \left(K_{\Xi} \cap
  \bigcup_{\substack{D \in \cald\\D \not \in G \cdot \cald_{\lambda}}}
  \Phi_{(-\epsilon,\epsilon)} (D^{\circ})\right).
\end{eqnarray*}
Since
$K_{\Xi} \subseteq \bigcup_{D \in \cald}
        \Phi_{(-\epsilon,\epsilon)} (D^{\circ})$ by assumption,
we have
\begin{eqnarray}
K_{\Xi}' & \subseteq & \bigcup_{D \in G \cdot \cald_{\lambda}}
\Phi_{(-\epsilon,\epsilon)} (D^{\circ});
\label{lem:can-shrink-and-still-cover:(a)}
\\
K_{\Xi} & \subseteq & K_{\Xi}' \cup
    \bigcup_{\substack{D \in \cald\\D \not
    \in G \cdot \cald_{\lambda}}}
\Phi_{(-\epsilon,\epsilon)} (D^{\circ}).
\label{lem:can-shrink-and-still-cover:(b)}
\end{eqnarray}
Since each $D \in \cald$ is not huge,
$\Phi_{[-\epsilon,\epsilon]}(D)$ is contained in
$C_{\lambda}^{\circ}$ for $D \in \cald_{\lambda}$.
Since
$gC_{\lambda} \cap C_{\lambda} \not=
          \emptyset \Rightarrow g \in G_{C_{\lambda}}$,
we get from~\eqref{lem:can-shrink-and-still-cover:(a)}
\begin{eqnarray}
K_{\Xi}' \cap C_{\lambda} & \subseteq &
  \bigcup_{D \in G_{C_{\lambda}} \cdot \cald_{\lambda}}
  \Phi_{(-\epsilon,\epsilon)} (D^{\circ});
  \label{lem:can-shrink-and-still-cover:(d)}
\\
K_{\Xi}' & = & \bigcup_{g \in G} g \cdot \left(K_{\Xi}' \cap
  C_{\lambda}\right).
\label{lem:can-shrink-and-still-cover:(c)}
\end{eqnarray}
Since $K_{\Xi}'$ is closed, $K_{\Xi}' \cap C_{\lambda}$ is compact.
Because $D^\circ = \bigcup_{n \geq 1} D_n^\circ$ and
$D_n^\circ \subset D_{n+1}^\circ$,
\eqref{lem:can-shrink-and-still-cover:(d)}
implies that there exists  a natural number $N$ with
$$K_{\Xi}' \cap C_{\lambda}  \subseteq \bigcup_{D \in G_{C_{\lambda}}
  \cdot \cald_{\lambda}}
\Phi_{(-\epsilon,\epsilon)} (D_N^{\circ}).$$
Since $D_N \subseteq  D\left(\sigma_D(\overline{V_D(N)}\right)$
we conclude
$$K_{\Xi}' \cap C_{\lambda}  ~ \subseteq ~
      \bigcup_{D \in G_{C_{\lambda}} \cdot
      \cald_{\lambda}} \Phi_{(-\epsilon,\epsilon)}
      \left(D\left(\sigma_D(\overline{V_D(N)}\right)
                                      ^{\circ}\right).$$
We conclude from~\eqref{lem:can-shrink-and-still-cover:(c)}
$$K_{\Xi}'  ~ \subseteq ~ \bigcup_{g \in G} \bigcup_{D \in
  \cald_{\lambda}} g \cdot \Phi_{(-\epsilon,\epsilon)}
\left(D\left(\sigma_D(\overline{V_D(N)}\right)^{\circ}\right).$$
We conclude from~\eqref{lem:can-shrink-and-still-cover:(b)}
\begin{eqnarray*}
K_{\Xi} & \subseteq &
\bigcup_{\substack{D \in \cald\\D \not\in G\cdot \cald_{\lambda}}}
  \Phi_{(-\epsilon,\epsilon)}  (D^{\circ})
  ~ \cup ~ \bigcup_{g \in G} \bigcup_{D \in \cald_{\lambda}}
  g \cdot \Phi_{(-\epsilon,\epsilon)}
  \left(D(\sigma_D(\overline{V_D(N)}))^{\circ}\right).
\end{eqnarray*}
So for every choice of natural numbers
$\{n_D \mid D \in G_{C_{\lambda}}\cdot \cald_{\lambda}\}$
satisfying $n_D \ge N$ and $n_{gD} = n_D$ for
$g \in G_{C_{\lambda}}$ and
$D \in G_{C_{\lambda}}\cdot \cald_{\lambda}$ the collection
$\{V_D(n_D) \mid D \in G_{C_{\lambda}}\cdot \cald_{\lambda}\}$
has the properties
\begin{itemize}

\item $V_D(n_D)$ is a $G_D$-invariant subset of
      the $G_D$-space $V_D$;

\item $V_D(n_D) = \left(\overline{V_D(n_D)}\right)^{\circ}$;

\item We have $V_D(n_D) \subseteq \overline{V_D(n_D)} \subseteq U_D$
      for $D \in G_{C_{\lambda}}\cdot \cald_{\lambda}$;

\item $V_{gD}(n_{gD}) = gV_D(n_D)$ holds for $g \in G_{C_{\lambda}}$
      and $D \in G_{C_{\lambda}}\cdot \cald_{\lambda}$;

\item $K_{\Xi} ~ \subseteq ~
        \bigcup_{\substack{D \in \cald\\D
               \not\in G\cdot \cald_{\lambda}}}
        \Phi_{(-\epsilon,\epsilon)} (D^{\circ}) ~ \cup ~
        \bigcup_{g \in G} \bigcup_{D \in \cald_{\lambda}}
        g \cdot \Phi_{(-\epsilon,\epsilon)}
        \left(D(\sigma_D(\overline{V_D(n_D)}))^{\circ}\right)$.

\end{itemize}

Next we show that for some choice of numbers
$\{n_D \mid D \in G_{C_{\lambda}}\cdot \cald_{\lambda}\}$
satisfying $n_D \ge N$ and $n_{gD} = n_D$ for
$g \in G_{C_{\lambda}}$ and
$D \in G_{C_{\lambda}}\cdot \cald_{\lambda}$ the collection
$\{V_D(n_D) \mid D \in G_{C_{\lambda}}\cdot \cald_{\lambda}\}$
also has property~\ref{lem:can-shrink-and-still-cover:(5)}.
We can write $G_{C_{\lambda}}\cdot \cald_{\lambda}$
as the disjoint union
$$G_{C_{\lambda}}\cdot \cald_{\lambda} =
     \calc_1 \amalg \calc_2 \amalg \ldots \amalg \calc_r$$
of its $G_{C_{\lambda}}$-orbits.
We show by induction that we can find numbers
$n_1$, $n_2$, $\dots$, $n_r$ with $n_k \geq N$
such that if we set $n_D = n_k$ for
$D \in \calc_k$ the collection
$\{ V_D = V_D( n_D) \; | \; D \in
           \calc_1 \cup \calc_2 \cup \dots \cup \calc_k \}$
satisfies property~\ref{lem:can-shrink-and-still-cover:(5)}.
The induction beginning $k = 1$ is trivial,
the induction step from $k-1$ to $k$ is done as follows.
For given $D \in \calc_k$ choose $n_0$ with
$\overline{V_D ( n_0 )} \neq \emptyset$.
Since $S_{C_{\lambda}}$ is connected, the
non-empty closed subset $\overline{V_D ( n )}$ and the open
subset $V_D( n +1) \neq S_{C_\lambda}$ cannot agree for $n \geq n_0$.
In particular $V_D( n ) \subset V_D( n+1 )$ for all
$n \geq n_0$ and since
$\calc_1 \cup \calc_2 \cup \dots \cup \calc_{k-1}$ is a finite
set we can find a number $n_k$ such that
$V_D( n_k ) \neq V_E ( n_E)$ for
$E \in \calc_1 \cup \calc_2 \cup \dots \cup \calc_{k-1}$.
By invariance under the
$G_{C_{\lambda}}$-action the same statement holds for all
$D \in \calc_k$ with this $n_k$.
If we have $D$, $gD \in \calc_k$ with
$g \in G_{C_{\lambda}}$ then since $\mu_{C_{\lambda}}$
from
Lemma~\ref{lem:basics_about_boxes}~
\ref{lem:basics_about_boxes:shape_of_B_in_terms_of_central_slice}
is $G_{C_{\lambda}}$-invariant and the action on the interval
trivial and $D$ is a $\Fin$-set, $U_D = g U_D$
already implies $D=gD$.
Therefore property~\ref{lem:can-shrink-and-still-cover:(5)}
holds for
$\calc_1 \cup \calc_2 \cup \dots \cup \calc_k$.
We have therefore verified
property~\ref{lem:can-shrink-and-still-cover:(5)}
for the collection
$V_D = V_D ( n_D)$ with $D \in G_{C_{\lambda}} \cald_{\lambda}$.

In order to achieve
property~\ref{lem:can-shrink-and-still-cover:(6)}
we repeat this construction replacing the collection of boxes
$G_{C_{\lambda}} \cald_{\lambda}$ with the collection of boxes
$\left\{D\left(\sigma_D(\overline{V_D(n_0)})\right)
  \mid D \in G_{C_{\lambda}}\cdot \cald_{\lambda}\right\}$.
Namely, put
$$W_D(n) ~ = ~ \left(\overline{V_D^{-1/n}}\right)^{\circ}.$$
Lemma~\ref{lem:properties_of_U_E_and_sigma_E(T)}
~\ref{lem:properties_of_U_E_and_sigma_E(T):identifyingE(T)^circ}
implies $U_{D\left(\sigma_D(\overline{V_D(n)}\right)} = V_D(n)$.
Thus we get open subsets $W_D(n)  \subseteq V_D$
and a natural number $N'$
such that for every choice of natural numbers
$\{n_D \mid D \in G_{C_{\lambda}}\cdot \cald_{\lambda}\}$
satisfying $n_D \ge N'$ and $n_{gD} = n_D$ for $g \in G_{C_\lambda}$
and $D \in G_{C_{\lambda}}\cdot \cald_{\lambda}$ the collection
$\{W_D(n_D) \mid D \in
G_{C_{\lambda}}\cdot \cald_{\lambda}\}$ has the properties
\begin{itemize}

\item $W_D(n_D)$ is a $G_D$-invariant subset of the
      $G_D$-space $V_D$;

\item $W_D(n_D) = \left(\overline{W_D(n_D)}\right)^{\circ}$;

\item We have $W_D(n_D) \subseteq \overline{W_D(n_D)} \subseteq V_D$
      for $D \in G_{C_{\lambda}}\cdot \cald_{\lambda}$;

\item $W_{gD}(n_{gD}) = gW_D(n_D)$ holds for $g \in G_{\lambda}$
      and $D \in G_{C_{\lambda}}\cdot \cald_{\lambda}$;

\item $K_{\Xi} ~ \subseteq ~
       \bigcup_{\substack{D \in \cald\\D \not\in G\cdot
           \cald_{\lambda}}}
       \Phi_{(-\epsilon,\epsilon)}  (D^{\circ})
       ~ \cup ~ \bigcup_{g \in G} \bigcup_{D \in \cald_{\lambda}}
       g \cdot \Phi_{(-\epsilon,\epsilon)}
       \left(D\left(\sigma_D(\overline{W_D(n_D)})\right)
                                             ^{\circ}\right)$.

\end{itemize}
Consider $D,E \in \calc$ with $V_D \not= V_E$.
Since $V_D = \bigcup_{n \ge 1} W_D(n)$ and
$V_E = \bigcup_{n \ge 1} W_E(n)$ holds,
we can find $N'(D,E)$ such that
$V_D \cap V_E$ does not contain both $W_D(n)$
and $W_E(n)$ for $n \ge N'(D,E)$.
Define $N''$ to be the maximum
over the numbers $N'(D,E)$ for
$D,E \in G_{C_{\lambda}}\cdot \cald_{\lambda} $
with $V_D \not= V_E$ and $N'$.
Put $W_D = W_D(N'')$ for
$D \in G_{C_{\lambda}}\cdot \cald_{\lambda}$.
Then  the collections
$\{V_D \mid D \in G_{C_{\lambda}}\cdot \cald_{\lambda}\}$ and
$\{W_D \mid D \in G_{C_{\lambda}}\cdot \cald_{\lambda}\}$
have all the desired properties.
This finishes the proof of
Lemma~\ref{lem:can-shrink-and-still-cover}.
\end{proof}

Now we can prove Lemma~\ref{lem:additional_properties_of_cald}.

\begin{proof}
In the sequel we will use the collections
$\{V_D \mid D \in G_{C_{\lambda}}\cdot \cald_{\lambda}\}$ and
$\{W_D \mid D \in G_{C_{\lambda}}\cdot \cald_{\lambda}\}$
appearing in
Lemma~\ref{lem:can-shrink-and-still-cover}.
We apply Proposition~\ref{pro:general-position-for-subsets}
in the case, where the space $Z$ is $S_{C_{\lambda}}$,
the finite group $F$ is $G_{C_\lambda}$,
the collection $\calu$ is
$\{V_D \mid D \in G_{C_{\lambda}}\cdot \cald_\lambda\}$
and we use the subsets $W_D \subseteq V_D$.
Since $S_{C_{\lambda}} \subseteq X$ is closed,
we have $\dim(S_{C_{\lambda}}) \le d_X$.
So from Proposition~\ref{pro:general-position-for-subsets}
we obtain for every $D \in G_{C_\lambda}\cdot \cald_{\lambda}$
an open subset $V_D'' \subseteq S_{C_{\lambda}}$
such that the following holds
\begin{itemize}
\item $W_D \subseteq V_D'' \subseteq V_D$;
\item If $\calu_0 \subseteq \{V_D'' \mid D \in G_{C_\lambda}
                                         \cdot \cald_\lambda\}$
      has more than $m = k_G \cdot (d_X + 1)$ elements, then
      \begin{equation*}
      \bigcap_{U'' \in \calu_0} \dd U'' = \emptyset;
      \end{equation*}
\item $(V_{gD})'' = (gV_D)'' = g(V_D'')$ for
      $g \in G_{C_{\lambda}}$ and
      $D \in G_{C_\lambda}\cdot \cald_\lambda$.
\end{itemize}

Now define
$$V_D' = \overline{V_D''}^{\circ}.$$
Since $V_D''$ is open, we conclude
$V_D' = \overline{V_D'}^{\circ}$.
Recall that $V_D = \overline{V_D}^{\circ}$ and
$W_D = \overline{W_D}^{\circ}$.
Notice that $V_D''$ is not necessarily $V_D'$.
Since $V_D''$ is open, we get
$V_D'' \subseteq V_D'$ and hence
$V_D'' \cap \partial V_D' = \emptyset$.
We have
$$\partial V_D' \subseteq \overline{V_D'}
    = \overline{\overline{V_D''}^{\circ}} =
    \overline{V_D''}.$$
Hence $\partial V_D' \subseteq \partial V_D''$.
Thus we have constructed for every
$D \in G_{C_\lambda}\cdot \cald_{\lambda}$ an open subset
$V_D' \subseteq S_{C_{\lambda}}$
such that the following holds
\begin{itemize}
\item
      $W_D \subseteq V_D' \subseteq V_D$;
\item If $\calv_0 \subseteq \{V_D' \mid D \in
           G_{C_\lambda}\cdot \calu(\cald_\lambda)\}$
      has more than $m = k_G \cdot (d_X + 1)$ elements, then
      \begin{equation*}
      \bigcap_{V' \in \calu_0} \dd V' = \emptyset;
      \end{equation*}
\item $(V_{gD})' = (gV_D)' = g(V_D')$ for $g \in G_{C_{\lambda}}$
      and $D \in G_{C_\lambda}\cdot \calu(\cald_\lambda)$;
\item $\overline{V_D'}^{\circ} = V_D'$.
\end{itemize}

We use next restriction of boxes to diminish some of the boxes in
$\cald$ as follows. Consider $D \in \cald$.
Suppose there exists
$g_0 \in G$ and $D_0 \in G_{C_{\lambda}} \cdot\cald_{\lambda}$
with $D = g_0 D_0$.
Then define
$\widehat{D} =
 g_0D_0\left(\sigma_{D_0}(\overline{V_{D_0}'})\right)$.
We have to check that this is well-defined.
Suppose we have $g_i \in G$ and
$D_i \in G_{C_{\lambda}} \cdot \cald_{\lambda}$
with $D = g_i D_i$ for $i = 0,1$.
We have $D_1 = g_1^{-1}g_0D_0$.
This implies $g_1^{-1}g_0 \in G_{C_{\lambda}}$
(see Lemma~\ref{lem:properties_of_U_E_and_sigma_E(T)}
~\ref{lem:properties_of_U_E_and_sigma_E(T):G_E_subsetG_{C_lambda}})
and hence $g_1^{-1}g_0 U_{D_0} = U_{D_1}$.
We conclude $g_1^{-1}g_0 (V_{D_0})' = (V_{D_1})'$ and hence
$g_1^{-1}g_0 D_0\left(\sigma_{D_0}(\overline{V_{D_0}'})\right) =
D_1\left(\sigma_{D_1}(\overline{V_{D_1}'})\right)$.
If there does not exist $g_0 \in G$ and $D_0 \in \cald_{\lambda}$
with $D = g_0 D_0$, we put $\widehat{D} = D$.
Define a new collection of boxes
$$\widehat{\cald} = \{\widehat{D} \mid D \in \cald\}.$$

Next we want to show that $\widehat{\cald}$ is a
$\delta_r$-good box cover of
$K_{\Xi} = \bigcup_{\xi \in \Xi} A_\xi$.
Since $W_D \subseteq V_D'$
for $D \in \cald_{\lambda}$, we conclude from
property~\ref{lem:can-shrink-and-still-cover:(4)} appearing in
Lemma~\ref{lem:can-shrink-and-still-cover}
$$K_{\Xi} ~ \subseteq ~
     \bigcup_{\widehat{D} \in \widehat{\cald}}
        \Phi_{(-\epsilon,\epsilon)}(\widehat{D}^{\circ}).$$
One easily checks that the other required properties of a
$\delta_r$-good box cover do pass from $\cald$ to $\widehat{\cald}$
since elements in $\widehat{\cald}$ are obtained from
those in $\cald$ by restriction in a $G$-equivariant way.

We conclude
$G_{C_{\lambda}}\cdot \calu(\widehat{\cald}_{\lambda}) ~ \subseteq ~
            \{V_D' \mid D \in G_{C_{\lambda}}\cdot \cald_{\lambda}\}$
from Lemma~\ref{lem:properties_of_U_E_and_sigma_E(T)}
~\ref{lem:properties_of_U_E_and_sigma_E(T):identifyingE(T)^circ}.
(We remark that $D \cap B_\lambda \neq \emptyset$ does not necessarily
imply $\widehat{D} \cap B_\lambda \neq \emptyset$ and the above
inclusion may be a strict inclusion.)
By construction $\widehat{\cald}$
satisfies~\eqref{eq:general-position-for-cald}.

Suppose for
$\widehat{D}, \widehat{D'}\in G_{C_{\lambda}} \cdot
                             \widehat{\cald}_{\lambda}$
that $U_{\widehat{D}} = U_{\widehat{D'}}$.
By construction $U_{\widehat{D}} = V_D'$ and
$U_{\widehat{D}} = V_{D'}'$.
We conclude $V_D' = V_{D'}'$ and hence both $W_D$ and $W_{D'}$ are
contained in $V_D \cap V_{D'}$.
Properties property~\ref{lem:can-shrink-and-still-cover:(5)}
and \ref{lem:can-shrink-and-still-cover:(6)} appearing in
Lemma~\ref{lem:can-shrink-and-still-cover}
imply $D = D'$.
We conclude that $\widehat{\cald}$
satisfies~\eqref{counting_cald_{lambda}_in_terms_of_cald_{lambda}}.
This finishes the proof of
Lemma~\ref{lem:additional_properties_of_cald}.
\end{proof}

Now we have finished our arrangement of $\cald$ and can now
construct the desired new collection $\cald'$ out of $\cald$
as demanded in Lemma~\ref{lem:induction_step_for_prop}.


\subsection{Carrying out the induction step}

Recall that we defined numbers $m$, $M$, $a$, $b$ and $c$
in the beginning of
Section~\ref{subsec:preliminaries-and-basic-induction-structure}.
Recall also that $N$ is the number of $G$-orbits of $\Lambda$
and that $\e$ is given in
Proposition~\ref{prop:cover-K-by-boxes-hyp}.
In the sequel we will abbreviate
\begin{eqnarray*}
a_{\pm} & = & \pm l_{A_{\lambda}}/2;
\\
b_{\pm} & = & \pm l_{B_{\lambda}}/2;
\\
c_{\pm} & = & \pm l_{C_{\lambda}}/2.
\end{eqnarray*}
We have
\begin{eqnarray*}
a_+ - a_-  & = & a;
\\
b_+ -a_+ & = & a_- - b_- = b;
\\
c_+ -b_+ & = & b_- - c_- = c.
\end{eqnarray*}
Put
\begin{equation}
\mu :=  \frac{\epsilon}{(N+1)(m+1)}
\quad \mbox{and} \quad \eta := \frac{\mu}{5}.
\label{def_of_mu}
\end{equation}

We will use the $G_{C_{\lambda}}$-homeomorphism
\begin{eqnarray}
\label{identification_of_C_lambda}
& & \mu_{C_{\lambda}} \colon
  S_{C_{\lambda}} \times [c_-,c_+] \to C_{\lambda},
   \quad (x,\tau) ~ \mapsto ~ \Phi_{\tau}(x)
\end{eqnarray}
from Lemma~\ref{lem:basics_about_boxes}
~\ref{lem:basics_about_boxes:shape_of_B_in_terms_of_central_slice}
as an identification.
Note that $A_\lambda = S_{A_\lambda} \times [a_-,a_+]$,
$B_\lambda = S_{B_\lambda} \times [b_-,b_+]$ and
$B_\lambda^{\circ} = \left(S_{B_\lambda} \cap B^{\circ}\right)
                                               \times (b_-,b_+)$
under this identification.
Note that for $g \in G_{C_\lambda}$ we have
$g \cdot (x,t) = (g x,t)$,
by Lemma~\ref{lem:basics_about_boxes}~
\ref{lem:basics_about_boxes:shape_of_B_in_terms_of_central_slice}.

Choose a $G_{C_{\lambda}}$-invariant metric
$d_{S_{C_{\lambda}}}$ on $S_{C_{\lambda}}$.
Consider $D \in \cald_{\lambda}$.
Recall that $D \subseteq C_{\lambda}$ and that the retraction
$\pi_{C_{\lambda}} \colon C_{\lambda} \to S_{C_{\lambda}}$
induces a $G_D$-homeomorphism
$S_D \to \pi_{C_{\lambda}}(S_D)$.
We have introduced the continuous $G_D$-invariant map
$$\tau_D \colon \pi_{C_{\lambda}}(S_D) \to [c_-,c_+]$$
in Definition~\ref{def:U_E}.
It is uniquely characterized by
$\mu_{C_{\lambda}}(x) =
     (\pi_{C_{\lambda}}(x),\tau_D \circ\pi_{C_{\lambda}}(x))$
for $x \in S_D$.
Since $\tau_D$ is continuous and $\pi_{C_{\lambda}}(S_D)$
is compact, we can find $\delta_D > 0$ such that
\begin{eqnarray}
\label{tau_D,uniform_cont.} \hspace{10mm}
|\tau_D(x) - \tau_D(y)| & ~ \le & \eta
  \quad \text { for } x,y \in \pi_{C_{\lambda}}(S_D)
  \text { with } d_{S_{C_{\lambda}}}(x,y) < \delta_D.
\end{eqnarray}
Because $\cald_{\lambda}$ is finite we can set
\begin{eqnarray}
& \delta   ~ :=  ~
   \min\left\{\delta_D \mid D \in \cald_{\lambda}\right\} .&
\label{lem:can-shrink-and-still-cover:def_of_delta}
\end{eqnarray}
Then $\delta > 0$.

In the sequel interiors and closures of subsets of
$S_{C_{\lambda}}$ are to be understood
with respect to $S_{C_{\lambda}}$.
One easily checks with this convention that
Lemma~\ref{lem:basics_about_boxes}
~\ref{lem:basics_about_boxes:shape_of_B_in_terms_of_central_slice}
implies $S_{B_{\lambda}}^{\circ} = B^{\circ} \cap S_{B_{\lambda}}$
since $B_{\lambda} \subseteq C_{\lambda}^{\circ}$.

Next we want to apply Proposition~\ref{prop:cover-Z-compatitible-with-calu}
to the locally connected compact metric space
$Z := S_{C_{\lambda}}$ with the obvious isometric
$F := G_{B_{\lambda}}$-action (note that $F \subseteq G_{C_\lambda}$
by Lemma~\ref{lem:properties_of_U_E_and_sigma_E(T)}
~\ref{lem:properties_of_U_E_and_sigma_E(T):G_E_subsetG_{C_lambda}}),
the  $F = G_{B_{\lambda}}$-invariant open subset
$Y := S_{C_{\lambda}} \cap C_{\lambda}^{\circ}$
which is locally connected by
Lemma~\ref{lem:basics_about_boxes}
~\ref{lem:basics_about_boxes:S_BcapB^circ_locally_connected},
the collection of sets
$\calu := \calu(\cald_{\lambda})$
and the number $\delta > 0$ in
\eqref{lem:can-shrink-and-still-cover:def_of_delta}.
Note that for $D \in \cald_\lambda$ we have by definition
$D \cap B_\lambda \neq \emptyset$
and therefore $D \subseteq C_\lambda^\circ$ since $D$ is not huge.
Therefore $\overline{U_D} \subseteq \pi_{C_\lambda}(D)
                \subseteq S_{C_\lambda} \cap C_\lambda^\circ = Y$
for $D \in \cald_\lambda$.
In the notation just introduced this means $\overline{U} \subset Y$
for $U \in \calu$.
Thus we can indeed apply
Proposition~\ref{prop:cover-Z-compatitible-with-calu}
in this situation.
By intersecting the resulting open covering with
$S_{B_{\lambda}}^{\circ} = S_{B_{\lambda}} \cap B^{\circ}$,
we obtain a collection
$\calv = \calv^0 \cup \calv^1 \cup \dots \cup \calv^m$
of open subsets of $S_{B_{\lambda}}^{\circ}$
which has the following properties:
\begin{numberlist}
\item[\label{eq:calv-covers-U}]
           $\calv$ is a  open covering of $S_{B_{\lambda}}^{\circ}$
           consisting of finitely many elements;
\item[\label{eq:calv-intersects-few-U}]
           For every $V \in \calv$ there are at most
           $k^2_G \cdot (d_X + 1)$
           different $U \in \calu(\cald_{\lambda})$
           such that $U \cap V \neq \emptyset$ and
           $V \not\subset U$;

\item[\label{eq:few-intersections-among-the-calv}]
           For fixed $j$ and $V_0 \in \calv^j$ we have
           $V_0 \cap V \neq \emptyset$ for at most
           $2^{j+1}-2 < 2^{m+1}$ different
           subsets $V \in \calv^0 \cup \dots \cup \calv^{j-1}$;

\item[\label{eq:calv-j-mutually-disjoint}]
           For fixed $j$ and $V_0,V_1 \in \calv^j$ we have
           either $V_0 = V_1$ or
           $\overline{V_0} \cap  \overline{V_1} = \emptyset$;

\item[\label{eq:calv-invariant}]
           Each $\calv^i$ is $G_{B_\lambda}$-invariant,
           i.e., $gV \in \calv^i$ if $g \in G_{B_\lambda}$,
           $V \in \calv^i$;

\item[\label{eq:calv-Fin-subset}]
           For $V \in \calv$ its closure $\overline{V}$ is a
           $\Fin$-subset of $S_{B_{\lambda}}$
           with respect
           to the $G_{B_{\lambda}}$-action;

\item[\label{eq:calv-overline{V}^circ=V}]
           We have $\overline{V}^{\circ} = V$ for $V \in  \calv$;

\item[\label{eq:calv-diameter}] For every
           $V \in \calv$ the diameter
           of $V$ is bounded by $\delta$;

\item[\label{eq:calv-has-no-doubles}]
           $\calv^j \cap \calv^k = \emptyset$ if $j \neq k$.

\end{numberlist}
Properties~\eqref{eq:calv-covers-U},
\eqref{eq:calv-intersects-few-U},
\eqref{eq:few-intersections-among-the-calv},
\eqref{eq:calv-j-mutually-disjoint},
\eqref{eq:calv-invariant} and~\eqref{eq:calv-diameter}
are direct consequences of
Proposition~\ref{prop:cover-Z-compatitible-with-calu}.
Property~\eqref{eq:calv-Fin-subset} follows from
properties~\eqref{eq:calv-j-mutually-disjoint}
and~\eqref{eq:calv-invariant}.
Property~\eqref{eq:calv-has-no-doubles} can be achieved by
replacing $\calv^j$ by a subset of $\calv^j$ if necessary.

Since for every subset $Y
\subseteq S_{C_{\lambda}}$ with $\overline{Y}^{\circ} = Y$ we have
$$\left(\overline{Y \cap S_{B_{\lambda}}^{\circ}}\right)^{\circ}
  \subseteq \left(\overline{Y} \cap S_{B_{\lambda}}\right)^{\circ}
   = \overline{Y}^{\circ} \cap S_{B_{\lambda}}^{\circ}
   = Y\cap S_{B_{\lambda}}^{\circ} \subseteq \left(\overline{Y \cap
       S_{B_{\lambda}}^{\circ}}\right)^{\circ}$$
and hence
$\left(\overline{Y \cap S_{B_{\lambda}}^{\circ}}\right)^{\circ}
 = Y \cap S_{B_{\lambda}}^{\circ}$,
property~\eqref{eq:calv-overline{V}^circ=V} holds.
We mention that because of \eqref{eq:calv-Fin-subset}
we can consider for $V \in \calv$ the restriction
$B_{\lambda}(\overline{V})$ and
property~\eqref{eq:calv-overline{V}^circ=V} ensures
$$
\begin{array}{lclcl}
S_{B_{\lambda}(\overline{V})} &  = & \overline{V}; & &
\\
S_{B_{\lambda}(\overline{V})}^{\circ}  & = &
S_{B_{\lambda}(\overline{V})} \cap B_{\lambda}^{\circ} & = & V.
\end{array}
$$

The collection $\cald'$ we are seeking will be of the form
$\cald \cup \{ g B_\lambda( \overline{W};a^W_-,a^W_+ ) \; | \;
                      W \in \calw, g \in G \}$,
where $\calw \subseteq \calv$.
In order to find suitable $\calw \subseteq \calv$ and $a^W_\pm$
we proceed by a subinduction over $n = -1,\dots,m$.
Using \eqref{eq:delta-r} and \eqref{def_of_mu} we set
$$\delta_{r,n} ~ := ~ \delta_r - (n+1) \cdot \mu
            \quad \text {for } n = -1,0,1, \ldots, m.$$
Clearly,
\begin{equation*}
\delta_r = \delta_{r,-1} > \delta_{r,0} > \dots >
  \delta_{r,n-1} > \delta_{r,n} > \dots >
  \delta_{r,m} = \delta_{r+1}
\end{equation*}
and
$\delta_{r,n-1} - \delta_{r,n} =  \mu$.
For $j = 0,\dots,m$ let
$$
K_{\lambda}^{(j)} :=
  \bigcup_{V \in \calv_0 \cup \dots \cup \calv_j} V \times [a_-,a_+].
$$
(Recall that we use \eqref{identification_of_C_lambda}
to identify $K_\lambda^{(j)}$ with a subset of $C_\lambda$.)
The induction step from $(n-1)$ to $n$ is formulated in the
following Lemma.

\begin{lemma}[Induction step : $n-1$ to $n$]
\label{lem:subinductionstep}
Assume that we have for $j=0,\dots,n-1$ subsets
$\calw^j \subseteq \calv^j$ and numbers
$\{a^{W}_\pm \mid W \in \calw^j\}$
satisfying $b_-  \le a_-^W < a_+^W\le b_+$
such that the collection of boxes
$\cald^{n-1} = \cald \cup \{ gD^{W} \; | \;
  W \in \calw^0 \cup \dots \cup \calw^{n-1}, g \in G \}$
is a $\delta_{r,n-1}$-good box cover of
$K_\Xi \cup G K_\lambda^{(n-1)}$,
where $D^W := B_\lambda(\overline{W};a_-^W,a^W_+)$ for
$W \in \calw^0 \cup \dots \cup \calw^{n-1}$.

Then there is  a subset
$\calw^{n} \subseteq \calv^{n}$
and numbers $a^W_\pm \in \IR$
with $b_-  \le a_-^W < a_+^W \le b_+$ for $W \in \calw^{n}$
such that
$\cald^{n} = \cald^{n-1} \cup \{g D^W \; | \;
          W \in \calw^{n}, g \in G \}$
is a $\delta_{r,n}$-good box cover of
$K_\Xi \cup G K_\lambda^{(n)}$,
where $D^W = B_\lambda(\overline{W};a^W_-,a^W_+)$ for
$W \in \calw^{n}$.
\end{lemma}

Since \eqref{eq:calv-covers-U} implies
$S_{A_{\lambda}} \subseteq S_{B_{\lambda}}^{\circ} =
      \bigcup_{V \in \calv} V$
we get $A_{\lambda} \subseteq \bigcup_{V \in \calv}
                                    V \times [a_-,a_+]$.
We conclude $K_{\Xi'} \subseteq K_{\Xi} \cup GK^{(m)}_{\lambda}$.
($K_{\Xi'}$ was defined in
Lemma~\ref{lem:induction_step_for_prop}.)
Hence $\cald' = \cald^{m}$ is the desired
$\delta_{r+1}$-good box cover of $K_{\Xi'}$.
Therefore Lemma~\ref{lem:induction_step_for_prop}
follows from Lemma~\ref{lem:subinductionstep}.

The proof of Lemma~\ref{lem:subinductionstep}
will occupy the remainder of
Section~\ref{sec:covering-G->gamma}.

\begin{definition}
\label{def:W-is-not-yet-covered}
Let $\calw^{n}$ be the set of all $W \in \calv^{n}$ for which
$$W \times [a_-,a_+]  ~ \not\subset ~
  \bigcup_{D \in \cald^{n-1}} \Phi_{(-\epsilon,\epsilon)}
    (D^{\circ}).$$
\end{definition}

\begin{lemma}
\label{lem:use-of-friction-control}
Let $V \in \calv$ and let $\partial_{\pm} D$ be the top or bottom
of a box $D \in \cald_{\lambda}$.
Consider $t \in (c_- ,c_+ )$.
Suppose that $\partial_{\pm} D \cap (V \x \{ t \}) \neq \emptyset$.
Then
$$\partial_{\pm} D \cap (V \times [c_-,c_+])
            \subseteq V \times (t - \eta,t + \eta).$$
\end{lemma}

\begin{proof}
Consider $v \in V$ with $(v,t) \in \partial_{\pm}D$.
Then $(v,t\mp l_D/2) \in S_D$.
Hence $\tau_D(v) = t \mp l_D/2$,
where $\tau_D$ is the function introduced in
Definition~\ref{def:U_E}.
Consider $w \in V$ and $s \in [c_-,c_+]$
with $(w,s) \in \partial_{\pm} D$.
Then $s \mp l_D/2= \tau_D(w)$.
{}From~\eqref{tau_D,uniform_cont.},
\eqref{lem:can-shrink-and-still-cover:def_of_delta}
and~\eqref{eq:calv-diameter} we conclude
$|\tau_D(v)- \tau_D(w)| < \eta$ and hence $|t - s| < \eta$.
This implies
$\partial_{\pm} D \cap (V \times [c_-,c_+])
      \subseteq V \times (t- \eta,t+ \eta)$.
\end{proof}

Recall that for $D \in \cald_{\lambda}$ we have
$D \subseteq C_{\lambda}^{\circ}$ and we have
associated to such $D$ an open subset
$U_D = \pi_{C_{\lambda}}(S_D \cap D^{\circ}) =
      \pi_{C_{\lambda}}(D^{\circ})$ of $S_{C_{\lambda}}$.

\begin{definition} \label{def:cald_W_and_bad_good}
For $W \in \calw^n$ define
\begin{eqnarray*}
\cald_W & := &
\{D \in \cald^{n-1} \mid
       D^\circ \cap W \x ( {b}_-, {b}_+ ) \neq \emptyset \};
\\
\cald_W^\good & := &
\{ D \in \cald_W \mid  W \subseteq U_D \};
\\
\cald_W^{\bad} & := &
\{ D \in \cald_W \mid W \cap U_D \neq \emptyset,
                               W \not\subseteq U_D \};
\\
J_W^{\good,\pm} & := &
\{ t \in (b_-,b_+) \; | \; \exists D \in \cald_W^\good
     \; | \; \dd_\pm D^\circ \cap W \x \{ t \} \neq \emptyset \};
\\
J_W^{\bad,\pm} & := &
\{ t \in (b_-,b_+) \; | \; \exists D \in \cald_W^\bad
     \; | \; \dd_\pm D^\circ \cap W \x \{ t \} \neq \emptyset \};
\\
J_W^{\good} & := &
\{ t \in (b_-,b_+) \; | \; \exists D \in \cald_W^\good
     \; | \; D^\circ \cap W \x \{ t \} \neq \emptyset \};
\\
J_W^{\dd} & := &
J_W^{\good,-} \cup J_W^{\good,+} \cup
 J_W^{\bad,-} \cup J_W^{\bad,+}.
\end{eqnarray*}
\end{definition}

Since $D \in \cald_W$ implies $W \cap U_D \neq \emptyset$, we
have $\cald_W = \cald_W^\good \cup \cald_W^\bad$.
The reason for the names of $\cald_W^{\good}$ and $\cald_W^{\bad}$
is this.
In the construction of $D^W$ for $W \in \calw^n$
we will be able to
allow top and bottom of $D^W$ to be very close to
top or bottom of a box in $\cald_W^{\good}$
(compare Lemma~\ref{lem:distance-for-points-in-J})
but will have
to make sure that top and bottom of $D^W$ will be far away
from top and bottom of every box in $\cald_W^{\bad}$.
Thus, both for choosing $a^W_-$ and $a^W_+$ there will be two
cases: either we find a suitable top (for $a^W_-$)
respectively bottom (for $a_+^W) $ of a box
in $\cald_W^{\good}$ to put $W \x \{ a^W_\mp \}$ close by,
or all boxes from $\cald_W^{\good}$ are far away and we
will only have to worry about the boxes from $\cald_W^{\bad}$.
The crucial point will then be, that the number of members of
$\cald_W^\bad$ is uniformly bounded,
see Lemma~\ref{lem:estimate_for_|D_W^bad|}.

Note that for $g \in G_{B_\lambda}$ we have
$J_W^{\good,\pm} = J_{gW}^{\good,\pm}$,
$J_W^{\bad,\pm} = J_{gW}^{\bad,\pm}$,
$J_W^{\good} = J_{gW}^{\good}$ and
$J_W^{\dd} = J_{gW}^{\dd}$,
because $g$ acts trivially on the second factor
of $C_\lambda = S_{C_\lambda} \x [c_-,c_+]$.

\begin{lemma} \label{lem:estimate_for_|D_W^bad|}
We have
\begin{eqnarray*}
|\cald_W^{\bad}| & \le & M = (k_G)^2 \cdot (d_X + 1) + 2^{m+1}.
\end{eqnarray*}
\end{lemma}

\begin{proof}
We conclude from the definitions,
Lemma~\ref{lem:properties_of_U_E_and_sigma_E(T)}
~\ref{lem:properties_of_U_E_and_sigma_E(T):identifyingE(T)^circ}
and~\eqref{eq:calv-overline{V}^circ=V}
\begin{eqnarray*}
\cald^{n-1}  & = &
\cald \cup \{ gD^{{W}} \mid
        W \in \calw^0 \cup \dots \cup \calw^{n-1}, g \in G \};
\\
D^{W} & = & B_{\lambda}(\overline{W};a_-^{W},a_+^{W}) ~ = ~
              \overline{W} \times [a_-^{W},a_+^{W}]~
      \quad \text{ for }  W \in \calw^0 \cup \dots \cup \calw^{n},
\\
U_{D^{W}} & = & W
\quad \text{ for }  W \in \calw^0 \cup \dots \cup \calw^{n},
\end{eqnarray*}
Recall that
$gB_{\lambda} \cap B_{\lambda} \not= \emptyset
               \Rightarrow g \in G_{B_{\lambda}}$ holds.
Hence we get
\begin{eqnarray*}
\cald_W^{\bad} & = &
   \left\{D \in \cald \mid D^\circ \cap W \times (b_-,b_+)
       \not= \emptyset, W \cap U_D \not= \emptyset,
       W \not\subset U_D\right\}
\\
& & \hspace{10mm} \cup
     \left\{gD^{W'} \mid W' \in \calw^0 \cup \dots
                           \cup \calw^{n-1},   g \in G, \right.
\\ & & \hspace{20mm}
         \left. \left( gD^{W'} \right)^\circ
                 \cap W \times (b_-,b_+)  \not= \emptyset,
   W \cap U_{gD^{W'}}\not= \emptyset, W \not\subset U_{gD^{W'}}
                                                        \right\}
\\
& = & \left\{D \in \cald\mid
         D^\circ \cap W \times (b_-,b_+) \not= \emptyset,
         W \cap U_D \not= \emptyset, W \not\subset U_D\right\}
\\
& & \hspace{5mm} \cup \left\{gD^{W'} \mid
           W' \in \calw^0 \cup \dots \cup \calw^{n-1},
           g \in G_{B_\lambda},  W \cap  gW' \not= \emptyset,
                           W \not\subset gW'\right\}
\\
& \subseteq & \left\{D \in \cald_{\lambda}\mid
         W \cap U_D \not= \emptyset, W \not\subset U_D\right\}
\\
& & \hspace{5mm} \cup \left\{gD^{W'}
             \mid W' \in \calw^0 \cup \dots \cup \calw^{n-1},
             g \in G_{B_{\lambda}},  W \cap gW' \not= \emptyset,
                           W \not\subset gW'\right\}.
\end{eqnarray*}
We have $gD^{W'} = D^{gW'}$ and
$D^{W'} = D^{W''} \Leftrightarrow W' = W''$.
Hence we conclude
using~\eqref{counting_cald_{lambda}_in_terms_of_cald_{lambda}},
\eqref{eq:calv-intersects-few-U},
\eqref{eq:few-intersections-among-the-calv}
and~\eqref{eq:calv-has-no-doubles}.
\begin{eqnarray*}
|\cald_W^{\bad}| & \le & \left|\left\{D \in \cald_{\lambda}\mid
W \cap U_D \not= \emptyset, W \not\subset U_D\right\}\right|
\\
& & + \left| \left\{gD^{W'} \mid W' \in \calw^0 \cup
    \dots \cup \calw^{n-1}, g \in G_{B_{\lambda}},
W \cap gW' \not= \emptyset,
                           W \not\subset gW'\right\}\right|
\\
& = & \left|\left\{U_D \in \calu(\cald_{\lambda})\mid W
    \cap U_D \not= \emptyset,
W \not\subset U_D\right\}\right|
\\
& & + \left|\left\{gW' \mid W' \in \calw^0 \cup \dots
    \cup \calw^{n-1},
g \in G_{B_{\lambda}},  W \cap  gW' \not= \emptyset,
                           W \not\subset gW'\right\}\right|
\\
& \le &  \left|\left\{U_D \in  \calu(\cald_{\lambda})\mid
W \cap U_D \not= \emptyset, W \not\subset U_D\right\}\right|
\\
& & + \left|\left\{V \mid V \in \calv^0 \cup \dots \cup
    \calv^{n-1},
W \cap V \not= \emptyset,
                           W \not\subset V\right\}\right|
\\
& \le & k^2_G \cdot (d_X + 1) + 2^{m+1}
\\
& = & M.
\end{eqnarray*}
This finishes the proof of Lemma~\ref{lem:estimate_for_|D_W^bad|}.
\end{proof}

\begin{lemma}
\label{lem:distance-for-points-in-J}
Let $W \in \calw^n$.
If $t_0 \in J^{\good,-}_W \cup J^{\good,+}_W$ and
$t_1 \in J^\dd_W$, then
\[
|t_0-t_1| \not\in [\e - \delta_{r,n-1} + \eta,
                            \alpha+\delta_{r,n-1}-\eta].
\]
\end{lemma}

\begin{proof}
There are $D_0 \in \cald^\good_W$, $D_1 \in \cald_W$,
$\sigma_0$, $\sigma_1 \in \{ -,+ \}$ and $w_0$, $w_1 \in W$
such that $(w_0,t_0) \in \dd_{\sigma_0} D_0^\circ$
and $(w_1,t_1) \in \dd_{\sigma_1} D_1^\circ$.
By definition of $\cald^\good_W$ we have $W \subseteq U_{D_0}$
and $D_0 \subseteq C_\lambda^\circ$.
Therefore there is $\tau \in \IR$ such that
$\Phi_\tau(w_1,t_1) = (w_1,t_1+\tau) \in \dd_{\sigma_0} D_0^\circ$
and $t_1 + \tau \in (c_-,c_+)$.
We conclude
from Lemma~\ref{lem:use-of-friction-control}
that $|t_0 - (t_1 + \tau)| < \eta$.
Because $\cald^{n-1}$ satisfies the induction assumption,
we know that $\cald_W \subseteq \cald^{n-1}$ is
$\delta_{r,n-1}$-overlong.
Therefore
\[
|\tau| \not\in [\e - \delta_{r,n-1}, \alpha+\delta_{r,n-1}].
\]
This implies our result, because
$| |t_0 - t_1| - |\tau| | \leq | (t_0 - t_1) - \tau | < \eta$.
\end{proof}

\begin{definition} \label{def:R_pm}
For $W \in \calw^n$ let
\begin{eqnarray*}
R_-^W & := \sup
    \left( (b_- + \alpha + \e, a_-) \cap J_W^{\good,+} \right)
                \cup \{ b_- + \alpha + \e \}
\\
R_+^W & := \inf
    \left( (a_+, b_+ - \alpha - \e) \cap J_W^{\good,-} \right)
                \cup \{ b_+ - \alpha - \e \}
\end{eqnarray*}
\end{definition}

\begin{lemma}
\label{lem:properties-of-R-W-pm}
Let $W \in \calw^n$. We have:
\begin{enumerate}
\item \label{lem:properties-of-R-W-pm:no-intersection}
      $(R^W_-,R^W_+) \cap J^{\good}_W = \emptyset$.
\item \label{lem:properties-of-R-W-pm:R-m-much-less-R-p}
      $R^W_- + \alpha + \delta_{r,n} + 2\eta < R^W_+$.
\end{enumerate}
\end{lemma}

\begin{proof}
We first show that
\begin{equation}
\label{eq:a-m-a-p-interval-disjoint-from-J-good}
[a_- - \e/2 + \eta, a_+ + \e/2 -\eta] \cap J^{\good}_W = \emptyset.
\end{equation}
We proceed by contradiction.
If \eqref{eq:a-m-a-p-interval-disjoint-from-J-good} fails then
there are
$t_0 \in [a_- - \e/2 + \eta, a_+ + \e/2 - \eta]$,
$w_0 \in W$ and a box
$D \in \cald^\good_W$ such that $(w_0,t_0) \in D^\circ$.
We have $W \subseteq U_{D}$.
For every $w \in W$ there exist  unique real numbers
$\tau_{\pm}(w)$ such that
$(w,\tau_{\pm}(w)) \in \partial_{\pm}D^{\circ}$.
{}From Lemma~\ref{lem:use-of-friction-control} we conclude
$$\tau_{\pm}(w) \in
        (\tau_{\pm}(w_0) - \eta,\tau_{\pm}(w_0) + \eta)
                      \quad \text { for } w \in W.$$
We have $\tau_-(w_0) \leq t_0 \leq \tau_+(w_0)$.
{}From $a_+ - a_- = \e/2$ we conclude
$t_0 \leq a_+ + \e/2 - \eta = a_- + \e - \eta$
and
$t_0 \geq a_- - \e/2 + \eta = a_+ - \e + \eta$.
We estimate
\begin{eqnarray*}
&
\tau_- (w) - \e
 <     \tau_- (w_0) + \eta - \e
 \leq  t_0 + \eta - \e
 \leq  a_-;
& \\ &
\tau_+ (w) + \e
 >     \tau_+ (w_0) - \eta + \e
 \geq  t_0 - \eta + \e
 \geq  a_+.
&
\end{eqnarray*}
This implies $W \x [a_-,a_+] \subset \Phi_{(-\e,\e)} D^\circ$
which contradicts the definition of $\calw^n$ in
Definition~\ref{def:W-is-not-yet-covered}.
This proves \eqref{eq:a-m-a-p-interval-disjoint-from-J-good}.

We give now the proof of
\ref{lem:properties-of-R-W-pm:no-intersection}.
Assume that there is $D \in \cald^\good_W$,
$t_0 \in (R^W_-,R^W_+)$
and $w_0 \in W$ such that $(w_0,t_0) \in D^\circ$.
Because $\eta = \mu/5 < \e/5 < \e/2$ we
conclude from \eqref{eq:a-m-a-p-interval-disjoint-from-J-good}
that either $t_0 < a_-$ or $t_0 > a_+$.
We treat the first case, in the second case there is an analogous
argument.
There is $\tau_+ \geq 0$ such that
$(w_0,t_0 + \tau_+) \in \dd_+ D^\circ$.
If $t_0 + \tau_+ > a_-$, then $(w_0,a_-) \in D^\circ$, i.e.,
$a_- \in J^\good_W$.
Because this contradicts
\eqref{eq:a-m-a-p-interval-disjoint-from-J-good}
we conclude $t_0 + \tau_+ \leq a_-$.
Clearly $t_0 + \tau_+ \in J^{\good,+}_W$
and $b_- + \alpha + \e \leq R_-^W < t_0 \leq t_0 + \tau_+$.
But this is in contradiction to the construction
of $R^W_-$ in Definition~\ref{def:R_pm}.
This proves \ref{lem:properties-of-R-W-pm:no-intersection}.

Next we prove \ref{lem:properties-of-R-W-pm:R-m-much-less-R-p}.
First we treat the case $R_-^W = b_- + \alpha + \e$.
Since $2\eta = 2\mu/5< 2\e/5<\e$, $\delta_{r,n} < \e$ and
$2 \alpha + 3\e < b = a_- - b_-$ we conclude
\begin{multline*}
R_-^W + \alpha + \delta_{r,n} + 2 \eta
~ = ~ b_- + 2\alpha + \e + \delta_{r,n} + 2 \eta
~ < ~ b_- + 2\alpha + 3\e
~ < ~ b_- + b
~ = ~ a_-
~ \leq R_+^W.
\end{multline*}
The case $R_+^W = b_+ - \alpha - \e$ can be treated similarly.
Therefore we may assume now $R_-^W \neq b_- + \alpha + \e$ and
$R_+^W \neq b_+ - \alpha - \e$.
{}From the construction of $R_\pm^W$ we conclude then
that there are $t_\pm \in J^{\good,\mp}_W$,
such that $R^W_- - \eta < t_- \leq R^W_-$ and $
R^W_+ \leq t_+ < R^W_+ + \eta$.
Clearly $t_- \leq R^W_- \leq a_- < a_+ \leq R^W_+ \leq t_+$.
Thus $t_+ - t_- > 0$.
By Lemma~\ref{lem:distance-for-points-in-J}
\[
t_+ - t_- \not\in [\e - \delta_{r,n-1} + \eta,
                            \alpha + \delta_{r,n-1} - \eta].
\]
On the other hand
\eqref{eq:a-m-a-p-interval-disjoint-from-J-good}
implies
$t_- < a_- - \e/2 + \eta$ and $t_+ > a_+ + \e/2 - \eta$.
Using $a_+ - a_- = \e/2$, $2\eta = 2\mu/5 < 2\e/5 < \e/2$ and
$\delta_{r,n-1} - \eta =  \delta_{r,n} + \mu - \eta
                       >  \delta_{r,n} \geq 0$
we estimate
\begin{multline*}
t_+ - t_- ~ > ~ (a_+ + \e/2 - \eta) - (a_- - \e/2 + \eta)
          ~ = ~ 3 \e/2 - 2\eta ~ > ~ \e
          ~ > ~ \e - \delta_{r,n-1} + \eta.
\end{multline*}
Therefore
$t_+ - t_- > \alpha+\delta_{r,n-1}-\eta
           = \alpha+\delta_{r,n}+\mu-\eta$.
This implies $t_- + \alpha + \delta_{r,n} < t_+ + \eta - \mu$.
Using this and $5 \eta = \mu$
we compute
\begin{multline*}
R^W_- + \alpha + \delta_{r,n} + 2 \eta
~ < ~ t_- + \alpha + \delta_{r,n} + 3 \eta
~ < ~ t_+  - \mu + 4 \eta
~ < ~ R^W_+ - \mu + 5 \eta
~ = ~ R^W_+
\end{multline*}
\end{proof}

We can now give the construction of $a_\pm^W$ for $W \in \calw^n$.
If $R_-^W > b_- + \alpha + \e$ then we set
$a_-^W := R_-^W + \eta$.
Otherwise $R_-^W = b_- + \alpha +\e$
and we will use the fact that we arranged $b = a_- - b_-$
to be very large.
It follows from Lemmas~\ref{lem:estimate_for_|D_W^bad|}
and \ref{lem:properties-of-R-W-pm}~
\ref{lem:properties-of-R-W-pm:no-intersection}
and $a_- \leq R_+^W$
that $J^\dd_W \cap [R_-^W,a_-]$
is contained in the union of $2M$
intervals of length $2\eta$.
Using
$\delta_{r,n} + \eta < \delta_{r,n} + \mu = \delta_{r,n-1} < \e$
we estimate
\begin{multline*}
a_-  - (b_- + \alpha + \e)
~ = ~ b - (\alpha + \e)
~ = ~ 4M(\alpha + 2\e) + 2(\alpha + \e) \\
~ > ~ (2M+1) (2\alpha + 2\e)
~ > ~ (2M+1) (2\alpha + 2\eta + 2\delta_{r,n}).
\end{multline*}
If from an interval $I$ of length strictly larger than $L$,
we take out $2M$ or less intervals, each of which
has length less than or equal to $l$,
then the remaining set contains an interval of
length $\tilde{l} :=  (L - 2Ml) / (2M + 1)$.
The center of such an interval, will have distance
$\tilde{l}/2$ from all points in the $2M$ intervals
and from the boundary of $I$.
Therefore we find
$a_-^W \in [(b_- + \alpha+\e) + (\alpha + \delta_{r,n}),
                               a_- - (\alpha + \delta_{r,n})]$
such that
\begin{equation}
\label{eq:a-m-W-is-far-from-things-case-2}
|a_-^W - t| > \alpha + \delta_{r,n}
       \; \mbox{for all} \; t \in J^\dd_W.
\end{equation}
This finishes the construction of $a_-^W$.
To construct $a_+^W$ we proceed similarly.
If $R_+^W < b_+ - \alpha - \e$ then we set
$a_+^W := R_+^W - \eta$.
Otherwise $R_+^W = b_+ - \alpha - \e$ and there is
$a_+^W \in [a_+ + (\alpha + \delta_{r,n}),
     (b_+ - \alpha + \e) - (\alpha + \delta_{r,n})]$
such that
\begin{equation}
\label{eq:a-p-W-is-far-from-things-cse-2}
|a_+^W - t| > \alpha + \delta_{r,n}
       \; \mbox{for all} \; t \in J^\dd_W.
\end{equation}
This finishes the construction of $a_+^W$.
We can arrange that $a_\pm^W = a_\pm^{gW}$
for $g \in G_{B_\lambda}$.

For $W \in \calw^n$ let now
$D^W := B_\lambda(\overline W; a_-^W, a_+^W)
      = \overline{W} \x [a^W_-,a^W_+]$.

\begin{lemma}
\label{lem:good-things-about-D-W}
Let $W \in \calw^n$.
\begin{enumerate}
\item \label{lem:good-things-D-W:in-B-lambda}
      $\Phi_{(-\alpha-\e,\alpha+\e)}  (D^W) \subset B_\lambda$;
\item \label{lem:good-things-D-W:long}
      If $x$ lies in the open bottom or open top of $D^W$ then
      \[
      \Phi_{[-\alpha - \delta_{r,n},-\epsilon + \delta_{r,n}]
         \cup [\epsilon - \delta_{r,n},\alpha + \delta_{r,n}]} (x)
      \]
      does not intersect the open bottom or top of a
      box $D \in \cald^{n-1}$;
\item \label{lem:good-things-D-W:each-long}
      $l_{D^W} = a_+^W - a_-^W > \alpha + \delta_{r,n}$;
\item \label{lem:good-things-D-W:cover}
      $W \x [a_-,a_+] \subseteq \Phi_{(-\e,\e)} ( D^W )^\circ$;
\item \label{lem:good-things-D-W:dim}
      $| \{ D \in \cald^{n-1} \mid
             D^\circ \cap (D^W)^\circ \neq \emptyset \} | \leq M$.
\end{enumerate}
\end{lemma}

\begin{proof}
\ref{lem:good-things-D-W:in-B-lambda}
By construction
$b_- + \alpha + \e < a_-^W < a_+^W < b_+ - \alpha - \e$
and \eqref{eq:calv-Fin-subset} implies
$\overline W \subseteq S_{B_\lambda}$.
\\[1mm]
\ref{lem:good-things-D-W:long}
We consider the open bottom first.
Let $w \in W$ and $x = (w,a_-^W)$.
By \ref{lem:good-things-D-W:in-B-lambda}
$\Phi_{[-\alpha - \delta_{r,n},-\epsilon + \delta_{r,n}]
         \cup [\epsilon - \delta_{r,n},\alpha + \delta_{r,n}]} (x)$
is contained in $B_\lambda$ and can therefore only intersect
boxes from $\cald_W$.
The claim follows thus if
$| a_-^W - t| \not\in [\e-\delta_{r,n},\alpha+\delta_{r,n}]$
for all $t \in J^\dd_W$.
If $R^W_- = b_- + \alpha + \e$ then
\eqref{eq:a-m-W-is-far-from-things-case-2}
holds and implies our claim.
Otherwise $a_-^W = R^W_- + \eta$ and there is
$t_0 \in J^{\good,+}_W \cap [R^W_- - \eta,R^W_-]$
by the construction of $R_-^W$ in Definition~\ref{def:R_pm}.
Now Lemma~\ref{lem:distance-for-points-in-J} implies
our claim since
$\delta_{r,n-1} - \delta_{r,n} = \mu > 3\mu/5 = 3\eta$.
The open top can be treated completely analogously.
\\[1mm]
\ref{lem:good-things-D-W:each-long}
Clearly $a_+^W - a_-^W$ is the length of $D^W$.
By construction $a_-^W \leq a_- + \eta < a_+$ and
$a_+^W \geq a_+ - \eta > a_-$ since
$a_+ - a_- = a = \e/2 > \mu /5 = \eta$.
If $R_-^W = b_- + \alpha + \e$, then by construction
$a_-^W \leq a_- - (\alpha + \delta_{r,n})$
and our claim follows.
Similarly the claim follows if $R_+^W = b_+ - \alpha - \e$.
Thus we are left with the case $a_\pm^W = R_\pm^W \mp \eta$
and the claim follows from
Lemma~\ref{lem:properties-of-R-W-pm}
~\ref{lem:properties-of-R-W-pm:R-m-much-less-R-p}.
\\[1mm]
\ref{lem:good-things-D-W:cover}
As noted above the construction of $a_\pm^W$ implies
that $a_-^W - \eta \leq a_-$ and $a_+^W + \eta \geq a_+$.
The claim follows therefore from $\eta < \e$.
\\[1mm]
\ref{lem:good-things-D-W:dim}
Because $(D^W)^\circ \subset W \x (b_-,b_+)$
\[
\{ D \in \cald^{n-1} \mid
             D^\circ \cap (D^W)^\circ \neq \emptyset \} =
\{ D \in \cald_W \mid
             D^\circ \cap (D^W)^\circ \neq \emptyset \}.
\]
By construction $R_-^W < a_-^W < a_+^W < R_+^W$.
Thus Lemma~\ref{lem:properties-of-R-W-pm}
~\ref{lem:properties-of-R-W-pm:no-intersection}
and Lemma~\ref{lem:estimate_for_|D_W^bad|}
imply
\[
|\{ D \in \cald_W \mid
             D^\circ \cap (D^W)^\circ \neq \emptyset \}| =
|\{ D \in \cald_W^\bad \mid
             D^\circ \cap (D^W)^\circ \neq \emptyset \}| \leq M.
\]
\end{proof}

We now define
\[
\cald^{n} := \cald^{n-1} \cup \{ gD^W \mid W \in \calw^n, g \in G \}.
\]
It remains to check that $\cald^n$ is a $\delta_{r,n}$-good
box cover of $K_\Xi \cup G K_\lambda^{(n)}$.
Recall the induction hypothesis that $\cald^{n-1}$ is a
$\delta_{r,n-1}$-good box cover of $K_\Xi \cup G K_\lambda^{(n-1)}$.

We begin with showing that $\cald^n$ is $\delta_{r,n}$-overlong.
So we have to show for every $x \in X$ which lies on the
open bottom or open top of a box $D_1$ in $\cald^n$, that the set
$\Phi_{[-\alpha - \delta_{r,n},-\epsilon + \delta_{r,n}]
         \cup [\epsilon - \delta_{r,n},\alpha + \delta_{r,n}]} (x)$
does not intersect the open bottom or the open top of
any box $D_2$ in $\cald^n$.

If $D_1$ and $D_2$ lie in $\cald^{n-1}$,
this follows from the induction hypothesis.

Suppose that $D_1 \not\in \cald^{n-1}$ and $D_2 \in \cald^{n-1}$.
Then we can assume without loss of generality that $D_1 = D^W$ for
some $W \in \calw^n$ since $\cald^{n-1}$ is $G$-invariant.
The claim follows then from
Lemma~\ref{lem:good-things-about-D-W}
~\ref{lem:good-things-D-W:long}.

The case $D_1 \in \cald^{n-1}$ and
$D_2 \notin \cald^{n-1}$ is treated analogously.

If $D_1 = D_2$ and $D_1 \not \in \cald^{n-1}$,
then the claim follows from
Lemma~\ref{lem:good-things-about-D-W}
~\ref{lem:good-things-D-W:in-B-lambda} and
\ref{lem:good-things-D-W:each-long}.

If $D_1 \not= D_2$ and $D_1,D_2 \not \in \cald^{n-1}$,
the claim follows from \eqref{eq:calv-j-mutually-disjoint}
and Lemma~\ref{lem:good-things-about-D-W}
~\ref{lem:good-things-D-W:in-B-lambda}
since $B_\lambda$ is a $\Fin$-subset of $X$.
Hence $\cald^n$ is $\delta_{r,n}$-overlong.

We conclude from Lemma~\ref{lem:good-things-about-D-W}
~\ref{lem:good-things-D-W:in-B-lambda} and
\ref{lem:good-things-D-W:each-long}
that $\cald^n$ satisfies \eqref{prop:cover-K-hyp:each-box-long}.

We derive the inclusion $K_\Xi \cup G K_\lambda^{(n)} \subset
     \bigcup_{D \in \cald^n} \Phi_{(-\epsilon,\epsilon)} (D^\circ)$
from
Definition~\ref{def:W-is-not-yet-covered} and
Lemma~\ref{lem:good-things-about-D-W}
~\ref{lem:good-things-D-W:cover}.
(The set $K_\lambda^{(n)}$ was defined before
Lemma~\ref{lem:subinductionstep}.)

By~\eqref{eq:calv-j-mutually-disjoint}
the $D^W$ are mutually disjoint.
Therefore Lemma~\ref{lem:good-things-about-D-W}
~\ref{lem:good-things-D-W:dim}
implies that \eqref{prop:cover-K-hyp:dim} holds for $\cald^{n}$.

It is clear that
\eqref{prop:cover-K-hyp:equivariant}
and \eqref{prop:cover-K-hyp:cofinite} hold for $\cald^{n}$.

Next we prove property~\eqref{prop:cover-K-hyp:Fin-subset}.
Because of the induction hypothesis it suffices to prove the
assertion for the boxes $gD^W$ for $g \in G$ and $W \in \calw^n$,
where $D^W = B_{\lambda}(\overline{W};a_-^W,a_+^W)$.
{}From Lemma~\ref{lem:good-things-about-D-W}
~\ref{lem:good-things-D-W:in-B-lambda} we conclude
$\Phi_{[-\alpha- \epsilon, \alpha + \epsilon]}(D^W)
                                  \subseteq B_{\lambda}(W)$.
Since $\overline{W}$ is a $\Fin$-subset of
$S_{B_{\lambda}}$ with respect to the $G_{B_{\lambda}}$-action
by~\eqref{eq:calv-Fin-subset}
and $B$ is a $\Fin$-subset of the $G$-space $X$,
$\Phi_{[-\alpha- \epsilon, \alpha + \epsilon]}(D^W)$ is a
$\Fin$-subset of the $G$-space $X$.

Finally we show that elements in $\cald^{(n)}$ are not huge.
For $g \in G$ and $W \in \calw^n$ the box $g D^W$ can be
obtained by restriction from $B_{g \lambda}$ and is
therefore not huge, compare Definition~\ref{def:not-huge}
and the subsequent comment.

We have shown that $\cald^n$ is the required $\delta_{r,n}$-good box
cover of $K_\Xi \cup G K_\lambda^{(n)}$.
This finishes the proof of Lemma~\ref{lem:subinductionstep}.

As was noted before,
Proposition~\ref{prop:cover-K-by-boxes-hyp}
follows from
Lemma~\ref{lem:induction_step_for_prop}
which follows from
Lemma~\ref{lem:subinductionstep}.
The proof of Proposition~\ref{prop:cover-K-by-boxes-hyp}
is therefore now completed.


\section{Construction of long $\VCyc$-covers of $X$}
\label{sec:construction-long-VCyc-cover-X}

At the end of this section we will give the
proof of Theorem~\ref{thm:long-thin-cover}.
Throughout this section we will work in the situation of
Convention~\ref{conv:for-long-and-thin-covers}.
In order to construct the long and thin cover
of $X$ we need to discuss
covers of $X^{\IR}$ and $X_{\leq \gamma} - X^\IR$.

\begin{lemma}
\label{lem:cover-X^IR}
There exists a collection $\calu_{X^\IR}$ of open $\Fin$-subsets
of $X$ such that $G\backslash \calu_{X^\IR}$ is finite,
$X^\IR \subset \bigcup_{U \in \calu_{X^\IR}} U$
and $\dim (\calu_{X^\IR}) < \infty$.
\end{lemma}

\begin{proof}
Because the action of $G$ on $X$ is proper there is
for $x \in X^\IR$  an open $\Fin$-neighborhood $W_x$ of $x$.
Because the action of $G$ on $X$ is cocompact and $X^\IR$
is closed, there is a finite subset $\Lambda \subset X^\IR$
such that
$X^\IR \subset \bigcup_{g \in G} \bigcup_{\lambda \in \Lambda}
         g W_\lambda$.
Let $\calu_{X^\IR} = \{ g W_\lambda \; | \;
        g \in G, \lambda \in \Lambda \}$.
Because the $W_\lambda$ are $\Fin$-sets we have
$g W_\lambda \neq W_\lambda \Rightarrow
         g W_\lambda \cap W_\lambda = \emptyset$.
Therefore $\dim(\calu_{X^\IR}) \leq | \Lambda | - 1$.
\end{proof}

\begin{lemma}
\label{covering_X_ge_gamma-X^R}
Fix $\gamma > 0$.
Let $\call_{ \le \gamma}$ be the set of orbits $L = \Phi_{\IR}(x)$
in $X$ whose $G$-period satisfies $0 < \per^G_{\Phi}(L) \le \gamma$.
Then there exists a  collection
$\calu_{\gamma} = \{U_L \mid L \in \call_{\le \gamma}\}$
of open $\VCyc$-subsets $U_L$ of the $G$-space $X$
such that $L \subseteq U_L$ for $L \in \call_{\gamma}$
and $\dim \calu = 0$.
\end{lemma}

\begin{proof}
By assumption we can find  finitely many pairwise
distinct elements $L_1$, $L_2$, $\ldots$,
$L_r$ in $\call_{\le \gamma}$ such that $\call_{\le \gamma} =
G \cdot \{L_1,L_2, \ldots, L_r\}$. We can arrange that $L_j = g
\cdot L_k$ for some $g \in G$ implies $j = k$.
Since the $G_{L_j}$-action on $L_j$ is proper and
cocompact and $L_j$ is homeomorphic to $\IR$ or $S^1$,
the group $G_{L_j}$ is virtually cyclic. 
(A group that acts cocompact and properly on $\IR$ has two ends and
is therefore virtually cyclic, 
\cite[Theorem I.8.32(2)]{Bridson-Haefliger-Buch}.)
We can choose
compact subsets $K_j \subseteq L_j$
with $L_j = G_{L_j} \cdot K_j$ and $G \backslash G L_j$ closed.

Since $G \backslash X$ is compact and 
$G\backslash G L_j \cap G\backslash G L_k \neq \emptyset \Rightarrow j = k$ 
holds,
we can find open subsets $V_1', V_2', \ldots ,V_r'$ in $G\backslash X$ such that
$G\backslash G L_j \subseteq V_j'$ and $V_j' \cap V_k' \not= \emptyset \Rightarrow j = k$ holds.
Let $V_j$ be the preimage of $V_j'$ under the projection $X \to
G\backslash X$.
Then $V_j$ is a $G$-invariant open neighborhood of $L_j$ and $V_j \cap
V_k \not= \emptyset \Rightarrow j = k$.

Fix $j \in \{1,2,\ldots r\}$.  Since the $G$-action on $X$ is proper, we can find an open
neighborhood $W_j'$ of $K_j$ and a finite subset $S \subseteq G$ such that $W_j' \cap
gW_j' \neq \emptyset\Rightarrow g \in S$.  Let $S_0 \subseteq S$ be the subset consisting
of those elements $g \in S$ for which $K_j \cap g K_j = \emptyset$.  Since $K_j$ is
compact, we can find for $s \in S_0$ an open neighborhood $W_s''$ of $K_j$ such that
$W_s'' \cap s W_s'' = \emptyset$.  Put $U_j' := W_j' \cap \bigcap_{s \in S_0} W_s''$.
Then $U_j'$ is an open neighborhood of $K_j$ such that $U_j' \cap g U_j' \not= \emptyset$
implies $K_j \cap gK_j \neq \emptyset$.  Put $U_j = G_{L_j} \cdot U_j'$. Then $U_j$ is a
$G_{L_j}$-invariant open subset containing $L_j = G_{L_j}K_j$. Next we prove $g U_j \cap
U_j \not= \emptyset \Rightarrow g \in G_{L_j}$. Suppose for $g \in G$ that $g U_j \cap U_j
\not= \emptyset$.  Then we can find $g_0,g_1 \in G_{L_j}$ such that $gg_0 U_j' \cap g_1U_j
\not= \emptyset$.  This implies $g_1^{-1}gg_0 K_j \cap K_j \not= \emptyset$. We conclude
$g_1^{-1}gg_0L_j \cap L_j \not= \emptyset$ and hence $g_1^{-1}gg_2L_j = L_j$. This shows
$g_1^{-1}gg_0 \in G_{L_j}$ and thus $g \in G_{L_j}$.  Hence $U_j$ is an open
$\VCyc$-subset of the $G$-space $X$ such that $G_{U_j} = G_{L_j}$ and $L_j \subseteq U_j$.

Define for any element $L \in \call_{\le \gamma}$
$$U_L = g \cdot (V_j \cap U_j) \quad \text{ for } g \in G 
 \text{ with } L = gL_j.$$
This is independent of the choice of $g$ and $j$ and $U_L$ is a
$\VCyc$-subset of the $G$-space $X$ with $G_{U_L} = G_L$ since
$(V_j \cap U_j)$ is a $\VCyc$-subset of the $G$-space $X$ with
$G_{V_j \cap U_j} = G_{L_j}$. We have by construction
$$U_{L_1} \cap U_{L_2} \not= \emptyset \Rightarrow L_1 = L_2.$$
\end{proof}

Finally we can give the proof of
Theorem~\ref{thm:long-thin-cover}.

\begin{proof}
Let $\alpha > 0$ be given.
Choose $\epsilon$ such that $0 < \epsilon < \alpha$.
Let $M = M(k_G,d_X)$ and
$\gamma = \gamma(4\alpha,\epsilon,M) > 0$
be as in Proposition~\ref{prop:cover-K-by-boxes-hyp}.

Let $\calu_{X^\IR}$ be the collection of open
$\Fin$-sets from  Lemma~\ref{lem:cover-X^IR}
and $\calu_\gamma$ be the collection of open $\VCyc$-sets
from Lemma~\ref{covering_X_ge_gamma-X^R}.
Note that
$\dim(\calu_{X^\IR} \cup \calu_\gamma) = \dim(\calu_{X^\IR}) + 1$
is finite and does not depend on $\alpha$, but only on an arbitrarily 
small neighborhood of $X^\IR$ as a $G$-space.

Put
\begin{equation*}
S = \{ x \in X \mid \exists U \in \calu_\gamma \cup \calu_{X^{\IR}}
       \mbox{ such that } \Phi_{[-\alpha,\alpha]} (x) \subset U \}.
\end{equation*}
Note that $S$ is $G$-invariant, because $\calu_\gamma$ and 
$\calu_{X^\IR}$ are.
Consider $x \in S$.
Choose $U_x \in \calu_\gamma \cup \calu_{X^\IR}$ with
$\Phi_{[-\alpha,\alpha]} (x) \subseteq U$.
Since $\{ x \} \x [-\alpha,\alpha]$ is compact and contained in
$\Phi^{-1}(U_x)$,
we can find an open neighborhood $V_x$ of $x$ such that
$V_x \times [-\alpha,\alpha] \subseteq \Phi^{-1}(U_x)$.
This implies $V_x \subseteq S$.
Hence $S$ is an open $G$-invariant subset of $X$
which contains $X_{\le \gamma}$.

Let $K$ be the closure of the complement of $S$ in $X$.
Since $G\backslash X$ is compact
$K \subseteq X$ is a cocompact $G$-invariant subset
which does not meet $X_{\leq \gamma}$.
Hence we can apply Proposition~\ref{prop:cover-K-by-boxes-hyp}
to $K$ with respect to $4\alpha$ instead of $\alpha$ and $\epsilon$
with $0 < \epsilon < \alpha$.
Recall that $M = M(k_G,d_X)$ and
$\gamma = \gamma(M,4\alpha,\epsilon)$
are the numbers appearing in
Proposition~\ref{prop:cover-K-by-boxes-hyp}.
So we get a collection of boxes $\cald$ with the properties
described in Proposition~\ref{prop:cover-K-by-boxes-hyp}.
Put
\begin{equation*}
\calu_{K} = \{ \Phi_{(-\alpha-\epsilon,+\alpha+\epsilon)}
                     D^{\circ} \mid D \in \cald \}.
\end{equation*}
Then for $x \in K$ there is $U \in \calu_K$ such that
$\Phi_{[-\alpha,\alpha]}(x) \subset U$ and every element
in $\calu_K$ is an open $\Fin$-subset of $X$.

Next we show $\dim ( \calu_K ) \leq 2M+1$.
Consider pairwise disjoint elements $D_1$, $D_2$, $\ldots$,
$D_{2M+3}$ of $\cald$.
We have to show that
$\bigcap_{k=1}^{2M+3}
  \Phi_{(-\alpha-\epsilon,\alpha + \epsilon)}(D_k^{\circ})
                                                 = \emptyset.$
Suppose the contrary, i.e., there exist $x \in X$ such that
$x \in  \Phi_{(-\alpha-\epsilon,\alpha + \epsilon)}(D_k^{\circ})$
holds for $k = 1,2, \ldots , (2M+3)$.
Obviously
$x \in \Phi_{(-\alpha-\epsilon,\alpha + \epsilon)}(D_k^{\circ})$
implies that $\Phi_{2\alpha}(x) \in D_k^{\circ}$ or
$\Phi_{-2\alpha}(x) \in D_k^{\circ}$ since
$\epsilon < \alpha$ and for every $y \in X$ the set
$\Phi_{[0,4\alpha]}(y)$ can not intersect both
the open bottom and the open top of $D_k$.
Hence we can find $(M+2)$ pairwise distinct elements
$k_1, k_2, \ldots, k_{M+2} \subseteq \{1,2, \ldots, 2M+3\}$ such
that $\Phi_{2\alpha}(x) \in D_{k_j}^{\circ}$ holds
for $j = 1,2, \ldots, M+2$ or that
$\Phi_{-2\alpha}(x) \in D_{k_j}^{\circ}$ holds for
$j = 1,2, \ldots, M+2$.
In both  cases we get a contradiction to
$\dim (\{D^{\circ} \mid D \in \cald\} ) \leq M$.
This shows  $\dim ( \calu_K ) \leq 2M+1$.
Note that this bound depends only on $k_G$ and $d_X$
and is independent of $\alpha$.
Because
$$
\dim ( \calu_K \cup \calu_{\gamma} \cup \calu_{X^\IR} )
~ \leq ~
\dim ( \calu_K ) + \dim ( \calu_{\gamma} \cup \calu_{X^\IR} ) + 1
$$
this implies that the dimension of
$$\calu ~ :=~ \calu_K \cup \calu _{\gamma} \cup \calu_{X^\IR}$$
is bounded by a number that depends only on $k_G$, $d_X$ and
the $G$-action on an arbitrary small neighborhood of $X^\IR$.
Thus $\calu$ is
the required $\VCyc$-cover of $X$.
This finishes the proof of
Theorem~\ref{thm:long-thin-cover}.
\end{proof}

\begin{remark}
\label{rem:weaken-assumption-on-isolated-short-orbits}
In Convention~\ref{conv:for-long-and-thin-covers} we
assumed that the number of closed orbits,
which are not stationary and
whose  period is $ \le C$,
of the flow induced on $G\backslash X$ is finite
for every $C > 0$.
This assumption can be replace with the following less
restrictive assumption.
\begin{quote}
There is a number $N$ such that for every $\gamma > 0$
there is a $G$-invariant collection $\calu$ of open $\VCyc$-subsets of
$X$ such that for each $x \in X_{\leq \gamma}$ there is
$U \in \calu$ such that $\Phi_{[-\gamma,\gamma]} \subset U$,
$\dim ( \calu ) \leq N$ and $G \backslash \calu$ is finite.
\end{quote}
The proof of Theorem~\ref{thm:long-thin-cover}
given above clearly also works under this less restrictive
assumption.
This might be useful in nonpositively curved situations.
\end{remark}


\section{Mineyev's flow space}
\label{sec:Mineyevs_flow_space}

\subsection{Hyperbolic complexes, double difference and Gromov product}
\label{subsec:Hyperbolic complexes, double difference and Gromov product}

We collect some basic concepts such as hyperbolic complexes,
double difference and Gromov product
which are all taken from the paper by
Mineyev~\cite{Mineyev-flows-and-joins}
and which we will need for our purposes.

A simplicial complex is called \emph{uniformly locally finite} if 
there exists a number $N$ such that any
element in the $0$-skeleton $X_0$ occurs as vertex of at most $N$ simplices.

Let $X$ be a simplicial complex.
Given any metric $d$ on its $0$-skeleton, one can extend it to a 
metric $\widetilde{d}$ on $X$
as follows. Given points $u_k$ for $k = 1,2$ in $X$, we can find vertices
$x_k[0]$, $x_k[1]$, $x_k[n_k]$ such that $u_k$ belongs to the simplex 
with vertices
$x_k[0]$, $x_k[1]$, $\ldots$,  $x_k[n_k]$. 
There are unique numbers $\alpha_k[0]$, 
$\alpha_k[1]$, $\ldots$, $\alpha_k[n_k]$
in $[0,1]$ with $\sum_{i_k = 0}^{n_k} \alpha_k[i_k] = 1$ 
such that $u_k$ is given by
$\sum_{i_k = 0}^{n_k} \alpha_k[i_k] \cdot x_k[n_k]$. 
Then  define
$$\widetilde{d}(u_1,u_2) ~ = ~ \sum_{i_1 = 0}^{n_1}  \sum_{i_2 = 0}^{n_2} \alpha_1[i_1] \cdot \alpha_2[i_2] \cdot d(x_1[i_1],x_2[i_2]).$$
This is a well-defined metric extending $d$ such that each simplex with the metric induced by $\widetilde{d}$
is homeomorphic to the standard simplex.

Given a connected simplicial complex, define a metric $d_0$ on its $0$-skeleton by
defining $d_0(x,y)$ as the minimum of the numbers $n \ge 0$ such that there is a sequence
of vertices $x = x_0$, $x_1$, $ \ldots $, $x_n = y$ with the property that $x_i$ and $x_{i+1}$ are joint by an edge
for $i= 0,1, \ldots, n-1$.
The \emph{word metric} $d_{\mathit{word}}$ on  a  connected simplicial complex $X$
is the metric $\widetilde{d_0}$.

A \emph{metric complex}  $(X,d)$ is a connected uniformly locally finite simplicial complex $X$
equipped with its word metric $d = d_{\mathit{word}}$.
A \emph{hyperbolic complex} $X$ is a
metric complex $(X,d)$ such that $(X,d)$ is $\delta$-hyperbolic in the sense of Gromov for some $\delta > 0$ (see~\cite{Gromov(1987)},
\cite[Definition III.H.1.1]{Bridson-Haefliger-Buch}).
Let $\partial X$ be the boundary and $\overline{X} = X \cup \partial X$ be the compactification
of the hyperbolic complex $X$ in the sense of Gromov (see~\cite{Gromov(1987)}, \cite[III.H.3]{Bridson-Haefliger-Buch}).

Mineyev~\cite[6.1]{Mineyev-flows-and-joins} constructs for a hyperbolic
metric complex $(X,d)$ a new metric $\widehat{d}$ with certain
properties (see~\cite[Lemma~2.7 on page 449 and Theorem 32 on page
446]{Mineyev-flows-and-joins}). For instance $\widehat{d}$ is quasiisometric
to the word metric $d_{\mathit{word}}$. For a simplicial map $f
\colon X \to X$ the following conditions are equivalent: (i) $f$
is a simplicial automorphism, (ii) $f$ is a simplicial
automorphism preserving the word metric $d_{\mathit{word}}$, (iii)
$f$ is a simplicial automorphism preserving the metric
$\widehat{d}$.
Define for $a,a',b,b' \in X$ the \emph{double
difference} to be the real number
\begin{eqnarray}
\label{double_difference}
\langle a,a'|b,b' \rangle  & := &  \frac{1}{2} \cdot \left(\widehat{d}(a,b) + \widehat{d}(a',b') - \widehat{d}(a',b) - \widehat{d}(a,b')\right).
\end{eqnarray}
Recall that the \emph{Gromov  product} for $a,b,c \in X$ is defined to be the
positive real number
\begin{eqnarray}
\label{Gromov_product}
\langle a|b \rangle_c  & := &  \frac{1}{2} \cdot \left(\widehat{d}(a,c) + \widehat{d}(b,c) - \widehat{d}(a,b)\right).
\end{eqnarray}

Define the subset
$$S(\overline{X}) ~ \subseteq  ~
\left\{(a,a',b,b') \in \overline{X} \times \overline{X} \times \overline{X} \times \overline{X}\right\} = \overline{X}^4$$
to consist of those quadruples $(a,a',b,b')$ satisfying
\begin{eqnarray*}
a,b \in \partial X & \Rightarrow & a \not= b;
\\
a,b' \in \partial X & \Rightarrow & a' \not= b;
\\
a',b \in \partial X & \Rightarrow & a \not= b';
\\
a',b' \in \partial X & \Rightarrow & a' \not= b'.
\end{eqnarray*}
Let
$$T(\overline{X}) ~ := ~ \{(a,b,c) \mid \overline{X} \times \overline{X} \times \overline{X} \mid
c \in \partial X \Rightarrow (a \not=c \text{ and } b \not= c)\}.$$
We equip $S(\overline{X}) \subseteq \overline{X}^4$ and $T(\overline{X}) \subseteq \overline{X}^3$ with the subspace topology.
The following result is a special case of \cite[Theorem~35~on page~448 and Theorem~36 on page~452]{Mineyev-flows-and-joins}.
(We only need and want to consider the case where the double difference takes values in $\IR$.)

\begin{theorem}[Mineyev]
\label{the:extending_the_double_difference}
Let $(X,d)$ be a hyperbolic complex. Then the double difference
of~\eqref{double_difference} extends to a continuous function invariant under simplicial
automorphisms of $X$
$$\langle -,- | -, - \rangle \colon S(\overline{X}) \to \IR$$
satisfying
\begin{enumerate}

\item $\langle a,a'|b,b' \rangle ~ = ~ \langle b,b'|a,a' \rangle$;

\item $\langle a,a'|b,b' \rangle ~ = ~ - \langle a',a|b,b' \rangle ~ = ~ - \langle a,a'|b',b \rangle$;

\item $\langle a,a|b,b' \rangle ~ = ~ \langle a,a'|b,b \rangle ~ = ~ 0$;

\item $\langle a,a'|b,b' \rangle + \langle a',a''|b,b' \rangle ~ = ~ \langle a,a''|b,b' \rangle$;

\item $\langle a,b|c,x \rangle + \langle c,a|b,x \rangle +  \langle b,c|a,x \rangle ~ = ~ 0$;
\end{enumerate}

The Gromov product of \eqref{Gromov_product} extends to a
continuous function invariant under simplicial automorphisms of
$X$
$$\langle -|-\rangle_- \colon T(\overline{X}) \to [0,\infty]$$
satisfying
\begin{enumerate}

\item $\langle a|b\rangle_c = \infty ~ \Leftrightarrow ~ (c \in \partial X) \text { or } \left(a,b \in \partial X \text{ and } a = b\right)$;

\item $\langle a,b|x,y\rangle ~ = ~ \langle b|x\rangle_a - \langle b|y\rangle_a$ for $a \in X$ and
      $(a,b,x,y) \in S(\overline{X})$.

\end{enumerate}

\end{theorem}
We will often use the rules appearing in
Theorem~\ref{the:extending_the_double_difference} tacitly.

Another important ingredient will be the following result due to Mineyev~\cite[Proposition~38 on page~453]{Mineyev-flows-and-joins}.
\begin{proposition}[Mineyev]\label{prop:Mineyevs_Proposition_38}
Let $(X,d)$ be a hyperbolic complex. Then there exist constants $\lambda \in (e^{-1},1)$ and $T \in [0,\infty)$
depending only on $X$ such that for all $a,b,c,u \in \overline{X}$ satisfying
\begin{eqnarray*}
(a,c,u,b) & \in & S(\overline{X});
\\
(b,c,u,a) & \in & S(\overline{X});
\\
\max\{\langle a,c|u,b\rangle, \langle b,c|u,a\rangle\} & \ge & T,
\end{eqnarray*}
we have
\begin{eqnarray*}
(u,c,a,b) & \in & S(\overline{X});
\\
\left|\langle u,c|a,b \rangle\right| & \le & \lambda^{\max\{\langle a,c|u,b\rangle, \langle b,c|u,a\rangle\}}.
\end{eqnarray*}
\end{proposition}


\subsection{Two auxiliary functions}
\label{subsec:Two_auxiliary_functions}

In the sequel we will use the following two functions for $\alpha,\beta \in \overline{\IR} := \IR \coprod \{-\infty,\infty\}$
with $\alpha \le \beta$.
\begin{eqnarray}
\theta_{[\alpha,\beta]} \colon \overline{\IR} & \to & [\alpha,\beta]
\label{theta_[alpha,beta]}
\\
\Theta_{[\alpha,\beta]} \colon \overline{\IR} & \to & [\alpha,\beta]
\label{Theta_[alpha,beta]}
\end{eqnarray}
which are defined by
$$\theta_{[\alpha,\beta]}(t) ~ := ~ \left\{ \begin{array}{lll}
\alpha & & \mbox{ if } - \infty \le t \le \alpha;
\\
t & & \mbox{ if } \alpha \le  t \le \beta;
\\
\beta   & & \mbox{ if } \beta \le  t \le \infty,
\end{array}\right.
$$
and
$$\Theta_{[\alpha,\beta]}(t) ~ := ~ \left\{ \begin{array}{lll}
-\infty & & \mbox{ if } - \infty = t = \alpha,
\\
\alpha + e^{t-\alpha}/2 - e^{t-\beta}/2 & & \mbox{ if } - \infty \le t \le \alpha, -\infty < \alpha < \infty;
\\
t + e^{\alpha-t}/2 - e^{t-\beta}/2 & & \mbox{ if } \alpha \le  t \le \beta , -\infty < t < \infty;
\\
\beta  + e^{\alpha-t}/2 - e^{\beta-t}/2 & & \mbox{ if } \beta \le  t \le \infty, - \infty < \beta < \infty;
\\
\infty & & \mbox{ if } t = \beta = \infty.
\end{array}\right.
$$
Here and in the sequel we use the convention that
for $r,s \in\IR$ the expressions $r + s$, $r \cdot s$ and
$e^r$ are defined as usual and furthermore
$$
\begin{array}{lcl}
r + \infty = \infty + r = \infty                 & &
                              \mbox{ for } r \in \IR;
\\
r - \infty = - \infty + r = - \infty                 & &
                              \mbox{ for } r \in \IR;
\\
e^{-\infty} = 0;                                 & &
\\
e^{\infty} = \infty;                             & &
\\
|\pm \infty| = \infty.                           & &
\end{array}
$$

The function $\Theta_{[\alpha,\beta]}$ agrees with the function
denoted by $\theta'[\alpha,\beta;\cdot]$
in Mineyev~\cite[Section~1.6]{Mineyev-flows-and-joins}.
We equip $\overline{\IR}$ with the topology uniquely
determined by the properties that
$\IR \subseteq \overline{\IR}$ is an open subset,
the subspace topology on $\IR \subseteq \overline{\IR}$
is the standard topology and  a fundamental system for
open neighborhoods of $\infty$ is
$\{(R,\infty)  \cup \{\infty\} \mid R \in \IR\}$
and of $- \infty$ is
$\{(-\infty,R)  \cup \{-\infty\} \mid R \in \IR\}$.

The elementary proof of the following basic
properties of $\Theta_{[\alpha,\beta]}$ is left to the reader.

\begin{lemma} \label{lem:properties_of_theta_[alpha,beta]}
Suppose that $\alpha < \beta$. Then:
\begin{enumerate}

\item \label{lem:properties_of_theta_[alpha,beta]:description_by_an_integral}
We have for $t \in \IR$
$$\Theta_{[\alpha,\beta]}(t) ~ = ~ \int_{-\infty}^{\infty} \frac{\theta_{[\alpha,\beta]}(t+s)}{2 \cdot e^{|s|}} ~ ds
~ = ~ \theta_{[\alpha,\beta]}(t) + e^{-|\alpha-t|}/2 - e^{-|\beta-t|}/2;$$
\item \label{lem:properties_of_theta_[alpha,beta]:C^1}
The restriction of $\Theta_{[\alpha,\beta]}$ to $\IR$ is a homeomorphism
$\IR \xrightarrow{\cong}  (\alpha,\beta)$
which is a $C^1$-function. Its first derivative is the continuous function
$$
t \in \IR ~\mapsto ~
\left\{ \begin{array}{lll}
e^{t-\alpha}/2 - e^{t-\beta}/2 & & \mbox{ if } - \infty < t \le \alpha;
\\
1 - e^{\alpha-t}/2 - e^{t-\beta}/2 & & \mbox{ if } \alpha \le  t \le \beta;
\\
-e^{\alpha-t}/2 + e^{\beta-t}/2 & & \mbox{ if } \beta \le  t < \infty;
\end{array}\right.
$$

\item \label{lem:properties_of_theta_[alpha,beta]:Monotonicity}
The function $\Theta_{[\alpha,\beta]}$ is strictly monotone increasing.
The function $\theta_{[\alpha,\beta]}$ is monotone increasing;

\item \label{lem:properties_of_theta_[alpha,beta]:non-expanding}
The function $\Theta_{[\alpha,\beta]}$ is non-expanding, i.e., $|\Theta_{[\alpha,\beta]}(t) - \Theta_{[\alpha,\beta]}(s)| \le |t-s|$ for
$t,s \in \overline{\IR}$. The same is true for $\theta_{[\alpha,\beta]}$;

\item \label{lem:properties_of_theta_[alpha,beta]:homeomorphism}
The map $\Theta_{[\alpha,\beta]} \colon \overline{\IR} \to [\alpha,\beta]$ is a homeomorphism;

\item \label{lem:properties_of_theta_[alpha,beta]:shift_invariance}
We have for $t \in \overline{\IR}$ and $s \in \IR$
\begin{eqnarray*}
\theta_{[\alpha + s,\beta + s]}(t + s) & = & \theta_{[\alpha ,\beta]}(t) + s;
\\
\Theta_{[\alpha + s,\beta + s]}(t + s) & = & \Theta_{[\alpha ,\beta]}(t) + s;
\end{eqnarray*}

\item  \label{lem:properties_of_theta_[alpha,beta]:dependency_on_alpha_and_beta_part1}
Consider $\alpha_0,\alpha_1,\beta_0,\beta_1 \in \overline{\IR}$ such that
$\alpha_i \le \beta_i$ for $i = 0,1$, 
$(\alpha_0 = - \infty \Leftrightarrow \alpha_1 = - \infty)$
and
$(\beta_0 =  \infty \Leftrightarrow \beta_1 =  \infty)$
holds. Put
$$C ~ := ~ \left\{
\begin{array}{lcl}
\max\{|\alpha_1 - \alpha_0|,|\beta_1 - \beta_0|\} & & \text{if } \alpha_0,\alpha_1,\beta_0,\beta_1 \in \IR;
\\
|\alpha_1 - \alpha_0| & & \text{if } \alpha_0,\alpha_1,\in \IR, \beta_0 = \beta_1 =  \infty;
\\
|\beta_1 - \beta_0| & & \text{if }  \alpha_0 = \alpha_1 =  - \infty, \beta_0,\beta_1,\in \IR;
\\
0  & & \text{if }  \alpha_0 = \alpha_1 =  - \infty, \beta_0 = \beta_1 = \infty.
\end{array}
\right.;
$$
Then we get for all $t \in \IR$ that
\begin{eqnarray*}
\left|\theta_{[\alpha_1,\beta_1]}(t) - \theta_{[\alpha_0,\beta_0]}(t)\right| & \le & C;
\\
\left|\Theta_{[\alpha_1,\beta_1]}(t) - \Theta_{[\alpha_0,\beta_0]}(t)\right| & \le & C;
\end{eqnarray*}

\item  \label{lem:properties_of_theta_[alpha,beta]:dependency_on_alpha_and_beta_part2}
If $\alpha \le t$, then
$$\Theta_{[\alpha,\beta]}(t) - \Theta_{[-\infty,\beta]}(t) ~ =~ e^{\alpha -t}/2.$$
If $\beta \ge t$, then
$$\Theta_{[\alpha,\infty]}(t) - \Theta_{[\alpha,\beta]}(t) ~ = ~ e^{t - \beta}/2;$$

\item \label{lem:properties_of_theta_[alpha,beta]:bound_on_t}
Consider $\alpha, \beta \in \overline{\IR}$ with $-\infty < \alpha <  \beta$ and $t \in \IR$.

If $\Theta_{[\alpha,\beta]}(t)  ~ \le ~ \frac{\alpha + \beta}{2}$, then
$$t ~ \le \min\{\beta,\Theta_{[\alpha,\beta]}(t) + 1/2\}.$$
If $\Theta_{[\alpha,\beta]}(t)  ~ \ge ~ \frac{\alpha + \beta}{2}$, then
$$t ~ \ge \max\{\alpha,\Theta_{[\alpha,\beta]}(t) - 1/2\}.$$

\end{enumerate}
\end{lemma}


\subsection{The construction of the flow space}
\label{subsec:The_construction_of_the_geodesic_flow_space}

Let $(X,d)$ be a hyperbolic complex.
We want to define the associated \emph{flow space},
i.e.\ a metric space $\FS(X)$ together with a flow,
following Mineyev~\cite{Mineyev-flows-and-joins}.
(It is the same as the
half open symmetric join $\sjX$ constructed by
Mineyev \cite[Section 8.3]{Mineyev-flows-and-joins}.)
The underlying set is
\begin{multline}
\FS(X) :=  ~
  \left\{(a,b,t) \in \overline{X} \times
              \overline{X} \times \overline{\IR} \mid
    (a \in \partial X \Rightarrow t \not= - \infty) \right.
\\
\left. \text{ and }
     (b \in \partial X \Rightarrow t \not=  \infty)
     \text{ and }
     (a,b \in \partial X \Rightarrow a \not=  b)
                                    \right\}/ \sim,
\label{FS(X)}
\end{multline}
where we identify $(a,b,-\infty) \sim (a,b',-\infty)$,
$(a,b,\infty) \sim (a',b,\infty)$, and $(a,a,t) \sim (a,a,t')$.
In the sequel we will denote for
$(a,b,t) \in \overline{X} \times \overline{X}
                          \times \overline{\IR} $
which satisfies
$a \in \partial X \Rightarrow t \not= - \infty$,
$b \in \partial X \Rightarrow t \not=  \infty$ and
$a,b \in \partial X \Rightarrow a \not=  b$
its class in $\FS{X}$ again by $(a,b,t)$.

{}From now on we fix a base point $x_0 \in X$.
The metric on $\FS(X)$ will depend on this choice.

Define the map
\begin{multline}
l_{x_0} \colon X \times \FS(X) ~ \to ~ \IR
\label{l_{x_0}}
\\
(u,(a,b,t)) ~ \mapsto ~ \langle a|b\rangle_u +
\left|\theta_{[-\langle b|x_0\rangle_a,\langle a|x_0\rangle_b]}
                    (t) - \langle a,b|u,x_0\rangle \right|.
\end{multline}

It is easy to check that it is compatible with the equivalence
relation appearing in the definition of $\FS(X)$.

\begin{definition} \label{def:pseudometric_d^{times}_{FS,x_0}}
Define a  pseudometric on $\FS(X)$
\begin{multline*}
d^{\times}_{\FS,x_0} = d^{\times}_{\FS(X),x_0} \colon \FS(X) \times \FS(X) \to \IR,
\\
(a,b,t),(c,d,s) ~ \mapsto \sup_{u \in X} \left|l_{x_0}(u, (a,b,t)) - l_{x_0}(u,(c,d,s))\right|.
\end{multline*}
\end{definition}

Recall that a pseudometric satisfies the same axioms as a metric except
that the condition $d(x,y) = 0 \Rightarrow x= y$ is dropped.
The proof that this definition makes sense and yields a
pseudometric $d^{\times}_{\FS,x_0}$  is given in
Mineyev~\cite[Theorem~44 on page~459]{Mineyev-flows-and-joins}.

\begin{lemma} \label{lem:X_to_FS(X)_is_isometric_for_d^{times}}
The inclusion $X \to \FS(X), x \mapsto (x,x,0)$ is an isometric embedding with
respect to the  metric $\widehat{d}$ on $X$ and  the pseudometric
$d^{\times}_{\FS,x_0}$ on $\FS(X)$.
\end{lemma}
\begin{proof}
We compute for $u \in X$ and $x \in X$ and $t \in \overline{\IR}$
\begin{eqnarray*}
l_{x_0}(u,(x,x,t))
& = &
\langle x|x\rangle_u + \left|\theta_{[-\langle x|x_0\rangle_x,\langle
    x,x_0\rangle_x]}(t) - \langle x,x|u,x_0\rangle \right|
\\
& = &
\widehat{d}(x,u) + \left|\theta_{[0,0]}(t) - 0 \right|
\\
& = & \widehat{d}(u,x).
\end{eqnarray*}
Consider $a,b \in X$. Since by the triangle inequality
$|\widehat{d}(u,a) - \widehat{d}(u,b)| \le \widehat{d}(a,b)$  and  $|\widehat{d}(b,a) - \widehat{d}(b,b)| = \widehat{d}(a,b)$ holds, we
conclude
\begin{eqnarray*}
d^{\times}_{\FS,x_0}((a,a,t),(b,b,s))
& = &
\sup_{u \in X} \left|l_{x_0}(u,(a,a,t)) - l_{x_0}(u,(b,b,s))\right|
\\
& = &
\sup_{u \in X} |\widehat{d}(u,a) -\widehat{d}(u,b)|
\\
& = & \widehat{d}(a,b).
\end{eqnarray*}
\end{proof}

The canonical $\IR$-action on $\overline{\IR}$
\begin{eqnarray}
\phi \colon \IR \times \overline{\IR}, \quad (\tau,t) \mapsto \phi_{\tau}(t)
\label{R-action_on_overline{R}}
\end{eqnarray}
is defined by $\phi_{\tau}(t) = t + \tau$, if $t \in \IR$, $\phi_{\tau}(-\infty) = - \infty$ and
$\phi_{\tau}(\infty) = \infty$. This $\IR$-action on $\overline{\IR}$
together with the trivial $\IR$-action on $\overline{X} \times \overline{X}$
yields an action of $\IR$ on $\overline{X} \times \overline{X} \times \overline{\IR}$ which in turn
induces an $\IR$-action
\begin{eqnarray}
\phi \colon \IR \times \FS(X) & \to & \FS(X) , \quad (\tau,(a,b,t)) ~ \mapsto (a,b,\phi_{\tau}(t)).
\label{flow_phi}
\end{eqnarray}

For $a,b \in \overline{X}$ we define the
\emph{line} $(a,b)_{\FS(X)}$ to be the set of points
$\{(a,b,t) \mid t \in \IR\}$.
Obviously
$(a,b)_{\FS(X)}$ is a transitive free $\IR$-set if $a \neq b$.

Next we construct the desired metric from the pseudometric above.

\begin{definition}
\label{def:metric_d_{FS,x_0}}
Define a  metric on $\FS(X)$
\begin{multline*} d_{\FS(X),x_0} \colon \FS(X) \times \FS(X) \to \IR,
\\
(a,b,t),(c,d,s) ~ \mapsto \int_{\IR} \frac{d^{\times}_{\FS,x_0}(\phi_{\tau}(a,b,t),\phi_{\tau}(c,d,s))}{2 \cdot e^{|\tau|}} d\tau.
\end{multline*}
\end{definition}

Obviously $d_{\FS,x_0}$ inherits from $d^{\times}_{\FS,x_0}$ the properties of a pseudometric.
The proof that  $d_{\FS,x_0}$ is a metric can be found in
\cite[Theorem~14 on page~426 and Theorem 45 on page~459]{Mineyev-flows-and-joins}.

\begin{lemma} \label{lem:X_to_FS(X)_is_isometric_for_d}
The inclusion $X \to \FS(X)$ is an isometric embedding with
respect to the  metric $\widehat{d}$ on $X$ and  the metric  $d_{\FS,x_0}$ on $\FS(X)$.
\end{lemma}
\begin{proof}
We compute for $u \in X$ and $x \in X$ and $t \in \overline{\IR}$
using Lemma~\ref{lem:X_to_FS(X)_is_isometric_for_d^{times}}.
\begin{eqnarray*}
d_{\FS,x_0}((a,a,t),(b,b,s))
& = &
\int_{\IR} \frac{d^{\times}_{\FS,x_0}(\phi_{\tau}(a,a,t),\phi_{\tau} (b,b,s))}{2 \cdot e^{|\tau|}}~ d\tau
\\
& = &
\int_{\IR} \frac{d^{\times}_{\FS,x_0}((a,a,t + \tau),(b,b,s + \tau))}{2 \cdot e^{|\tau|}}~ d\tau
\\
& = &
\int_{\IR} \frac{\widehat{d}(a,b)}{2 \cdot e^{|\tau|}}~ d\tau
\\
\\
& = &
\widehat{d}(a,b) \cdot \int_{\IR} \frac{1}{2 \cdot e^{|\tau|}}~ d\tau
\\
& = &
\widehat{d}(a,b).
\end{eqnarray*}
\end{proof}

In the sequel we will consider $X$ as a subspace of $\FS(X)$ by the isometric embedding
appearing in Lemma~\ref{lem:X_to_FS(X)_is_isometric_for_d}.

Let $f \colon X \to X$ be an \emph{isometry},
i.e. a bijection respecting the word metric $d$ on $X$.
It extends uniquely to a  homeomorphism
$\overline{f} \colon \overline{X} \to \overline{X}$
and induces an isometry with respect to the metric
$d_{\FS(X),x_0}$
\begin{eqnarray*}
\FS(f) \colon \FS(X) & \to & \FS(X), \quad (a,b,t) ~
    \mapsto (\overline{f}(a),\overline{f}(b),
             t + \langle a,b|x_0,f^{-1}(x_0)\rangle.
\end{eqnarray*}
(Already $d^{\times}_{\FS,x_0}$ is invariant under this map by a
straight-forward calculation that
uses Lemma~\ref{lem:properties_of_theta_[alpha,beta]}~
\ref{lem:properties_of_theta_[alpha,beta]:shift_invariance}.)
We have $\FS(g \circ f) = \FS(g) \circ \FS(f)$ and
$\FS(\id) = \id$.
In particular a $G$-action on $X$ by isometries with
respect to the word metric extends to a $G$-action on
$\FS(X)$ by isometries with respect to
the metric $d_{\FS(X),x_0}$.

Next we compute the pseudometric $d^{\times}_{\FS}$ and the
metric $d_{\FS}$ on a line.

\begin{lemma}\label{lem:d^times_{x_0}((a,b,t),(a,b,s)right)}
We get for $(a,b,t), (a,b,s) \in \FS(X)$
and a given base point $x_0 \in X$
\begin{eqnarray*}
d^{\times}_{\FS,x_0}\left((a,b,t),(a,b,s)\right)
& = &
\left|\theta_{[-\langle b|x_0\rangle_{a},\langle a|x_0\rangle_b]}(t)
- \theta_{[-\langle b|x_0\rangle_{a},\langle a|x_0\rangle_b]}(s)\right|;
\\
d_{\FS,x_0}\left((a,b,t),(a,b,s)\right)
& = &
\left|\Theta_{[-\langle b|x_0\rangle_{a},\langle a|x_0\rangle_b]}(t)
- \Theta_{[-\langle b|x_0\rangle_{a},\langle a|x_0\rangle_b]}(s)\right|.
\end{eqnarray*}
\end{lemma}

\begin{proof}
We have
\begin{eqnarray*}
\lefteqn{d^{\times}_{\FS,x_0}\left((a,b,t),(a,b,s)\right)}
& &
\\
& = &
\sup_{u \in X} \left|l_{x_0}(u,(a,b,t)) - l_{x_0}(u,(a,b,s))\right|
\\
& = &
\sup_{u \in X} \left|\left(\langle a|b \rangle_{u} + \left|\theta_{[-\langle b|x_0\rangle_a,\langle a|x_0\rangle_b]}(t) -
    \langle a,b|u,x_0\rangle \right|\right) \right.
\\
& & \hspace{40mm}
 - \left.\left(\langle a|b \rangle_{u} + \left|\theta_{[-\langle b|x_0\rangle_a,\langle a|x_0\rangle_b]}(s) -
    \langle a,b|u,x_0\rangle \right|\right)\right|
\\
& = &
\sup_{u \in X} \left|\left|\theta_{[-\langle b|x_0\rangle_a,\langle a|x_0\rangle_b]}(t) -
    \langle a,b|u,x_0\rangle \right|\right.
- \left.\left|\theta_{[-\langle b|x_0\rangle_a,\langle a|x_0\rangle_b]}(s) -
    \langle a,b|u,x_0\rangle \right|\right|
\end{eqnarray*}
We conclude from the triangle inequality
\begin{eqnarray*}
\lefteqn{\left|\left|\theta_{[-\langle b|x_0\rangle_a,\langle a|x_0\rangle_b]}(t) -
    \langle a,b|u,x_0\rangle \right|
- \left|\theta_{[-\langle b|x_0\rangle_a,\langle a|x_0\rangle_b]}(s) -
    \langle a,b|u,x_0\rangle \right|\right|}
\\
& \le &
\left|\theta_{[-\langle b|x_0\rangle_a,\langle a|x_0\rangle_b]}(t) -
    \langle a,b|u,x_0\rangle
    - \left(\theta_{[-\langle b|x_0\rangle_a,\langle a|x_0\rangle_b]}(s) -
    \langle a,b|u,x_0\rangle \right)\right|
\\
& \le &
\left|\theta_{[-\langle b|x_0\rangle_a,\langle a|x_0\rangle_b]}(t)
    - \theta_{[-\langle b|x_0\rangle_a,\langle a|x_0\rangle_b]}(s) \right|,
\end{eqnarray*}
and we get for $u = a$
\begin{eqnarray*}
\lefteqn{\left|\left|\theta_{[-\langle b|x_0\rangle_a,\langle a|x_0\rangle_b]}(t) -
    \langle a,b|a,x_0\rangle \right|
- \left|\theta_{[-\langle b|x_0\rangle_a,\langle a|x_0\rangle_b]}(s) -
    \langle a,b|a,x_0\rangle \right|\right|}
\\
& = &
\left|\left|\theta_{[-\langle b|x_0\rangle_a,\langle a|x_0\rangle_b]}(t) -
    (-\langle b|x_0\rangle_a ) \right|
    - \left|\theta_{[-\langle b|x_0\rangle_a,\langle a|x_0\rangle_b]}(s) -
    (-\langle b|x_0\rangle_a ) \right|\right|
\\
& = &
\left|\theta_{[-\langle b|x_0\rangle_a,\langle a|x_0\rangle_b]}(t)
    - \theta_{[-\langle b|x_0\rangle_a,\langle a|x_0\rangle_b]}(s) \right|.
\end{eqnarray*}
This implies
\begin{eqnarray*}
d^{\times}_{\FS,x_0}\left((a,b,t),(a,b,s)\right)
& = &
\left|\theta_{[-\langle b|x_0\rangle_{a},\langle a|x_0\rangle_b]}(t)
- \theta_{[-\langle b|x_0\rangle_{a},\langle a|x_0\rangle_b]}(s)\right|.
\end{eqnarray*}

We prove the claim for $d_{\FS,x_0}$ only in the case $t \ge s$, the case $t \le s$ is analogous.
Then $t + \tau \ge s + \tau $ holds for all $\tau \in \IR$.
Since both $\theta_{[-\langle b|x_0\rangle_a,\langle a|x_0\rangle_b]}$ and
$\Theta_{[-\langle b|x_0\rangle_a,\langle a|x_0\rangle_b]}$ are monotone increasing,
we conclude for all $\tau \in \IR$
\begin{eqnarray*}
\lefteqn{\theta_{[-\langle b|x_0\rangle_a,\langle a|x_0\rangle_b]}(t+ \tau)
    - \theta_{[-\langle b|x_0\rangle_a,\langle a|x_0\rangle_b]}(s+\tau)}
\\
& \hspace{20mm} = &
\left|\theta_{[-\langle b|x_0\rangle_a,\langle a|x_0\rangle_b]}(t+ \tau)
    - \theta_{[-\langle b|x_0\rangle_a,\langle a|x_0\rangle_b]}(s+\tau)\right|;
\\
\lefteqn{\Theta_{[-\langle b|x_0\rangle_{a},\langle a|x_0\rangle_b]}(t+\tau)
- \Theta_{[-\langle b|x_0\rangle_{a},\langle a|x_0\rangle_b]}(s + \tau)}
\\
& \hspace{20mm} = &
\left|\Theta_{[-\langle b|x_0\rangle_{a},\langle a|x_0\rangle_b]}(t + \tau)
- \Theta_{[-\langle b|x_0\rangle_{a},\langle a|x_0\rangle_b]}(s + \tau)\right|.
\end{eqnarray*}
Now we get
\begin{eqnarray*}
\lefteqn{d_{\FS,x_0}\left((a,b,t),(a,b,s)\right)}
\\
& = &
\int_{\IR} \frac{d^{\times}_{\FS,x_0}(\phi_{\tau}((a,b,t)),\phi_{\tau}(a,b,s))}{2 \cdot e^{|\tau|}} d\tau
\\ & = &
\int_{\IR} \frac{d^{\times}_{\FS,x_0}((a,b,t+ \tau)),(a,b,s+\tau))}{2 \cdot e^{|\tau|}} d\tau
\\
& = &
\int_{\IR} \frac{\left|\theta_{[-\langle b|x_0\rangle_a,\langle a|x_0\rangle_b]}(t+ \tau)
    - \theta_{[-\langle b|x_0\rangle_a,\langle a|x_0\rangle_b]}(s+\tau)\right|}{2 \cdot e^{|\tau|}} d\tau
\\
& = &
\int_{\IR} \frac{\theta_{[-\langle b|x_0\rangle_a,\langle a|x_0\rangle_b]}(t+ \tau)
    - \theta_{[-\langle b|x_0\rangle_a,\langle a|x_0\rangle_b]}(s+\tau)}{2 \cdot e^{|\tau|}} d\tau
\\
& = &
\int_{\IR} \frac{\theta_{[-\langle b|x_0\rangle_a,\langle a|x_0\rangle_b]}(t+ \tau) }{2 \cdot e^{|\tau|}} d\tau
~ - ~
\int_{\IR} \frac{\theta_{[-\langle b|x_0\rangle_a,\langle a|x_0\rangle_b]}(s+\tau)}{2 \cdot e^{|\tau|}} d\tau
\\
& = &
\Theta_{[-\langle b|x_0\rangle_{a},\langle a|x_0\rangle_b]}(t)
- \Theta_{[-\langle b|x_0\rangle_{a},\langle a|x_0\rangle_b]}(s)
\\
& = &
\left|\Theta_{[-\langle b|x_0\rangle_{a},\langle a|x_0\rangle_b]}(t)
- \Theta_{[-\langle b|x_0\rangle_{a},\langle a|x_0\rangle_b]}(s)\right|.
\end{eqnarray*}
This finishes the proof of Lemma~\ref{lem:d^times_{x_0}((a,b,t),(a,b,s)right)}.
\end{proof}

\begin{remark} \label{rem:base_point_free_formulation}
We have fixed a base point $x_0 \in X$.
For a different base point $x_1 \in X$ there is a
canonical isometry
$(\FS(X),d_{\FS,x_0}) \to (\FS(X),d_{\FS,x_1})$ defined by
$(a,b,t) \mapsto (a,b,t+\langle a,b | x_0,x_1 \rangle)$.
(Of course  this isometry appeared already when we defined
$\FS(f)$ for an isometry $f \colon X \to X$.)
Using these isometries and a colimit over all choices
of base points it is possible to give a canonical construction
of the metric space $\FS(X)$ without choosing a base point.
However, then we do no longer have canonical coordinates
in $\FS(X)$, i.e.\ to make sense out of $(a,b,t) \in \FS(X)$
we would still need to pick a base point.
Since the base-point free formulation is not directly relevant for our applications,
we do not give any details here.
\end{remark}

\begin{remark}
As pointed out before, $\FS(X)$ and $\sjX$ agree
as topological spaces.
But it should be noted that
the construction of the metric $d_{\FS,x_0}$ on $\FS(X)$
differs slightly from the metric $d_*$ constructed by Mineyev
on $\sjX$.
The definitions of $l_{x_0}$ in \eqref{l_{x_0}} and
of $l(u,x)$ in
\cite[Definition~10 on page 422]{Mineyev-flows-and-joins}
do not quite agree.
There Mineyev uses a different parametrization to
the effect that his formula translates to
replacing $\theta$ by $\Theta$ in \eqref{l_{x_0}}.
(The point $[\![ a,b;s ]\!]_{x_0}  \in \sjX$
corresponds to
$(a,b,\Theta_{[-\langle b|x_0\rangle_a,
                 \langle a|x_0\rangle_b]}(s) )$
in our parametrization, because
$[\![ a,b;s ]\!]'_{x_0} =
 [\![ a,b;\Theta_{[-\langle b|x_0\rangle_a,
                    \langle a|x_0\rangle_b]}(s) ]\!]_{x_0}$
is used in \cite[Section 2.3]{Mineyev-flows-and-joins}
to identify the models $\sjnobarX$ and $\sjnobarXx$.)
However, this difference is not important.
All results of \cite{Mineyev-flows-and-joins}
that we will use are also valid with this minor variation.
Moreover we remark that since $|\theta_{[-\langle b|x_0\rangle_a,
                    \langle a|x_0\rangle_b]}  (t)  -
             \Theta_{[-\langle b|x_0\rangle_a,
                    \langle a|x_0\rangle_b]}(t)| \le 1$
for all $t \in \IR$, it is easy to check that
the identification of $\FS(X)$ and $\osjX$
is a quasi-isometry with respect to
$d_{\FS,x_0}$ and $d_*$.
\end{remark}


\section{Flow estimates}
\label{sec:flow_estimates}

In this section we prove the main exponential flow estimate
for $\FS(X)$.
Recall that we have fixed a base point $x_0 \in X$.

\begin{theorem}[Exponential Flow estimate]
\label{the:exponential_estimate_for_d_{overline{FS}_x_0}}
Let $\lambda \in (e^{-1},1)$ and $T \in [0,\infty)$ be the constants depending only on $X$ which appear in
Proposition~\ref{prop:Mineyevs_Proposition_38}. Consider $a,b,c \in \overline{X}$ such that
$a,c \in \partial X \Rightarrow a \not= c$ and $b,c \in \partial X \Rightarrow b \not= c$ holds.
Let $t,s, \tau \in \IR$. Put
\begin{eqnarray*}
\tau_0  & = & t - s -  \langle a,b|c,x_0 \rangle;
\\
N & = & 2 + \frac{2}{\lambda^T \cdot (-\ln(\lambda))}.
\end{eqnarray*}

Then we get
$$d_{\FS,x_0}(\phi_{\tau}(a,c,t),\phi_{\tau + \tau_0}(b,c,s)) ~ \le ~
\frac{N}{1 - \ln(\lambda)^2} \cdot \lambda^{\left(t - \langle a,c|b,x_0\rangle\right)} \cdot \lambda^{\tau}.$$
\end{theorem}

For the sphere bundle of the universal cover of
a strictly negatively curved manifold estimates as above
are classical results and have been used in algebraic
$K$-theory in \cite{Farrell-Jones-dynamics-I}.
Compare also
\cite[Proposition~14.2]{Bartels-Farrell-Jones-Reich(2004)}.
There only $c \in \dd X$ is considered and $\tau_0$ is chosen to ensure that
$\phi_\tau(a,c,t)$ and $\phi_{\tau+\tau_0}(b,c,s)$
both lie on the same horosphere around $c$.

As mentioned in the introduction, the proof of
Theorem~\ref{the:exponential_estimate_for_d_{overline{FS}_x_0}}
is strongly based on ideas due to
Mineyev~\cite[Theorem~57 on p.~468]{Mineyev-flows-and-joins}.

We also will use the following basic flow estimate.

\begin{lemma} \label{lem:_exponential_estimate_on_flow}
We get for
$(a,b,t), (c,d,s) \in \FS(X)$ and $\tau \in \IR$
$$d_{\FS,x_0}(\phi_{\tau}(a,b,t),\phi_{\tau}(c,d,s)) ~ \le ~ e^{|\tau|} \cdot d_{\FS,x_0}((a,b,t),(c,d,s)).$$
\end{lemma}
\begin{proof}
We compute
\begin{eqnarray*}
\lefteqn{d_{\FS,x_0}(\phi_{\tau}(a,b,t),\phi_{\tau}(c,d,s))}
\\
& = & \int_{\IR}\frac{d_{\FS,x_0}^{\times}(\phi_{\sigma}(\phi_{\tau}(a,b,t)),\phi_{\sigma}(\phi_{\tau}(c,d,s)))}{2 \cdot e^{|\sigma|}}~ d\sigma
\\
& = & \int_{\IR}\frac{d_{\FS,x_0}^{\times}(\phi_{\sigma + \tau}(a,b,t),\phi_{\sigma + \tau}(c,d,s))}{2 \cdot e^{|\sigma|}}~ d\sigma
\\
& = & \int_{\IR}\frac{d_{\FS,x_0}^{\times}(\phi_{\sigma}(a,b,t),\phi_{\sigma}(c,d,s))}{2 \cdot e^{|\sigma-\tau|}}~ d\sigma
\\
& \le & \int_{\IR}\frac{d_{\FS,x_0}^{\times}(\phi_{\sigma}(a,b,t),\phi_{\sigma}(c,d,s))}{2 \cdot e^{|\sigma| -|\tau|}}~ d\sigma
\\
& \le & e^{|\tau|} \cdot \int_{\IR}\frac{d_{\FS,x_0}^{\times}(\phi_{\sigma}(a,b,t),\phi_{\sigma}(c,d,s))}{2 \cdot e^{|\sigma|}}~ d\sigma
\\
& = & e^{|\tau|} \cdot d_{\FS,x_0}((a,b,t),(c,d,s)).
\end{eqnarray*}
\end{proof}

We record the following result due to Mineyev\cite[Proposition 48 on page~460]{Mineyev-flows-and-joins}.

\begin{theorem} \label{the:topology of FS(X) - X}
The map
$$\left(\overline{X} \times \overline{X} - \Delta(\overline{X}) \right) \times \IR ~ \xrightarrow{\cong} \FS(X) - X,
\quad ((a,b),t) \mapsto (a,b,t)$$
is a homeomorphism, where $\Delta(\overline{X}) \subseteq \overline{X} \times \overline{X}$ is the diagonal.
\end{theorem}


\subsection{Flow estimates for the pseudo metric}
\label{subsec:flow_estimates_for_the_pseudo_metric}

The goal of this subsection is to prove the version of Theorem~\ref{the:exponential_estimate_for_d_{overline{FS}_x_0}}
for the pseudometric $d^{\times}_{\FS,x_0}$.

\begin{theorem} \label{the:exponential_estimate_for_d^times}
Let $\lambda \in (e^{-1},1)$ and $T \in [0,\infty)$ be the constants depending only on $X$ which appear in
Proposition~\ref{prop:Mineyevs_Proposition_38}. Consider $a,b,c \in \overline{X}$ such that
$a,c \in \partial X \Rightarrow a \not= c$ and $b,c \in \partial X \Rightarrow b \not= c$ holds.
Let $t,s, \tau \in \IR$. Put
\begin{eqnarray*}
\tau_0  & = & t - s -  \langle a,b|c,x_0 \rangle;
\\
N & = & 2 + \frac{2}{\lambda^T \cdot (-\ln(\lambda))}.
\end{eqnarray*}
Then we get
\begin{equation}
\label{eq:exponential-flow-estimate-for-d^x}
d^{\times}_{\FS,x_0}(\phi_{\tau}(a,c,t),\phi_{\tau + \tau_0}(b,c,s)) ~ \le ~
N \cdot \lambda^{t+\tau  - \langle a,c|b,x_0\rangle}.
\end{equation}
\end{theorem}

Its proof  needs some preparations and is then done in several steps.

We begin with the trivial case $a = b$.

\begin{lemma}
\label{lem:exponential_estimate_for_d^times:a=b_in_partial_X}
Consider the situation appearing in
Theorem~\ref{the:exponential_estimate_for_d^times}.
If $ a = b$
then \eqref{eq:exponential-flow-estimate-for-d^x} holds.
\end{lemma}
\begin{proof}
Since
$$\tau_0 ~ = ~ t - s - \langle a,b|c,x_0\rangle ~ = ~ t - s - \langle a,a|c,x_0\rangle ~ = ~ t-s,$$
we get for all $\tau \in \IR$
$$\phi_{\tau}(a,c,t) ~ = ~ (a,c,t + \tau) = ~ (b,c, s + \tau + \tau_0)
~ = ~ \phi_{\tau + \tau_0}(b,c,s)$$
and hence
$$d^{\times}_{\FS,x_0}(\phi_{\tau}(a,c,t),\phi_{\tau +
  \tau_0}(b,c,s)) ~ = ~ 0.$$
\end{proof}

So we can make in the sequel the additional assumption that
$a \not = b$. This has the advantage
that the expressions
$\langle a,c|b,x_0\rangle$, $\langle b,c|a,x_0\rangle$, $\langle
a,c|u,b\rangle$ and $\langle b,c|u,a\rangle$
for $u \in X$ which will appear below are well defined elements in $\IR$.

\begin{lemma} \label{lem:relating_t+tau_and_s+tau+tau_0}
Define $\tau, \tau_0,a,b,c, x_0$ as in
Theorem~\ref{the:exponential_estimate_for_d^times}.
Suppose that $a,b \in \partial X \Rightarrow a \not= b$.
Then

\begin{enumerate}
\item \label{lem:relating_t+tau_and_s+tau+tau_0:(1)}
We have
\begin{eqnarray*}
\langle a|x_0\rangle_c - \langle a,c|b,x_0\rangle
& = &
\langle b|x_0\rangle_c - \langle b,c|a,x_0\rangle ~ = ~ \langle a|b \rangle_c;
\\
t + \tau - \langle a,c|b,x_0\rangle
& = &
s + \tau + \tau_0 - \langle b,c|a,x_0\rangle ;
\\
\langle a|x_0\rangle_c - (t + \tau)
& = &
\langle b|x_0 \rangle_c - (s + \tau + \tau_0);
\\
-\langle c|x_0\rangle_a - \langle a,c|b,x_0\rangle & = & -\langle c|b\rangle_a;
\\
-\langle c|x_0\rangle_b - \langle b,c|a,x_0\rangle & = & -\langle c|a\rangle_b;
\\
\max\{\langle a,c|u,b\rangle, \langle b,c|u,a\rangle\}
& \le  &
\langle a|x_0 \rangle_c -  \langle a,c|b,x_0\rangle \quad \text { for } u \in X;
\end{eqnarray*}

\item \label{lem:relating_t+tau_and_s+tau+tau_0:(2)}
If
\begin{eqnarray*}
t + \tau & \ge & - \langle c|x_0\rangle_a;
\\
s + \tau + \tau_0 & \ge & - \langle c|x_0\rangle_b,
\end{eqnarray*}
then
\begin{equation*}
\theta_{[-\langle c|x_0\rangle_a,\langle a|x_0\rangle_c]}(t+\tau)  - \langle a,c|b,x_0\rangle
 = 
\theta_{[-\langle c|x_0\rangle_b,\langle b|x_0\rangle_c]}(s+\tau_0 + \tau)  - \langle b,c|a,x_0\rangle;
\end{equation*}

\item \label{lem:relating_t+tau_and_s+tau+tau_0:(3)}
If we assume
\begin{eqnarray*}
t + \tau & \ge & - \langle c|x_0\rangle_a;
\\
\theta_{[-\langle c|x_0\rangle_a,\langle a|x_0\rangle_c]}(t+\tau)  - \langle a,c|b,x_0\rangle
& < &
\max\{\langle a,c|u,b\rangle, \langle b,c|u,a\rangle\},
\end{eqnarray*}
then
\begin{eqnarray*}
\theta_{[-\langle c|x_0\rangle_a,\langle a|x_0\rangle_c]}(t+\tau) & = & t + \tau.
\end{eqnarray*}

\end{enumerate}
\end{lemma}
\begin{proof}
\ref{lem:relating_t+tau_and_s+tau+tau_0:(1)}
Notice for the sequel that $\langle a|a'\rangle_c = \infty$ if $c \in
\partial X$.
Hence the first equality, the third equality and the last inequality are true
for trivial reasons if $c \in \partial X$.

We compute for $c \in X$
\begin{eqnarray*}
\langle a|x_0\rangle_c - \langle a,c|b,x_0\rangle
& = & \langle a|x_0\rangle_c + \langle c,a|b,x_0\rangle
\\
& = &
\langle a|x_0\rangle_c + \left(\langle a|b \rangle_{c} - \langle a|x_0 \rangle_{c}\right)
\\
& = &
\langle a|b \rangle_{c},
\end{eqnarray*}
and analogously
$$
\langle b|x_0\rangle_c - \langle b,c|a,x_0\rangle
~ = ~
\langle b|a \rangle_{c}
~ = ~
\langle a|b\rangle_c.
$$
This proves the first equation.

The second follows from
\begin{eqnarray*}
\lefteqn{s + \tau + \tau_0 - \langle b,c|a,x_0\rangle}
\\
& = &
s + \tau + t - s - \langle a,b|c,x_0\rangle - \langle b,c|a,x_0\rangle
\\
& = &
\tau + t - \langle a,b|c,x_0\rangle - \langle b,c|a,x_0\rangle -
\langle c,a|b,x_0\rangle + \langle c,a|b,x_0\rangle
\\
& = &
\tau + t  - \langle a,c|b,x_0\rangle.
\end{eqnarray*}

The third equation is a direct consequence of the first two if $c \in
X$, and hence true for all $c \in \overline{X}$.

The proof of the fourth and fifth equation is analogous to the one
of the first one.

Since for $c \in X$
$$ \langle a,c|u,b\rangle ~ = ~\langle c,a|b,u\rangle ~ = ~ \langle a|b\rangle_c - \langle a|u \rangle_c
~ \le ~ \langle a|b\rangle_c $$
and
$$
\langle b,c|u,a\rangle ~ = ~ \langle c,b|a,u\rangle ~ = ~ \langle b|a\rangle_c - \langle b|u \rangle_c
~ \le ~ \langle b|a\rangle_c ~ = ~ \langle a|b\rangle_c
$$
holds, the last inequality  follows from the first equality.
\\[1mm]
\ref{lem:relating_t+tau_and_s+tau+tau_0:(2)}
We begin with the case
$t + \tau \ge \langle a|x_0\rangle_c$. Then we get from the third equation of
assertion~\ref{lem:relating_t+tau_and_s+tau+tau_0:(1)}
that also
$s + \tau + \tau_0 \ge \langle b|x_0\rangle_c$ holds. We conclude from the definitions
and the first equation of assertion~\ref{lem:relating_t+tau_and_s+tau+tau_0:(1)}
\begin{eqnarray*}
\theta_{[-\langle c|x_0\rangle_a,\langle a|x_0\rangle_c]}(t+\tau)  - \langle a,c|b,x_0\rangle
& = &
\langle a|x_0\rangle_c  - \langle a,c|b,x_0\rangle
\\
& = &  \langle b|x_0\rangle_c - \langle b,c|a,x_0\rangle
\\
& = &
\theta_{[-\langle c|x_0\rangle_b,\langle b|x_0\rangle_c]}(s+\tau_0 + \tau)  - \langle b,c|a,x_0\rangle.
\end{eqnarray*}
Next we treat the case $t + \tau \le \langle a|x_0\rangle_c$. Then we get from the third equation
of assertion~\ref{lem:relating_t+tau_and_s+tau+tau_0:(1)} that also
$s + \tau + \tau_0 \le \langle b|x_0\rangle_c$ holds. Since by assumption
$t + \tau  \ge  - \langle c|x_0\rangle_a$ and $s + \tau + \tau_0  \ge  -\langle c|x_0\rangle_b$ holds,
we conclude from the definitions
and the second equation of assertion~\ref{lem:relating_t+tau_and_s+tau+tau_0:(1)}
\begin{eqnarray*}\theta_{[-\langle c|x_0\rangle_a,\langle a|x_0\rangle_c]}(t+\tau)  - \langle a,c|b,x_0\rangle
& = & t + \tau  - \langle a,c|b,x_0\rangle
\\
& = & s + \tau + \tau_0  - \langle b,c|a,x_0\rangle
\\
& = &
\theta_{[-\langle c|x_0\rangle_b,\langle b|x_0\rangle_c]}(s+\tau_0 + \tau)  - \langle b,c|a,x_0\rangle.
\end{eqnarray*}
\\
\ref{lem:relating_t+tau_and_s+tau+tau_0:(3)}
Since the inequality in assertion~\ref{lem:relating_t+tau_and_s+tau+tau_0:(1)} implies
\begin{eqnarray*}
\theta_{[-\langle c|x_0\rangle_a,\langle a|x_0\rangle_c]}(t+\tau)
& <  &
\max\{\langle a,c|u,b\rangle, \langle b,c|u,a\rangle\} +  \langle a,c|b,x_0\rangle
\\
& \le &
\langle a|x_0 \rangle_c,
\end{eqnarray*}
we get $t+\tau < \langle a|x_0 \rangle_c$. Since we assume $-\langle c|x_0\rangle_{a} \le t + \tau$, we conclude
\begin{eqnarray*}
\theta_{[-\langle c|x_0\rangle_a,\langle a|x_0\rangle_c]}(t+\tau)
& = & t + \tau.
\end{eqnarray*}

This finishes the proof of Lemma~\ref{lem:relating_t+tau_and_s+tau+tau_0}.
\end{proof}

The elementary proof of the next lemma is left to the reader.

\begin{lemma}\label{lem:Choice_of_N}
Consider numbers $\lambda \in (e^{-1},1)$ and $T \in [0,\infty)$. Put
$$N ~ := ~ 2 + \frac{2}{\lambda^T \cdot (-\ln(\lambda))}.$$
Then we get
$$\begin{array}{lclcl}
2 & \le & N; & &
\\
2 \cdot (T - v + \lambda^v) & \le & N \cdot \lambda^v & & \text{ for all } v \le T;
\\
v + N \cdot \lambda^{v - w} & \le & N \cdot \lambda^{-w} & &  \text{ for } 0 \le v \le w.
\end{array}
$$
\end{lemma}

\begin{lemma} \label{lem:exponential_estimate_for_d^times:prel_computation}
Consider the situation appearing in
Theorem~\ref{the:exponential_estimate_for_d^times}.
Suppose that $a,b \in \partial X \Rightarrow a \not= b$.
Suppose that
\begin{eqnarray*}
t + \tau & \ge & - \langle c|x_0\rangle_a;
\\
s + \tau + \tau_0 & \ge & - \langle c|x_0\rangle_b.
\end{eqnarray*}
Then we get for all $u \in X$
\begin{eqnarray*}
\lefteqn{\left|l_{x_0}(u,\phi_{\tau}(a,c,t)) - l_{x_0}(u,\phi_{\tau + \tau_0}(b,c,s))\right|}
\\
\nonumber
& = &
\left|\langle u,c|a,b\rangle + \left|\theta_{[-\langle c|x_0\rangle_a,\langle a|x_0\rangle_c]}(t+\tau)  - \langle a,c|b,x_0\rangle
- \langle a,c|u,b \rangle \right|\right.
\\\nonumber
& &
\hspace{10mm} - \left.\left|\theta_{[-\langle c|x_0\rangle_a,\langle a,x_0\rangle_c]}(t+\tau)  - \langle a,c|b,x_0\rangle
- \langle b,c|u,a \rangle \right|\right|.
\end{eqnarray*}
\end{lemma}

\begin{proof}
We compute
\begin{eqnarray*}
\lefteqn{\left|l_{x_0}(u,\phi_{\tau}(a,c,t)) - l_{x_0}(u,\phi_{\tau + \tau_0}(b,c,s))\right|}
\\
& = &
\left|l_{x_0}(u,(a,c,t+\tau)) - l_{x_0}(u,(b,c,s+\tau+\tau_0))\right|
\\
& = &
\left|\langle a|c\rangle_u + \left|\theta_{[-\langle c|x_0\rangle_a,\langle a|x_0\rangle_c]}(t+\tau) - \langle a,c|u,x_0 \rangle \right|
-\langle b|c\rangle_u \right.
\\
& &
\hspace{5mm} - \left.\left|\theta_{[-\langle c|x_0\rangle_b,\langle b|x_0\rangle_c]}(s+\tau_0 + \tau) - \langle b,c|u,x_0 \rangle \right|\right|
\\
& = &
\left|\langle c|a\rangle_u -\langle c|b\rangle_u \right.
\\
& &
\hspace{5mm}
+ \left. \left|\theta_{[-\langle c|x_0\rangle_a,\langle a|x_0\rangle_c]}(t+\tau)  - \langle a,c|b,x_0\rangle +  \langle a,c|b,x_0\rangle
- \langle a,c|u,x_0 \rangle \right|
 \right.
\\
& &
\hspace{5mm} - \left.\left|\theta_{[-\langle c|x_0\rangle_b,\langle b|x_0\rangle_c]}(s+\tau_0 + \tau) - \langle b,c|a,x_0\rangle + \langle b,c|a,x_0\rangle
-\langle b,c|u,x_0 \rangle \right|\right|
\\
& = &
\left|\langle u,c|a,b\rangle + \left|\theta_{[-\langle c|x_0\rangle_a,\langle a|x_0\rangle_c]}(t+\tau)  - \langle a,c|b,x_0\rangle
- \langle a,c|u,b \rangle \right|\right.
\\
& &
\hspace{5mm} - \left.\left|\theta_{[-\langle c|x_0\rangle_b,\langle b|x_0\rangle_c]}(s+\tau_0 + \tau) - \langle b,c|a,x_0\rangle
- \langle b,c|u,a \rangle \right|\right|
\end{eqnarray*}
Since we get from Lemma~\ref{lem:relating_t+tau_and_s+tau+tau_0}~\ref{lem:relating_t+tau_and_s+tau+tau_0:(2)}
$$
\theta_{[-\langle c|x_0\rangle_a,\langle a|x_0\rangle_c]}(t+\tau)  - \langle a,c|b,x_0\rangle
~ = ~
\theta_{[-\langle c|x_0\rangle_b,\langle b|x_0\rangle_c]}(s+\tau_0 + \tau)  - \langle b,c|a,x_0\rangle.
$$
Lemma~\ref{lem:exponential_estimate_for_d^times:prel_computation} follows.
\end{proof}

\begin{lemma} \label{lem:exponential_estimate_for_d^times:case_with_d^{times}=0}
Consider the situation appearing in
Theorem~\ref{the:exponential_estimate_for_d^times}.
Suppose that $a,b \in \partial X \Rightarrow a \not= b$.
Suppose that
\begin{eqnarray*}
t + \tau & \ge & - \langle c|x_0\rangle_a;
\\
s + \tau + \tau_0 & \ge & - \langle c|x_0\rangle_b,
\end{eqnarray*}
and
$$\theta_{[-\langle c|x_0\rangle_a,\langle a|x_0\rangle_c]}(t+\tau)  - \langle a,c|b,x_0\rangle
~ \ge ~
\max\{\langle a,c|u,b\rangle, \langle b,c|u,a\rangle\}.$$
Then we get for all $u \in X$
$$\left|l_{x_0}(u,\phi_{\tau}(a,c,t)) - l_{x_0}(u,\phi_{\tau + \tau_0}(b,c,s))\right| ~ = ~ 0.$$
\end{lemma}

\begin{proof}
We compute
\begin{eqnarray*}
\lefteqn{\left|\theta_{[-\langle c|x_0\rangle_a,\langle a|x_0\rangle_c]}(t+\tau)  - \langle a,c|b,x_0\rangle
- \langle a,c|u,b \rangle \right|}
\\
& & \hspace{30mm}
- \left|\theta_{[-\langle c|x_0\rangle_a,\langle a|x_0\rangle_c]}(t+\tau)  - \langle a,c|b,x_0\rangle
- \langle b,c|u,a \rangle \right|
\\
& = &
\left(\theta_{[-\langle c|x_0\rangle_a,\langle a|x_0\rangle_c]}(t+\tau)  - \langle a,c|b,x_0\rangle
- \langle a,c|u,b \rangle \right)
\\
& & \hspace{30mm}
- \left(\theta_{[-\langle c|x_0\rangle_a,\langle a|x_0\rangle_c]}(t+\tau)  - \langle a,c|b,x_0\rangle
- \langle b,c|u,a \rangle \right)
\\
& = & - \langle a,c|u,b \rangle + \langle b,c|u,a \rangle.
\end{eqnarray*}
This implies together with Lemma~\ref{lem:exponential_estimate_for_d^times:prel_computation}
\begin{eqnarray*}
\lefteqn{\left|l_{x_0}(u,\phi_{\tau}(a,c,t)) - l_{x_0}(u,\phi_{\tau + \tau_0}(b,c,s))\right|}
\\
& = &
\left|\langle u,c|a,b\rangle - \langle a,c|u,b \rangle + \langle b,c|u,a \rangle\right|
\\
& = &
\left|\langle b,a|c,u\rangle + \langle a,c|b,u \rangle + \langle c,b|a,u \rangle\right|
\\
& = & 0.
\end{eqnarray*}
\end{proof}

\begin{lemma} \label{lem:exponential_estimate_for_d^times:other_case_and_ge_T}
Consider the situation appearing in
Theorem~\ref{the:exponential_estimate_for_d^times}.
Suppose that $a,b \in \partial X \Rightarrow a \not= b$.
Suppose that
\begin{eqnarray*}
t + \tau & \ge & - \langle c|x_0\rangle_a;
\\
s + \tau + \tau_0 & \ge & - \langle c|x_0\rangle_b,
\end{eqnarray*}
and that
\begin{eqnarray*}
\theta_{[-\langle c|x_0\rangle_a,\langle a|x_0\rangle_c]}(t+\tau)  - \langle a,c|b,x_0\rangle
& < &
\max\{\langle a,c|u,b\rangle, \langle b,c|u,a\rangle\};
\\
\theta_{[-\langle c|x_0\rangle_a,\langle a|x_0\rangle_c]}(t+\tau)  - \langle a,c|b,x_0\rangle  & \ge & T.
\end{eqnarray*}
Then we get for all $u \in X$
\begin{eqnarray*}
\left|l_{x_0}(u,\phi_{\tau}(a,c,t)) - l_{x_0}(u,\phi_{\tau + \tau_0}(b,c,s))\right|
&\le ~ &
2 \cdot \lambda^{t+\tau - \langle a,c|b,x_0\rangle}.
\end{eqnarray*}
\end{lemma}

\begin{proof}
We get from Lemma~\ref{lem:exponential_estimate_for_d^times:prel_computation}
\begin{eqnarray*}
\lefteqn{\left|l_{x_0}(u,\phi_{\tau}(a,c,t)) - l_{x_0}(u,\phi_{\tau + \tau_0}(b,c,s))\right|}
\\
& = &
\left|\langle u,c|a,b\rangle + \left|\theta_{[-\langle c|x_0\rangle_a,\langle a|x_0\rangle_c]}(t+\tau)  - \langle a,c|b,x_0\rangle
- \langle a,c|u,b \rangle \right|\right.
\\
& &
\hspace{10mm} - \left.\left|\theta_{[-\langle c|x_0\rangle_a,\langle a|x_0\rangle_c]}(t+\tau)  - \langle a,c|b,x_0\rangle
- \langle b,c|u,a \rangle \right|\right|
\\
& \le  &
\left|\langle u,c|a,b\rangle\right| + \left|\left(\theta_{[-\langle c|x_0\rangle_a,\langle a|x_0\rangle_c]}(t+\tau)  - \langle a,c|b,x_0\rangle
- \langle a,c|u,b \rangle \right)\right.
\\
& &
\hspace{10mm} - \left.\left(\theta_{[-\langle c|x_0\rangle_a,\langle a|x_0\rangle_c]}(t+\tau)  - \langle a,c|b,x_0\rangle
- \langle b,c|u,a \rangle \right)\right|
\\
& = &
\left|\langle u,c|a,b\rangle\right| + \left|- \langle a,c|u,b \rangle +  \langle b,c|u,a \rangle \right|
\\
& = &
\left|\langle u,c|a,b\rangle\right| + \left|\langle c,a|u,b \rangle + \langle a,u|c,b \rangle \right|
\\
& = &
\left|\langle u,c|a,b\rangle\right| + \left|\langle u,c |a,b\rangle\right|
\\
& = & 2 \cdot \left|\langle u,c|a,b\rangle\right|.
\end{eqnarray*}
Our assumptions imply $\max\{\langle a,c|u,b\rangle, \langle b,c|u,a\} ~ \ge ~ T$.
We conclude from Proposition~\ref{prop:Mineyevs_Proposition_38}
$$\left|\langle u,c|a,b\rangle\right| ~ \le ~ \lambda^{\max\{\langle a,c|u,b\rangle, \langle b,c|u,a\rangle\}}.$$
Thus we get using
Lemma~\ref{lem:relating_t+tau_and_s+tau+tau_0}~\ref{lem:relating_t+tau_and_s+tau+tau_0:(3)}
\begin{eqnarray*}
\left|l_{x_0}(u,\phi_{\tau}(a,c,t)) - l_{x_0}(u,\phi_{\tau + \tau_0}(b,c,s))\right|
& \le ~ &
2 \cdot \lambda^{\max\{\langle a,c|u,b\rangle, \langle b,c|u,a\rangle\}}
\\
& \le &
2 \cdot \lambda^{\theta_{[-\langle c|x_0\rangle_a,\langle a|x_0\rangle_c]}(t+\tau)  - \langle a,c|b,x_0\rangle}
\\
& = &
2 \cdot \lambda^{t+\tau - \langle a,c|b,x_0\rangle}.
\end{eqnarray*}
\end{proof}

\begin{lemma} \label{lem:exponential_estimate_for_d^times:other_case_and_le_T}
Consider the situation appearing in
Theorem~\ref{the:exponential_estimate_for_d^times}.
Suppose that $a,b \in \partial X \Rightarrow a \not= b$.
Suppose that
\begin{eqnarray*}
t + \tau & \ge & - \langle c|x_0\rangle_a;
\\
s + \tau + \tau_0 & \ge & - \langle c|x_0\rangle_b,
\end{eqnarray*}
and that
\begin{eqnarray*}
\theta_{[-\langle c|x_0\rangle_a,\langle a|x_0\rangle_c]}(t+\tau)  - \langle a,c|b,x_0\rangle
& < &
\max\{\langle a,c|u,b\rangle, \langle b,c|u,a\rangle\};
\\
\theta_{[-\langle c|x_0\rangle_a,\langle a|x_0\rangle_c]}(t+\tau)  - \langle a,c|b,x_0\rangle  & < & T.
\end{eqnarray*}
Then we get for all $u \in X$
$$\left|l_{x_0}(u,\phi_{\tau}(a,c,t)) - l_{x_0}(u,\phi_{\tau + \tau_0}(b,c,s))\right|
 ~ \le ~
N \cdot \lambda^{t+\tau  - \langle a,c|b,x_0\rangle}.$$
\end{lemma}

\begin{proof}
Since $\theta_{[-\langle c|x_0\rangle_a,\langle a|x_0\rangle_c]}$ is
monotone increasing and by  
Lemma~\ref{lem:relating_t+tau_and_s+tau+tau_0}~\ref{lem:relating_t+tau_and_s+tau+tau_0:(1)}
and by assumption 
\begin{eqnarray*}
\theta_{[-\langle c|x_0\rangle_a,\langle a|x_0\rangle_c]}(t+\tau)  - \langle a,c|b,x_0\rangle
& < & \min\left\{\max\{\langle a,c|u,b\rangle, \langle
  b,c|u,a\rangle\},T\right\};
\\
\max\{\langle a,c|u,b\rangle, \langle b,c|u,a\rangle\} & \le & \langle a|x_0\rangle_c - \langle a,c|b,x_0\rangle,
\end{eqnarray*}
holds, we can choose $\tau' \in \IR$ satisfying
\begin{eqnarray*}
\theta_{[-\langle c|x_0\rangle_a,\langle a|x_0\rangle_c]}(t+\tau')  - \langle a,c|b,x_0\rangle
& = & \min\left\{\max\{\langle a,c|u,b\rangle, \langle  b,c|u,a\rangle\},T\right\};\\
\tau & \le & \tau'.
\end{eqnarray*}
In particular we have
\begin{eqnarray*}
t + \tau' & \ge & - \langle c|x_0\rangle_a;
\\
s + \tau' + \tau_0 & \ge & - \langle c|x_0\rangle_b.
\end{eqnarray*}
Hence Lemma~\ref{lem:exponential_estimate_for_d^times:case_with_d^{times}=0} and
Lemma~\ref{lem:exponential_estimate_for_d^times:other_case_and_ge_T} imply
\begin{eqnarray}
\label{lem:exponential_estimate_for_d^times:other_case_and_le_T(1)}
\left|l_{x_0}(u,\phi_{\tau'}(a,c,t)) - l_{x_0}(u,\phi_{\tau' + \tau_0}(b,c,s))\right|
 & \le  &
2 \cdot \lambda^{t+\tau' - \langle a,c|b,x_0\rangle}
\\
& \le & 2 \cdot \lambda^{t+\tau - \langle a,c|b,x_0\rangle}.
\nonumber
\end{eqnarray}
We compute
\begin{eqnarray}
\label{lem:exponential_estimate_for_d^times:other_case_and_le_T(2)}
\\
\lefteqn{\left|l_{x_0}(u,\phi_{\tau}(a,c,t)) - l_{x_0}(u,\phi_{\tau + \tau_0}(b,c,s))\right|}
\nonumber
\\
& = &
\left|l_{x_0}(u,\phi_{\tau}(a,c,t)) - l_{x_0}(u,\phi_{\tau'}(a,c,t)) \right.
\nonumber
\\
& & \hspace{20mm} 
+ l_{x_0}(u,\phi_{\tau'}(a,c,t)) -
l_{x_0}(u,\phi_{\tau' + \tau_0}(b,c,s))
\nonumber
\\
& & \hspace{20mm}  + \left.  l_{x_0}(u,\phi_{\tau' + \tau_0}(b,c,s)) -
l_{x_0}(u,\phi_{\tau + \tau_0}(b,c,s))\right|
\nonumber
\\
& \le &
\left|l_{x_0}(u,\phi_{\tau}(a,c,t)) - l_{x_0}(u,\phi_{\tau'}(a,c,t))\right| 
\nonumber
\\
& & \hspace{20mm} 
+ \left|l_{x_0}(u,\phi_{\tau'}(a,c,t)) -
l_{x_0}(u,\phi_{\tau' + \tau_0}(b,c,s))\right|
\nonumber
\\
& & \hspace{20mm}  + \left|  l_{x_0}(u,\phi_{\tau' + \tau_0}(b,c,s)) -
l_{x_0}(u,\phi_{\tau + \tau_0}(b,c,s))\right|.
\nonumber
\end{eqnarray}
We estimate using Lemma~\ref{lem:relating_t+tau_and_s+tau+tau_0}~\ref{lem:relating_t+tau_and_s+tau+tau_0:(3)}
\begin{eqnarray}
\label{lem:exponential_estimate_for_d^times:other_case_and_le_T(3)}
\\
\lefteqn{\left|l_{x_0}(u,\phi_{\tau}(a,c,t)) - l_{x_0}(u,\phi_{\tau'}(a,c,t))\right|}
\nonumber
\\
& = &
\left|\left(\langle a|c\rangle_u  + 
\left|\theta_{[-\langle c|x_0\rangle_a,\langle a|x_0\rangle_c]}(t+\tau) -
\langle a,c|u,x_0\rangle\right|\right)\right.
\nonumber
\\ & & \hspace{30mm} - \left.\left(\langle a|c\rangle_u  + \left|\theta_{[-\langle c|x_0\rangle_a,\langle a|x_0\rangle_c]}(t+\tau') -
\langle a,c|u,x_0\rangle\right|\right)\right|
\nonumber
\\
& = &
\left|\left|\theta_{[-\langle c|x_0\rangle_a,\langle a|x_0\rangle_c]}(t+\tau) -
\langle a,c|u,x_0\rangle\right|
\nonumber  - \left|\theta_{[-\langle c|x_0\rangle_a,\langle a|x_0\rangle_c]}(t+\tau') -
\langle a,c|u,x_0\rangle\right|\right|
\nonumber
\\
& \le &
\left|\theta_{[-\langle c|x_0\rangle_a,\langle a|x_0\rangle_c]}(t+\tau) -
\theta_{[-\langle c|x_0\rangle_a,\langle a|x_0\rangle_c]}(t+\tau')\right|
\nonumber
\\
& = &
\left|\left(\theta_{[-\langle c|x_0\rangle_a,\langle a|x_0\rangle_c]}(t+\tau) - \langle a,c|b,x_0\rangle\right)
- \left(\theta_{[-\langle c|x_0\rangle_a,\langle a|x_0\rangle_c]}(t+\tau')- \langle a,c|b,x_0\rangle\right)\right|
\nonumber
\\
& = &
\left|\left(\theta_{[-\langle c|x_0\rangle_a,\langle a|x_0\rangle_c]}(t+\tau)  - \langle a,c|b,x_0\rangle\right) -
\min\left\{\max\{\langle a,c|u,b\rangle, \langle b,c|u,a\rangle\},T\right\}\right|
\nonumber
\\
& = & \min\left\{\max\{\langle a,c|u,b\rangle, \langle b,c|u,a\rangle\},T\right\} -
\left(\theta_{[-\langle c|x_0\rangle_a,\langle a|x_0\rangle_c]}(t+\tau) - \langle a,c|b,x_0\rangle\right)
\nonumber
\\
& = & \min\left\{\max\{\langle a,c|u,b\rangle, \langle b,c|u,a\rangle\},T\right\} -
\left(t+\tau - \langle a,c|b,x_0\rangle\right)
\nonumber
\\
& \le  & T - (t + \tau - \langle a,c|b,x_0\rangle).
\nonumber
\end{eqnarray}
Analogously we get using Lemma~\ref{lem:relating_t+tau_and_s+tau+tau_0}~\ref{lem:relating_t+tau_and_s+tau+tau_0:(2)}
\begin{eqnarray}
\label{lem:exponential_estimate_for_d^times:other_case_and_le_T(4)}
\\
\lefteqn{\left|l_{x_0}(u,\phi_{\tau+\tau_0}(b,c,s)) - l_{x_0}(u,\phi_{\tau'+ \tau_0}(b,c,s))\right|}
\nonumber
\\
& = &
\left|\left(\langle b|c\rangle_u  + 
\left|\theta_{[-\langle c|x_0\rangle_b,\langle b|x_0\rangle_c]}(s+\tau + \tau_0) -
\langle b,c|u,x_0\rangle\right|\right)\right.
\nonumber
\\ & & \hspace{30mm} - \left.\left(\langle b|c\rangle_u  + \left|\theta_{[-\langle c|x_0\rangle_b,\langle b|x_0\rangle_c]}(s+\tau' +\tau_0) -
\langle b,c|u,x_0\rangle\right|\right)\right|
\nonumber
\\
& = &
\left|\left|\theta_{[-\langle c|x_0\rangle_b,\langle b|x_0\rangle_c]}(s+\tau + \tau_0)
 - \langle b,c|u,x_0\rangle\right| \right.
\nonumber
\\ & & \hspace{30mm}  
- \left. 
\left|\theta_{[-\langle c|x_0\rangle_b,\langle b|x_0\rangle_c]}(s+\tau' +\tau_0) -
\langle b,c|u,x_0\rangle\right|\right|
\nonumber
\\
& \le  &
\left|\theta_{[-\langle c|x_0\rangle_b,\langle b|x_0\rangle_c]}(s+\tau + \tau_0) -
\theta_{[-\langle c|x_0\rangle_b,\langle b|x_0\rangle_c]}(s+\tau' +\tau_0)\right|
\nonumber
\\
& =  &
\left|\left(\theta_{[-\langle c|x_0\rangle_b,\langle  b|x_0\rangle_c]}
       (s+\tau + \tau_0) -
\langle b,c|a,x_0\rangle\right) \right.
\nonumber
\\ & & \hspace{30mm} 
-
\left.
\left(\theta_{[-\langle c|x_0\rangle_b,\langle
    b,x_0\rangle_c]}(s+\tau' +\tau_0) - \langle b,c|a,x_0\rangle\right)\right|
\nonumber
\\
& = &
\left|\left(\theta_{[-\langle c|x_0\rangle_a,\langle a|x_0\rangle_c]}(t+\tau) -
 \langle a,c|b,x_0\rangle\right) \right.
\nonumber
\\ & & \hspace{30mm} 
- \left. 
\left(\theta_{[-\langle c|x_0\rangle_a,\langle a|x_0\rangle_c]}(t+\tau')- \langle a,c|b,x_0\rangle\right)\right|
\nonumber
\\
& \le &  T - (t + \tau - \langle a,c|b,x_0\rangle).
\nonumber
\end{eqnarray}
Lemma~\ref{lem:relating_t+tau_and_s+tau+tau_0}~\ref{lem:relating_t+tau_and_s+tau+tau_0:(3)}
 implies
$$t+\tau  - \langle a,c|b,x_0\rangle ~ = ~ \theta_{[-\langle c|x_0\rangle_a,\langle a|x_0\rangle_c]}(t+\tau) -
\langle a,c|b,x_0\rangle ~ \le ~ T.$$
Hence we conclude from Lemma~\ref{lem:Choice_of_N} for $v = t + \tau - \langle a,c|b,x_0\rangle$
\begin{eqnarray}
\label{lem:exponential_estimate_for_d^times:other_case_and_le_T(5)}
\hspace{10mm} 2 \cdot \left(T - (t + \tau - \langle a,c|b,x_0\rangle) + \lambda^{t+\tau - \langle a,c|b,x_0\rangle}\right)
& \le &
N \cdot \lambda^{t+\tau  - \langle a,c|b,x_0\rangle}.
\end{eqnarray}

Combining \eqref{lem:exponential_estimate_for_d^times:other_case_and_le_T(1)},
 \eqref{lem:exponential_estimate_for_d^times:other_case_and_le_T(2)},
 \eqref{lem:exponential_estimate_for_d^times:other_case_and_le_T(3)},
 \eqref{lem:exponential_estimate_for_d^times:other_case_and_le_T(4)} and
 \eqref{lem:exponential_estimate_for_d^times:other_case_and_le_T(5)}
yields
\begin{eqnarray*}
\lefteqn{\left|l_{x_0}(u,\phi_{\tau}(a,c,t)) - l_{x_0}(u,\phi_{\tau + \tau_0}(b,c,s))\right|}
\nonumber
\\
& \le &
2 \cdot \left(T - (t + \tau - \langle a,c|b,x_0\rangle) + \lambda^{t+\tau - \langle a,c|b,x_0\rangle}\right)
\\
& \le &
N \cdot \lambda^{t+\tau  - \langle a,c|b,x_0\rangle}.
\end{eqnarray*}
\end{proof}

\begin{lemma} \label{lem:exponential_estimate_for_d^times:t+tau_and_s+tau_0+tau_large}
Consider the situation appearing in
Theorem~\ref{the:exponential_estimate_for_d^times}.
Suppose that $a,b
\in \partial X \Rightarrow a \not= b$. Suppose that
\begin{eqnarray*}
t + \tau & \ge & - \langle c|x_0\rangle_a;
\\
s + \tau + \tau_0 & \ge & - \langle c|x_0\rangle_b,
\end{eqnarray*}
Then \eqref{eq:exponential-flow-estimate-for-d^x} holds.
\end{lemma}

\begin{proof}
This follows from
Lemma~\ref{lem:exponential_estimate_for_d^times:case_with_d^{times}=0},
Lemma~\ref{lem:exponential_estimate_for_d^times:other_case_and_ge_T}
(note that $2 \leq N$ by Lemma~\ref{lem:Choice_of_N}) and
Lemma~\ref{lem:exponential_estimate_for_d^times:other_case_and_le_T}
since by definition
$$d^{\times}_{\FS,x_0}(\phi_{\tau}(a,c,t),\phi_{\tau + \tau_0}(b,c,s))
~ = ~ \sup_{u \in X} \left|l_{x_0}(u,\phi_{\tau}(a,c,t)) -
  l_{x_0}(u,\phi_{\tau + \tau_0}(b,c,s))\right|.$$
\end{proof}

\begin{lemma} \label{lem:exponential_estimate_for_d^times:small_t+tau}
Consider the situation appearing in
Theorem~\ref{the:exponential_estimate_for_d^times}.
Suppose that $a,b \in \partial X \Rightarrow a \not= b$. Suppose that
at least one of the following inequalities is true
\begin{eqnarray*}
t + \tau & \le & - \langle c|x_0\rangle_a;
\\
s + \tau + \tau_0 & \le & - \langle c|x_0\rangle_b.
\end{eqnarray*}
Then \eqref{eq:exponential-flow-estimate-for-d^x} holds.
\end{lemma}

\begin{proof}
Put
\begin{eqnarray*}
\tau'' & := & \max\{-\langle c|x_0\rangle_a - t, -\langle c|x_0\rangle_b - s - \tau_0\}.
\end{eqnarray*}
Since by assumption
$\tau ~ \le ~ - \langle c|x_0\rangle_a - t $ or $\tau ~ \le ~ - \langle c|x_0\rangle_b - s - \tau_0$
holds, we must have
$$\tau \le \tau''.$$
We estimate
\begin{eqnarray}
\label{lem:exponential_estimate_for_d^times:small_t+tau:inequality(1)}
\\
\lefteqn{d^{\times}_{\FS,x_0}(\phi_{\tau}(a,c,t),\phi_{\tau + \tau_0}(b,c,s))}
\nonumber
\\
& = &
d^{\times}_{\FS,x_0}((a,c,t+\tau),(b,c,s + \tau + \tau_0))
\nonumber
\\
& \le &
d^{\times}_{\FS,x_0}((a,c,t+\tau),(a,c,t+\tau'')) + d^{\times}_{\FS,x_0}((a,c,t+\tau''),(b,c,s + \tau'' + \tau_0))
\nonumber
\\
& & \hspace{30mm} + d^{\times}_{\FS,x_0}((b,c,s + \tau'' + \tau_0),(b,c,s + \tau + \tau_0)).
\nonumber
\end{eqnarray}
Since
\begin{eqnarray*}
t + \tau'' & \ge & - \langle c|x_0\rangle_a;
\\
s + \tau'' + \tau_0 & \ge & - \langle c|x_0\rangle_b,
\end{eqnarray*}
holds by definition of $\tau''$, we get from Lemma~\ref{lem:exponential_estimate_for_d^times:t+tau_and_s+tau_0+tau_large}
\begin{eqnarray}
\label{lem:exponential_estimate_for_d^times:small_t+tau:inequality(2)}
d^{\times}_{\FS,x_0}((a,c,t+\tau''),(b,c,s + \tau'' + \tau_0))
& \le & N \cdot \lambda^{t + \tau'' - \langle a,c|b,x_0\rangle}.
\end{eqnarray}
Next  we want to show
\begin{eqnarray}
\label{lem:exponential_estimate_for_d^times:small_t+tau:inequality(3)}
\lefteqn{d^{\times}_{\FS,x_0}((a,c,t+\tau),(a,c,t+\tau''))}
\\ & & \hspace{5mm}  + d^{\times}_{\FS,x_0}((b,c,s + \tau'' + \tau_0),(b,c,s + \tau + \tau_0))
 ~ \le  ~  \tau'' - \tau.
\nonumber
\end{eqnarray}
Inspecting the definition of $ \tau''$ we see that we have to consider two cases, namely,
$$
t + \tau'' = - \langle c|x_0\rangle_a  ~ \text{ and } ~
s + \tau_0 + \tau'' \ge - \langle c | x_0\rangle_b,
$$
and
$$
t + \tau'' \ge - \langle c | x_0\rangle_a  ~ \text{ and } ~
s + \tau_0 + \tau'' =  - \langle c|x_0\rangle_b.
$$
We only treat the first one, the second is completely analogous.
{}From $t + \tau ~ \le ~ t + \tau'' ~ = ~  - \langle c|x_0\rangle_a$  we conclude
$\theta_{[-\langle c|x_0\rangle_{a},\langle a|x_0\rangle_c]}(t + \tau) ~ = ~
\theta_{[-\langle c|x_0\rangle_{b},\langle b|x_0\rangle_c]}(t + \tau'')$.
Lemma~\ref{lem:d^times_{x_0}((a,b,t),(a,b,s)right)} implies
$$d^{\times}_{\FS,x_0}((a,c,t+\tau),(a,c,t+\tau'')) = 0.$$
We conclude from Lemma~\ref{lem:d^times_{x_0}((a,b,t),(a,b,s)right)}
\begin{eqnarray*}
\lefteqn{d^{\times}_{\FS,x_0}((b,c,s + \tau'' + \tau_0),(b,c,s + \tau + \tau_0))}
\\
& = & \left|\theta_{[-\langle c|x_0\rangle_{b},\langle b|x_0\rangle_c]}(s + \tau'' + \tau_0)
- \theta_{[-\langle c|x_0\rangle_{b},\langle b|x_0\rangle_c]}(s + \tau + \tau_0)\right|
\\
& \le &
\left|(s + \tau'' + \tau_0) - (s + \tau + \tau_0)\right|
\\
& = &
\tau'' - \tau.
\end{eqnarray*}
This finishes the proof of \eqref{lem:exponential_estimate_for_d^times:small_t+tau:inequality(3)}.

If we combine \eqref{lem:exponential_estimate_for_d^times:small_t+tau:inequality(1)},
\eqref{lem:exponential_estimate_for_d^times:small_t+tau:inequality(2)} and
\eqref{lem:exponential_estimate_for_d^times:small_t+tau:inequality(3)}, we get
\begin{eqnarray}
\label{lem:exponential_estimate_for_d^times:small_t+tau:inequality(4)}
d^{\times}_{\FS,x_0}(\phi_{\tau}(a,c,t),\phi_{\tau + \tau_0}(b,c,s))
& \le &
\tau'' - \tau + N \cdot \lambda^{t + \tau'' - \langle a,c|b,x_0\rangle}.
\end{eqnarray}
We estimate
$$- \langle c|x_0\rangle_{a} -t  ~ \le ~ \langle c|b\rangle_{a} - \langle c|x_0\rangle_{a} - t~ = ~ \langle a,c|b,x_0\rangle - t,$$
and for $b \in X$
\begin{eqnarray*}
\lefteqn{-\langle c|x_0\rangle_b - s - \tau_0}
\\
& = &
 -\langle c|x_0\rangle_b - t + \langle a,b|c,x_0\rangle
\\
& = &
\langle a,c|b,x_0\rangle - t - \langle a,c|b,x_0\rangle + \langle a,b|c,x_0\rangle -\langle c|x_0\rangle_b
\\
& = &
\langle a,c|b,x_0\rangle - t + \langle c,a|b,x_0\rangle + \langle a,b|c,x_0\rangle  + \langle b,c|a,x_0\rangle
-  \langle b,c|a,x_0\rangle-\langle c|x_0\rangle_b
\\
& = &
\langle a,c|b,x_0\rangle - t
-  \langle b,c|a,x_0\rangle - \langle c|x_0\rangle_b
\\
& = &
\langle a,c|b,x_0\rangle - t  - \langle c|a \rangle_b
\\
& \le  &
\langle a,c|b,x_0\rangle - t.
\end{eqnarray*}
This inequality holds for $b \in \partial X$ for trivial reasons.
The last two inequalities imply
$$\tau'' ~ \le ~ \langle a,c|b,x_0\rangle - t$$
and hence
$$0 ~ \le ~ \tau'' - \tau ~ \le ~ -\left(t + \tau - \langle a,c|b,x_0\rangle\right).$$
Lemma~\ref{lem:Choice_of_N} applied to $v = \tau'' - \tau$ and $w = -\left(t + \tau - \langle a,c|b,x_0\rangle\right)$ yields
\begin{eqnarray}
\label{lem:exponential_estimate_for_d^times:small_t+tau:inequality(5)}
\tau'' - \tau + N \cdot \lambda^{t + \tau'' - \langle
  a,c|b,x_0\rangle} & \le &
N \cdot \lambda^{t + \tau - \langle a,c|b,x_0\rangle}.
\end{eqnarray}
If we  combine \eqref{lem:exponential_estimate_for_d^times:small_t+tau:inequality(4)},
\eqref{lem:exponential_estimate_for_d^times:small_t+tau:inequality(5)},
we get the desired inequality
\begin{eqnarray*}
d^{\times}_{\FS,x_0}(\phi_{\tau}(a,c,t),\phi_{\tau + \tau_0}(b,c,s))
& \le &
 N \cdot \lambda^{t + \tau - \langle a,c|b,x_0\rangle}.
\end{eqnarray*}
\end{proof}

Now Theorem~\ref{the:exponential_estimate_for_d^times} follows from
Lemma~\ref{lem:exponential_estimate_for_d^times:a=b_in_partial_X},
Lemma~\ref{lem:exponential_estimate_for_d^times:t+tau_and_s+tau_0+tau_large}
and
Lemma~\ref{lem:exponential_estimate_for_d^times:small_t+tau}.


\subsection{Flow estimates for the metric}
\label{subsec:flow_estimates_for_the_metric}

Next we prove Theorem\ref{the:exponential_estimate_for_d_{overline{FS}_x_0}}.
\begin{proof}
We estimate using $e \cdot \lambda > 1$, $0 <  e^{-1} \cdot \lambda  < 1$
and Theorem~\ref{the:exponential_estimate_for_d^times}
\begin{eqnarray*}
\\
\lefteqn{d_{\FS,x_0}(\phi_{\tau}(a,c,t),\phi_{\tau + \tau_0}(b,c,s))}
\\
& = & \int_{-\infty}^{\infty} \frac{d_{\FS,x_0}^{\times}\left(\phi_{\sigma}(\phi_{\tau}(a,c,t)),\phi_{\sigma}(\phi_{\tau + \tau_0}(b,c,s))\right)}
{2 \cdot e^{|\sigma|}} ~d\sigma
\\
& = & \int_{-\infty}^{\infty} \frac{d_{\FS}^{\times}\left(\phi_{\sigma + \tau}((a,c,t)), \phi_{\sigma+\tau+\tau_0}((b,c,s))\right)}
{2 \cdot e^{|\sigma|}} ~d\sigma
\\
& \le & \int_{-\infty}^{\infty}  \frac{N \cdot \lambda^{t + \sigma + \tau  - \langle a,c|b,x_0\rangle}}
{2 \cdot e^{|\sigma|}} ~d\sigma
\\
& =  & \int_{-\infty}^{0} \frac{N \cdot \lambda^{t + \sigma + \tau  - \langle a,c|b,x_0\rangle}}
{2 \cdot e^{|\sigma|}} ~d\sigma
~  + ~
\int_{0}^{\infty} \frac{N \cdot \lambda^{t + \sigma + \tau  - \langle a,c|b,x_0\rangle}}
{2 \cdot e^{|\sigma|}} ~d\sigma.
\\
& =  & \int_{-\infty}^{0} \frac{N \cdot  \lambda^{t + \sigma + \tau  - \langle a,c|b,x_0\rangle}}
{2 \cdot e^{-\sigma}} ~d\sigma
~  + ~
\int_{0}^{\infty} \frac{N \cdot \lambda^{t + \sigma + \tau  - \langle a,c|b,x_0\rangle}}
{2 \cdot e^{\sigma}} ~d\sigma.
\\
& =  & \frac{N}{2} \cdot \lambda^{t+\tau  - \langle a,c|b,x_0\rangle} \cdot \left(
\int_{-\infty}^{0} (e \cdot\lambda)^{\sigma}  ~d\sigma
~  + ~
\int_{0}^{\infty} (e^{-1} \cdot\lambda)^{\sigma}  ~d\sigma\right)
\\
& =  & \frac{N}{2} \cdot \lambda^{t+\tau  - \langle a,c|b,x_0\rangle} \cdot \left(
\left[ \frac{(e \cdot\lambda)^{\sigma}}{\ln(e \cdot \lambda)}\right]_{-\infty}^{0}
~  + ~
\left[\frac{(e^{-1} \cdot\lambda)^{\sigma}}{\ln(e^{-1} \cdot \lambda)}\right]_{0}^{\infty} \right)
\\
& =  & \frac{N}{2} \cdot \lambda^{t+\tau  - \langle a,c|b,x_0\rangle} \cdot \left(
\frac{1}{\ln(e \cdot \lambda)}
~  + ~
\frac{1}{-\ln(e^{-1} \cdot \lambda)} \right)
\\
& =  & \frac{N}{2} \cdot \lambda^{t+\tau  - \langle a,c|b,x_0\rangle} \cdot \left(
\frac{1}{1 + \ln(\lambda)}
~  + ~
\frac{1}{1 - \ln(\lambda)}\right)
\\
& =  &
\frac{N}{1 - \ln(\lambda)^2} \cdot \lambda^{\left(t - \langle a,c|b,x_0\rangle\right)} \cdot \lambda^{\tau}.
\end{eqnarray*}
\end{proof}


\section{The flow estimates for the map $\iota$}
\label{sec:The_flow_estimates_for_the_map_iota}

Let $X$ be a hyperbolic complex and $x_0 \in X$
be a base point. We define a map
 \begin{eqnarray}
 \iota_{x_0} \colon X \times \overline{X} & \to & \FS(X)
 \label{iota_{x_0}:X_times_overline{X}_to_FS(X)}
 \end{eqnarray}
 by
 $\iota_{x_0}(a,c) := $
$$   
 \left\{\begin{array}{lll}
  \left(a,c,\Theta_{[-\langle c|x_0\rangle_a,\langle
      a|x_0\rangle_c]}^{-1}\left(\min\left\{2,\frac{\widehat{d}(a,c)}{2}\right\} - \langle c|x_0\rangle_a\right)\right) & & c \in X \text{ and } a \not= c;
 \\
 \left(a,c,\Theta_{[-\langle c|x_0\rangle_a, \infty]}^{-1}(2 - \langle
   c|x_0\rangle_a)\right) & & c \in \partial X;
\\
c = (c,c,0) & & \text{ if } a = c.
 \end{array} \right.
 $$

We remind the reader that for $a$, $c \in X$ we have
$a = (a,a,0) = (a,a,-\infty) = (a,c,-\infty)$ in $\FS(X)$.

\begin{remark}
\label{rem:what-is-iota}
 Because of Lemma~\ref{lem:d^times_{x_0}((a,b,t),(a,b,s)right)} the
 point $\iota_{x_0}(a,c)$ is $c$ if $ a = c$, is  the unique point
 on the line $(a,c)_{\FS(X)}$ whose distance with respect to
 $d_{\FS,x_0}$ from $a$ is
 $\min\left\{2,\frac{\widehat{d}(a,c)}{2}\right\}$ if $c \in X$ and $a
 \not = c$, and is the unique point
 on the line $(a,c)_{\FS(X)}$ whose distance with respect to
 $d_{\FS,x_0}$ from $a$ is  $2$ if $c \in \partial X$.
\end{remark}

\begin{lemma}\label{lem:t-langle a,c|b,x_0rangle_le_5/2}
Consider $a,b\in X$ and $c \in \overline{X}$ with $a \not= c$ and $t \in \IR$.
Suppose
$$d_{\FS,x_0}((a,c,t),(a,a,0)) ~ = ~ \left\{
\begin{array}{lll}
\min\{2,\frac{\widehat{d}(a,c)}{2}\} & & c \in X;
\\
2                          & & c \in \partial X.
\end{array}\right.$$
Then
$$-\widehat{d}(a,b) ~ \le ~ t - \langle a,c|b,x_0\rangle ~ \le ~  5/2.$$
\end{lemma}
\begin{proof}
Note that $(a,a,0) = (a,c,-\infty) \in \FS(X)$.
We conclude from Lemma~\ref{lem:d^times_{x_0}((a,b,t),(a,b,s)right)},
Lemma~\ref{lem:relating_t+tau_and_s+tau+tau_0}~\ref{lem:relating_t+tau_and_s+tau+tau_0:(1)}
and Lemma~\ref{lem:properties_of_theta_[alpha,beta]}~\ref{lem:properties_of_theta_[alpha,beta]:shift_invariance}
\begin{eqnarray}
\lefteqn{\Theta_{[-\langle c|b\rangle_a,\langle a|b\rangle_c]}(t  - \langle a,c|b,x_0\rangle)}
\label{lem:t-langle a,c|b,x_0rangle_le_5/2:(1)}
\\
& = &
\Theta_{[-\langle c|x_0\rangle_a  - \langle a,c|b,x_0\rangle,\langle a|x_0\rangle_c  - \langle a,c|b,x_0\rangle]}(t  - \langle a,c|b,x_0\rangle)
\nonumber
\\
& = &
\Theta_{[-\langle c|x_0\rangle_a,\langle a|x_0\rangle_c]}(t) - \langle a,c|b,x_0\rangle
\nonumber
\\
& = &
\Theta_{[-\langle c|x_0\rangle_a,\langle a|x_0\rangle_c]}(t) - (-\langle c|x_0\rangle_a) + (-\langle c|x_0\rangle_a) - \langle a,c|b,x_0\rangle
\nonumber
\\
& = &
\Theta_{[-\langle c|x_0\rangle_a,\langle a|x_0\rangle_c]}(t) - (-\langle c|x_0\rangle_a) - \langle c|b\rangle_a
\nonumber
\\
& = &
d_{\FS,x_0}((a,c,t),(a,a,0)) - \langle c|b\rangle_a.
\nonumber\end{eqnarray}

If $c \in X$, we get by assumption
$$d_{\FS,x_0}((a,c,t),(a,a,0)) ~ \le ~ \frac{\widehat{d}(a,c)}{2} ~ = ~\frac{\langle a|b\rangle_c - (-\langle c|b\rangle_a)}{2},$$
and hence
$$d_{\FS,x_0}((a,c,t),(a,a,0)) - \langle c|b\rangle_a ~ \le ~
\frac{\langle a|b\rangle_c + (-\langle c|b\rangle_a)}{2}.$$
This inequality is true for $c \in \partial X$ for trivial reasons.
Hence we get from~\eqref{lem:t-langle a,c|b,x_0rangle_le_5/2:(1)}
$$\Theta_{[-\langle c|b\rangle_a,\langle a|b\rangle_c]}(t  - \langle a,c|b,x_0\rangle) ~ \le ~
\frac{-\langle c|b\rangle_a + \langle a|b\rangle_c}{2}.$$
Lemma~\ref{lem:properties_of_theta_[alpha,beta]}~\ref{lem:properties_of_theta_[alpha,beta]:bound_on_t}
together with \eqref{lem:t-langle a,c|b,x_0rangle_le_5/2:(1)} implies
\begin{eqnarray*}
t  - \langle a,c|b,x_0\rangle & \le &
\Theta_{[-\langle c|b\rangle_a,\langle a|b\rangle_c]}(t  - \langle a,c|b,x_0\rangle) + 1/2
\\
& = &
d_{\FS,x_0}((a,c,t),(a,a,0)) - \langle c|b\rangle_a + 1/2
\\
& \le  &
d_{\FS,x_0}((a,c,t),(a,a,0)) + 1/2
\\
& \le & 2 + 1/2
\\
& = & 5/2.
\end{eqnarray*}

Thus we have proven the upper bound $t  - \langle a,c|b,x_0\rangle \le 5/2$.
It remains to show the lower bound
$ - \widehat{d}(a,b) ~ \le ~ t - \langle a,c|b,x_0\rangle$.

We conclude from the assumptions that
$$d_{\FS,x_0}((a,c,t),(a,a,0)) ~ = ~ \frac{\widehat{d}(a,c)}{2}  ~ = ~
\frac{\langle a|b\rangle_c - (-\langle c|b\rangle_a)}{2} ~ \text { and } c \in X$$
or
$$d_{\FS,x_0}((a,c,t),(a,a,0)) ~ = ~ 2$$
holds. We begin with the first case. Then
$$d_{\FS,x_0}((a,c,t),(a,a,0)) - \langle c|b\rangle_a ~ \ge ~
\frac{\langle a|b\rangle_c + (-\langle c|b\rangle_a)}{2}.$$
Lemma~\ref{lem:properties_of_theta_[alpha,beta]}~\ref{lem:properties_of_theta_[alpha,beta]:bound_on_t}
together with \eqref{lem:t-langle a,c|b,x_0rangle_le_5/2:(1)} implies
$$t  - \langle a,c|b,x_0\rangle ~ \ge ~ -\langle c|b\rangle_a ~ \ge ~ - \widehat{d}(a,b).$$
Finally we treat the second case. Then \eqref{lem:t-langle a,c|b,x_0rangle_le_5/2:(1)} implies
$$\Theta_{[-\langle c|b\rangle_a,\langle a|b\rangle_c]}(t  - \langle a,c|b,x_0\rangle)
~ \ge ~ 2 - \langle c|b\rangle_a.$$
Since for $u \le -\langle c|b\rangle_a$ we have
$$\Theta_{[-\langle c | b\rangle_a,\langle a|b\rangle_c]}(u) ~ = ~
 -\langle c | b\rangle_a + e^{u -(-\langle c|b\rangle_a)}/2 - e^{u -\langle a|b\rangle_c}/2
~ < ~ 2 - \langle c|b\rangle_a,$$
we conclude
$$t  - \langle a,c|b,x_0\rangle ~ \ge ~  -\langle c|b\rangle_a ~ \ge ~ - \widehat{d}(a,b).$$
This finishes the proof of
Lemma~\ref{lem:t-langle a,c|b,x_0rangle_le_5/2}.
\end{proof}

\begin{lemma} \label{lem:iota_is_continuous}
The map $ \iota_{x_0} \colon X \times \overline{X} \to  \FS(X)$
from~\eqref{iota_{x_0}:X_times_overline{X}_to_FS(X)} is continuous. It is $\Isom(X)$-equivariant with respect
to the diagonal $\Isom(X)$-action on the source. For $x \in X$ the map
$\iota_{x_0}(x,-) \colon \overline{X} \to \FS(X), \quad y \mapsto \iota_{x_0}(x,y)$ is injective.
\end{lemma}

\begin{proof}
We only prove continuity, the other claims are straight-forward to check
using Remark~\ref{rem:what-is-iota}.
Recall that $X$ and $\FS(X)$ are metric spaces. 
By \cite[Exercise~III.H.3.18(4)]{Bridson-Haefliger-Buch} 
the space $\overline{X}$ is metrizable.
Hence it suffices to check continuity
for sequences.
Consider sequences $(a_n)_{n \ge 0}$ in $X$ and
$(c_n)_{n \ge 0}$ in $\overline{X}$ and points $a \in X$ and $c \in \overline{X}$ such that
$\lim_{n \to \infty} a_n = a$ in $X$ and $\lim_{n \to \infty} c_n = c$ in $\overline{X}$ hold. We have to show
$\lim_{n \to \infty} \iota_{x_0}(a_n,c_n) = \iota_{x_0}(a,c)$ in $\FS(X)$.

Suppose that $a = c$. Then we can assume $c_n \in X$ and $\widehat{d}(a_n,c_n)/2 < 2$ for $n \ge 0$
and $\lim_{n \to \infty} \widehat{d}(a_n,c_n) = 0$. 
This implies by the construction of $\iota_{x_0}$
that  $\lim_{n \to \infty}d_{\FS,x_0}(a_n,\iota_{x_0}(a_n,c_n)) 
     = \lim_{n \to \infty}d_{\FS,x_0}(a_n,c_n)/2 ~ = ~ 0$
and hence
$$\lim_{n \to \infty} \iota_{x_0}(a_n,c_n) = \lim_{n \to \infty} a_n = a = c =    \iota_{x_0}(a,c).$$
Hence we can assume without loss of generality that $a \not= c$ and $a_n \not= c_n$ for all $n \ge 0$ holds.

For $n \ge 0$ put $\alpha_n = -\langle c_n|x_0\rangle_{a_n}$  and
$\beta_n = \langle a_n|x_0\rangle_{c_n}$. Put $\alpha = -\langle c|x_0\rangle_{a}$  and
$\beta = \langle a|x_0\rangle_{c}$.
Then the continuity of the Gromov product
(see Theorem~\ref{the:extending_the_double_difference}) implies
\begin{eqnarray*}
\lim_{n \to \infty} \alpha_n  & = & \alpha;
\\
\lim_{n \to \infty} \beta_n   & = & \beta.
\end{eqnarray*}
Define $t_n$ to be the real number satisfying
\begin{equation}
\label{eq:define-t_n}
\begin{array}{lclll}
\Theta_{[\alpha_n,\beta_n]}(t_n)
      & = &
      \min\left\{2,\frac{\widehat{d}(a_n,c_n)}{2}\right\} - \langle c_n|x_0\rangle_{a_n}
      & & c_n \in X;
\\
\Theta_{[\alpha_n,\beta_n]}(t_n)
      & = &
      2  - \langle c_n|x_0\rangle_{a_n}
      & & c_n \in \partial X.
\end{array}
\end{equation}
Define $t$ to be the real number satisfying
\begin{equation}
\label{eq:define-t}
\begin{array}{lclll}
\Theta_{[\alpha,\beta]}(t)
      & = &
      \min\left\{2,\frac{\widehat{d}(a,c)}{2}\right\} - \langle c|x_0\rangle_a
      & & c \in X;
\\
\Theta_{[\alpha,\beta]}(t)
      & = &
      2  - \langle c|x_0\rangle_a
      & & c \in \partial X.
\end{array}
\end{equation}

Then $\iota_{x_0}(a_n,c_n) = (a_n,c_n,t_n)$ and
$\iota_{x_0}(a,c) = (a,c,t)$.
Because of Theorem~\ref{the:topology of FS(X) - X}
it suffices to show that $\lim_{n \to \infty} t_n = t$
holds.
{}From Lemma~\ref{lem:t-langle a,c|b,x_0rangle_le_5/2}
applied in the case $b = x_0$ we conclude
$-\widehat{d}(a_n,x_0) \le t_n \le 5/2$ and
$-\widehat{d}(a,x_0) \le t \le 5/2$.
Hence we can assume without loss of generality for
all $n \ge 0$
$$-\widehat{d}(a,x_0) - 1 \le t_n,t \le 5/2.$$
We conclude from the Mean Value Theorem for some number $\rho \in [0,1]$
$$\Theta_{[\alpha,\beta]}(t_n) - \Theta_{[\alpha,\beta]}(t) ~ = ~
\Theta_{[\alpha,\beta]}'(\rho \cdot t_n + (1-\rho) t) \cdot (t_n -t).$$
Choose a constant $C > 0$ such that $\Theta_{[\alpha,\beta]}'(s) \ge C^{-1} $ holds for
$s \in [-\widehat{d}(a,x_0) - 1,5/2]$.
Since $\rho \cdot t_n + (1-\rho) t$ lies in
$[-\widehat{d}(a,x_0) - 1,5/2]$, we conclude for all $n \ge 0$
$$|t_n - t| \le C \cdot |\Theta_{[\alpha,\beta]}(t_n) - \Theta_{[\alpha,\beta]}(t)|.$$
We get from the triangle inequality
\begin{eqnarray*}
 |\Theta_{[\alpha,\beta]}(t_n) - \Theta_{[\alpha,\beta]}(t)|
& \le &
 |\Theta_{[\alpha_n,\beta_n]}(t_n) -  \Theta_{[\alpha,\beta]}(t_n)|
+ |\Theta_{[\alpha_n,\beta_n]}(t_n) - \Theta_{[\alpha,\beta]}(t)|
\end{eqnarray*}
and hence
\begin{eqnarray}
\hspace{8mm} |t_n -t|
& \le & C \cdot \left(
 |\Theta_{[\alpha_n,\beta_n]}(t_n) -  \Theta_{[\alpha,\beta]}(t_n)|
+  |\Theta_{[\alpha_n,\beta_n]}(t_n) - \Theta_{[\alpha,\beta]}(t)|\right).
\label{lem:iota_is_continuous:|t_n-t|}
\end{eqnarray}
Since $a_n,a \in X$, we have $- \infty < \alpha_n$ and $- \infty < \alpha$.
We have $\beta < \infty$ if and only if $c \in X$ and
$\beta_n < \infty$ if and only if $c_n \in X$. If $\beta = \infty$,
we can assume without loss of generality $5/2 \le \beta_n$ and hence
$t_n \le \beta_n$ for all $n \ge 0$. We conclude
$$\begin{array}{lclll}
|\Theta_{[\alpha_n,\beta_n]}(t_n) -  \Theta_{[\alpha,\beta]}(t_n)| & \le &
\max\{|\alpha_n - \alpha|,|\beta_n - \beta|\} & & \beta_n, \beta < \infty;
\\
|\Theta_{[\alpha_n,\beta_n]}(t_n) -  \Theta_{[\alpha,\beta]}(t_n)| & \le &
|\alpha_n - \alpha| & & \beta_n = \beta =\infty;
\\
|\Theta_{[\alpha_n,\beta_n]}(t_n) -  \Theta_{[\alpha,\beta]}(t_n)| & \le &
|\alpha_n - \alpha|  + e^{t_n -\beta_n}& & \beta_n < \infty,  \beta = \infty.
\end{array}
$$
from
Lemma~\ref{lem:properties_of_theta_[alpha,beta]}~\ref{lem:properties_of_theta_[alpha,beta]:dependency_on_alpha_and_beta_part1}
and~\ref{lem:properties_of_theta_[alpha,beta]:dependency_on_alpha_and_beta_part2} and (the triangle inequality
in the last case).
Since $t_n \le 5/2$, $\lim_{n \to \infty} \alpha_n = \alpha$ and
$\lim_{n \to \infty} \beta_n = \beta$, we conclude
$$\lim_{n \to \infty} |\Theta_{[\alpha_n,\beta_n]}(t_n) -  \Theta_{[\alpha,\beta]}(t_n)| ~ = ~ 0.$$
Since $\widehat{d}$ and the Gromov product are continuous, we get
using \eqref{eq:define-t_n} and \eqref{eq:define-t}
$$\lim_{n \to \infty}  |\Theta_{[\alpha_n,\beta_n]}(t_n) - \Theta_{[\alpha,\beta]}(t)| ~ = ~ 0.$$
Now~\eqref{lem:iota_is_continuous:|t_n-t|} implies $\lim_{n \to \infty} t_n = t$.
This finishes the proof of Lemma~\ref{lem:iota_is_continuous}.\end{proof}

 \begin{theorem}[Flow estimate for $\iota$]
 \label{the:flow_estimate_for_d_{overline{FS}_x_0}_for_iota}
 Let $\lambda \in (e^{-1},1)$ and $T \in [0,\infty)$
 be the constants depending only on $X$ which appear in
 Proposition~\ref{prop:Mineyevs_Proposition_38}.
 Consider $a,b \in X$ and $c \in \overline{X}$.
 Put
 \begin{eqnarray*}
 N & = & 2 + \frac{2}{\lambda^T \cdot (-\ln(\lambda))}.
 \end{eqnarray*}

 Then there exists a real number $\tau_0$ such that
 $$|\tau_0| ~ \le 2 \cdot \widehat{d}(a,b) + 5$$
 and for all $\tau \in \IR$
 $$d_{\FS,x_0}(\phi_{\tau}\iota_{x_0}(a,c),\phi_{\tau + \tau_0}\iota_{x_0}(b,c)) ~ \le ~
 \frac{N}{1 - \ln(\lambda)^2} \cdot \lambda^{-\widehat{d}(a,b)} \cdot \lambda^{\tau}.$$
 \end{theorem}

 \begin{proof}
Put
 \begin{eqnarray*}
 t & := &
 \left\{
 \begin{array}{lll}
 \left(\Theta_{[-\langle c|x_0\rangle_a,\langle a|x_0\rangle_c]}\right)^{-1}
\left(\min\left\{2,\frac{\widehat{d}(a,c)}{2}\right\} - \langle c|x_0\rangle_a\right)
 & & \text{ if } c \in X  \text{ and } a \not= c;
 \\
 \left(\Theta_{[-\langle c|x_0\rangle_a,\infty]}\right)^{-1}\left(2 - \langle c|x_0\rangle_a\right)
 & & \text{ if } c \in \partial X;
\\
\langle a,c|b,x_0\rangle = 0
& & \text{ if }  a = c;
 \end{array}\right.
 \\
 s & := &
 \left\{
 \begin{array}{lll}
 \left(\Theta_{[-\langle c | x_0\rangle_b,\langle b|x_0\rangle_c]}\right)^{-1}
\left(\min\left\{2,\frac{\widehat{d}(b,c)}{2}\right\} - \langle c|x_0\rangle_b\right)
 & & \text{ if } c \in X   \text{ and } b \not= c;
 \\
 \left(\Theta_{[-\langle c|x_0\rangle_b,\infty]}\right)^{-1}\left(2 - \langle c|x_0\rangle_b\right)
 & & \text{ if } c \in \partial X;
\\
\langle b,c|a,x_0 \rangle = 0
&  &
\text{ if } b = c.
 \end{array}\right.
 \end{eqnarray*}
Put
$$ \tau_0 ~ = ~ t - s  - \langle a,b|c,x_0\rangle.$$

We have by definition
 \begin{eqnarray*}
 \iota_{x_0}(a,c) & = & (a,c,t);
 \\
 \iota_{x_0}(b,c) & = & (b,c,s);
 \\
 \tau_0 & = &  t - s -  \langle a,b|c,x_0 \rangle ~ = ~ (t - \langle a,c|b,x_0\rangle) - (s - \langle b,c|a,x_0\rangle).
 \end{eqnarray*}

 We conclude from Theorem~\ref{the:exponential_estimate_for_d_{overline{FS}_x_0}}
 $$d_{\FS,x_0}(\phi_{\tau}\iota_{x_0}(a,c),\phi_{\tau + \tau_0}\iota_{x_0}(b,c)) ~ \le ~
 \frac{N}{1 - \ln(\lambda)^2} \cdot \lambda^{\left(t - \langle a,c|b,x_0\rangle\right)} \cdot \lambda^{\tau}.$$
 We have
 $$d_{\FS,x_0}((a,c,t),(a,a,0)) ~ = ~ \left\{
 \begin{array}{lll}
 \min\{2,\frac{\widehat{d}(a,c)}{2}\} & & c \in X;
 \\
 2                          & & c \in \partial X.
 \end{array}\right.$$
 and
 $$d_{\FS,x_0}((b,c,s),(b,b,0)) ~ = ~ \left\{
 \begin{array}{lll}
 \min\{2,\frac{\widehat{d}(b,c)}{2}\} & & c \in X;
 \\
 2                          & & c \in \partial X.
 \end{array}\right.$$
 We conclude from
 Lemma~\ref{lem:t-langle a,c|b,x_0rangle_le_5/2} and the
 definition of $t$ and $s$ respectively that
 $$\begin{array}{lclcl}
 -\widehat{d}(a,b) & \le & t - \langle a,c|b,x_0 \rangle & \le &  5/2;
 \\
 -\widehat{d}(a,b) & \le & s - \langle b,c|a,x_0 \rangle & \le &  5/2,
 \end{array}
 $$
 holds. This finishes the proof of  Theorem~\ref{the:flow_estimate_for_d_{overline{FS}_x_0}_for_iota}.
 \end{proof}

We can now prove the flow estimate from the introduction.

\begin{proof}[Proof of Theorem~\ref{thm:final-flow-estimate}]
Define $j \colon G \x \overline{X} \to \FS(X)$
by $j(g,c) := \iota_{x_0} (g x_0,c)$.
It follows from Lemma~\ref{lem:iota_is_continuous}
that $j$ is continuous and $G$-equivariant with respect to the
diagonal $G$-action on the source.
Let $C := \max \{ \widehat{d}(x_0, s x_0) \mid s \in S \}$
where $S$ is the generating set $S$ of $G$ used to define
the word metric $d_G$ on $G$.
Then
\[
\widehat{d} (g x_0, h x_0) \leq C d_G(g, h) \quad \forall \;
                    g, h \in G.
\]
Let $\alpha > 0$ be given.
Let $\lambda \in (e^{-1},1)$ and $T \in [0,\infty)$
be the constants depending only on $X$ which appear in
Proposition~\ref{prop:Mineyevs_Proposition_38}.
Let $N$ be the number defined in Theorem~
\ref{the:flow_estimate_for_d_{overline{FS}_x_0}_for_iota}.
Define
\begin{eqnarray*}
\beta(\alpha) & :=  & 2 C \alpha + 5 \\
f_\alpha(\tau) & := &  \frac{N}{1 - \ln(\lambda)^2}
     \cdot \lambda^{-C \alpha} \cdot \lambda^{\tau}.
\end{eqnarray*}
It follows from Theorem~
\ref{the:flow_estimate_for_d_{overline{FS}_x_0}_for_iota}
that $\beta(\alpha)$ and $f_\alpha$ satisfy the
assertion of Theorem~\ref{thm:final-flow-estimate}.
\end{proof}


\section{Further properties of the flow space}
\label{sec:Properness_of_the_induced_action}

\begin{theorem}
\label{the:Properness_of_the_geodesic_flow_space}
Let $X$ be a hyperbolic complex with base point $x_0 \in X$. Suppose that $G$ acts on $X$ by
simplicial automorphisms such that every isotropy group is finite. Then

\begin{enumerate}

\item \label{the:Properness_of_the_geodesic_flow_space:properness_as_metric_space}
The metric space $(\FS(X),d_{\FS,x_0}(X))$ is proper;

\item \label{the:Properness_of_the_geodesic_flow_space:properness_of_action}
The induced $G$-action on the flow space $(\FS(X),d_{\FS,x_0})$ is proper;

\item \label{the:Properness_of_the_geodesic_flow_space:cocompact}
If $G$ acts cocompactly on $X$, then $G$ acts cocompactly on
$\FS(X)$.
\end{enumerate}
\end{theorem}

A $G$-space $Y$ is called \emph{proper} if for every $y \in Y$ there exists
an open neighborhood $U$ such that the set
$\{g \in G \mid g \cdot U \cap U \not= \emptyset\}$
is finite. Proper implies that all isotropy groups are
finite. The converse is not true in general but is true if $Y$ is a
$G$-$CW$-complex (\cite[Theorem~1.23 on page~18]{Lueck(1989)}).
If $(Y,d_Y)$ is a metric space and $G$ acts by isometries, the
$G$-action is proper if and only if for every $y \in Y$ there exists
an $\epsilon$ such that the set
$\{g \in G \mid g \cdot B_{\epsilon}(y) \cap B_{\epsilon}(y)  \not
= \emptyset\}$ is finite, where $B_{\epsilon}(y) = \{z \in Y \mid
d_Y(y,z) < \epsilon\}$. A metric space $(Y,d_Y)$ is called \emph{proper} if
and only if $\overline{B_{\epsilon}(y)} = \{z \in Y \mid d_Y(y,z) \le
\epsilon\}$ is compact for all $y \in Y$ and $\epsilon \ge 0$.

The elementary proof of the next result is left to the reader.

\begin{lemma} \label{lem:criterion_for_properness_for_metric_spaces}
Let $(Y,d_Y)$ be a proper metric space. Let $G$ act on $Y$ by
isometries.
Then the $G$-action on $Y$ is proper if and only if
for every $C > 0$ and $y \in Y$ the set
$\left\{g\mid d_Y(g\cdot y,y) \le C\right\}$ is finite.
\end{lemma}

Now we are ready to prove
Theorem~\ref{the:Properness_of_the_geodesic_flow_space}.
\begin{proof}
\ref{the:Properness_of_the_geodesic_flow_space:properness_as_metric_space}
This is proven in \cite[Proposition~54~on~page~464]{Mineyev-flows-and-joins}.
\\[1mm]
\ref{the:Properness_of_the_geodesic_flow_space:properness_of_action}
Mineyev~\cite[page~457]{Mineyev-flows-and-joins} constructs a map
$$\Psi \colon \FS(X) \to X$$
such that there exists constants $K_1$ and $K_2$ depending only on $X$ such that
for all $v,w \in \FS(X)$ and $g \in G$ we have
\begin{eqnarray*}
\left|d_{\FS,x_0}(v,w) - \widehat{d}(\Psi(v),\Psi(w))\right| & \le & K_1;
\\
\widehat{d}(\Psi (g \cdot v),g \cdot \Psi(v)) & \le & K_2.
\end{eqnarray*}
Compare also \cite[Proposition 43 on page~458]
{Mineyev-flows-and-joins}.

We conclude for $v \in \FS(X)$ and $g \in G$
\begin{eqnarray*}
\lefteqn{\widehat{d}(g \cdot \Psi(v), \Psi(v)) }
\\
& \le &
\widehat{d}(\Psi(g \cdot v), \Psi(v)) +  \widehat{d}(\Psi(g \cdot v),
g \cdot \Psi(v))
\\
& = &
d_{\FS,x_0}(g \cdot v,v) - d_{\FS,x_0}(g \cdot v,v) +
\widehat{d}(\Psi(g \cdot v), \Psi(v)) +  \widehat{d}(\Psi(g \cdot v),
g \cdot \Psi(v))
\\
& \le &
d_{\FS,x_0}(g \cdot v,v) + \left|d_{\FS,x_0}(g \cdot v,v) -
  \widehat{d}(\Psi(g \cdot v), \Psi(v))\right| +
\widehat{d}(\Psi(g \cdot v), g \cdot \Psi(v))
\\
& \le &
d_{\FS,x_0}(g \cdot v,v) + K_1 + K_2.
\end{eqnarray*}
There exist real numbers  $A \ge 1$ and $B \ge 0$ depending only on $X$ such that for all
$x_1,x_2 \in X$ we have
\begin{eqnarray*}
& A^{-1} \widehat{d}(x_1,x_2) - B ~ \le d(x_1,x_2) ~ \le ~ A \cdot \widehat{d}(x_1,x_2) + B &
\end{eqnarray*}
where $d$ is the word metric, compare the beginning of Section~\ref{subsec:Hyperbolic complexes, double difference and Gromov product}.
Hence we get
\begin{eqnarray*}
d(g \cdot \Psi(v), \Psi(v))  & \le  &  A \cdot \widehat{d}(g \cdot \Psi(v), \Psi(v)) + B
\\
& \le &
A \cdot \left(d_{\FS,x_0}(g \cdot v,v) + K_1 + K_2\right) + B.
\end{eqnarray*}
Consider $v \in \FS(X)$ and $C \ge 0$. Since $G$ acts properly on $X$
and $(X,d)$ is a proper metric space,
Lemma~\ref{lem:criterion_for_properness_for_metric_spaces}
implies that
the set
$$\left\{g \in G \mid d(g \cdot \Psi(v),\Psi(v)) \le A \cdot (C + K_1 + K_2) + B\right\}$$
is finite. Since this set contains
$\left\{g \in G \mid d_{\FS,x_0}(g \cdot v,v) \le C\right\}$, also the latter set
is finite. Hence the $G$-action on $(\FS(X),d_{\FS,x_0})$ is proper by
Lemma~\ref{lem:criterion_for_properness_for_metric_spaces}.
\\[1mm]
\ref{the:Properness_of_the_geodesic_flow_space:cocompact}
Since $G$ acts simplicially and cocompactly on $X$, we can find a compact subset $C \subseteq X$ such that
$G \cdot C = X$. Consider $D = \{x \in \FS(X) \mid \widehat{d}(\Psi(x),C) \le K_2\}$.
Since $C$ is compact, its diameter $\diam(C)$ is finite.
Since for $y,z \in D$ we get
$$d_{\FS,x_0}(y,z) \le \widehat{d}(\Psi(y),\Psi(z)) + K_1 \le \diam(C) + 2 K_2 + K_1$$
the set $D$ has finite diameter. Since $\FS(X)$ is proper as a metric space by
assertion~\ref{the:Properness_of_the_geodesic_flow_space:properness_as_metric_space},
the closure of $D$ is a compact subset of $\FS(X)$. Next we show $G\cdot D = \FS(X)$.
Consider $x \in \FS(X)$. Choose $g \in G$ such that $g^{-1}\Psi(x) \in C$.
{}From $\widehat{d}(\Psi(g^{-1}x),g^{-1} \Psi(x)) \le K_2$ 
we conclude $g^{-1}x \in D$.
This implies $x \in g\cdot D \subseteq G\cdot D$.
Hence $\overline{D}$ is a compact subset of $\FS(X)$ with $G\cdot \overline{D} = \FS$.
Therefore $G$ acts on $\FS$ cocompactly.
This finishes the proof of
Theorem~\ref{the:Properness_of_the_geodesic_flow_space}.
\end{proof}

The following facts are well-known.
We include a proof for the convenience of the reader.

\begin{lemma}
\label{lem:on-the-topology-of-overline-X}
Let $X$ be a $\delta$-hyperbolic complex in the sense of
Section~
\ref{subsec:Hyperbolic complexes, double difference and Gromov product}.
Let $\overline{X}$ be the compactification of $X$ in the
sense of Gromov. Then
\begin{enumerate}
\item \label{lem:on-overline-X:locally-connected}
      $\overline{X}$ is locally connected;
\item \label{lem:on-overline-X:finite-dim}
      $\overline{X}$ has finite covering dimension.
\end{enumerate}
\end{lemma}

\begin{proof}
We start by reviewing the topology of $\overline{X}$
following \cite[p.429]{Bridson-Haefliger-Buch}.
Recall that  $X^{(1)}$ denotes the $1$-skeleton of $X$.
A generalized ray $c \colon I \to X^{(1)}$ is a geodesic
with respect to the word metric $d_{\mathit{word}}$,
where $I = [0,R]$ for $R \geq 0$ or $I = [0,\infty)$.
In the later case $c$ will be called a geodesic ray.
If $I = [0,R]$ it is convenient to write
$c(t) = c(R)$ for $t \geq r$.
Two geodesic rays $c$, $c'$ are called equivalent if
there is $C > 0$ such that
${d}_{\mathit{word}}(c(t),c'(t)) < C$ for all
$t \in [0,\infty)$.
$\dd X = \dd X^{(1)}$ is the set of all equivalence classes
of such geodesic rays.
For a geodesic ray $c$ we denote by $c(\infty)$ the point in
$\dd X$ determined by $c$.
Fix a base point $x_0 \in X^{(0)}$ and $k > 2 \delta$.
Every point in $\dd X$ can be written as $c(\infty)$
where $c$ is a geodesic ray starting at $x_0$,
\cite[Lemma~3.1 on p.427]{Bridson-Haefliger-Buch}.
For a geodesic ray $c$ starting at $x_0$ and $n \in \IN$
let $V^{(1)}_n(c)$ denote the set of all  $c'(\infty)$
where $c'$ is a generalized ray
starting at $x_0$ with $d_{\mathit{word}}(c(n),c'(n)) < k$.
(Such a generalized ray may end in $X$.)
Let $V_n(c)$ be the union of $V^{(1)}_n(c)$ with the smallest
subcomplex of $X$ containing $V^{(1)}_n(c) \cap X$.
The topology on $\overline{X} = X \cup \dd X$
is now defined as follows: a $U \subset \overline{X}$
is open if and only if the following two condition hold
\begin{itemize}
\item $U \cap X$ is open in $X$.
\item If $c$ is a geodesic ray starting at $x_0$ and
      $c(\infty) \in U$ then there is $n \in \IN$
      such that $V_n(c) \subset U$.
\end{itemize}

We now prove
\ref{lem:on-overline-X:locally-connected}.
Obviously $X$ is locally connected as it is  a $CW$-complex.
It suffices to show for every
$x \in \partial X = \overline{X} - X$
that there is a (not necessarily open) connected neighborhood
(see~\cite[Exercise~10~on~page~163]{Munkres(1975)}).
But this follows because the $V_n(c)$ from above are connected.

Next we prove
\ref{lem:on-overline-X:finite-dim}.
Let $D \in \IN$ be the maximal number of points in
$X^{(0)}$ that are contained in a ball of radius
$k+ 2 \delta$.
This is a finite number because $X$ is uniformly locally finite.
Because $X$ has finite covering dimension and $\dd X$
is compact, it suffices to show the following:
For every finite collection $\calu$ of open subsets of
$\overline{X}$ that covers $\dd X$ there is an
$(D-1)$-dimensional refinement of $\calu$ that still covers
$\dd X$.

We will need the following observation:
Let $c$, $c'$ be two geodesics starting at $x_0$.
Let $x$ be the endpoint of $c$ and
$x'$ be the endpoint of $c'$.
If $N$ is such that
$d(c(N),x) > d(x,x') + \delta$,
then $d(c(N),c'(N)) \leq 2\delta$.
This is an easy application of the condition that all geodesic
triangles in $X^{(1)}$ are $\delta$-thin.

For $N \in \IN$ consider the set $S_N$ of all $c(N)$ where
$c$ is a geodesic ray in $X^{(1)}$ starting at $x_0$.
This is a finite subset of $X^{(0)}$ because $X$
is locally finite.
For every $x \in S_N$ pick a geodesic ray $c_x$ starting at
$x_0$ such that $c_x(N) = x$.
Denote by $B$ the close ball of radius $n+\delta$ around $x_0$.
Let $U_x := (V_N(c_x))^\circ - B$.
Then the collection
$\calu_N := \{ U_x    \mid x \in S_N \}$
covers $\dd X$.
Let $z \in \overline{X} - B$ and choose
$x \in S_N$ such that $x$ lies on a geodesic from $x_0$ to $z$.
Using the above observation it is not hard to show that
for every $x' \in S_N$ with $z \in U_{x'}$ we have
$d(x,x') < k + 2 \delta$.
Therefore the dimension of $\calu_N$ is bounded by $D-1$.

It remains to show that for sufficiently large $N$,
the collection $\calu_N$ will be a refinement of
the given collection $\calu$.
Let $R > k + 3\delta$, $R \in \IN$.
By the definition of the topology of $\dd X$
and because $\dd X$ is compact, there are
$M_1,\dots,M_n \in \IN$ and
geodesic rays $c_1, \dots, c_n$ starting at $x_0$
such that
\begin{itemize}
\item $\dd X \subset V_{M_1 + R}(c_1) \cup
                    \dots \cup V_{M_n + R}(c_n)$;
\item for every $i=1,\dots,n$ there is $U \in \calu$
      such that $V_{M_i}(c_i) \subset U$.
\end{itemize}
Let $M := \max\{M_1,\dots,M_n\}$.
Two applications of the above observation give
the following:
If $c$ is a geodesic ray starting at $x_0$ such that
$c(\infty) \in V_{M_i + R}(c_i)$,
then $V_{M+2R}(c) \subset V_{M_i}(c_i)$.
Thus $\calu_{M+2R}$ is a refinement of $\calu$.
\end{proof}

We can now check the additional properties stated in
Section~\ref{subsec:intro-flow}.

\begin{proof}
[Proof of Proposition~\ref{prop:further-properties-of-FS}]
\ref{prop:further-properties-of-FS:finite-subgroups}
This follows from
\cite[Theorem~3.2 on p.459]{Bridson-Haefliger-Buch}.
\\[1mm]
\ref{prop:further-properties-of-FS:connected+dimension}
By Lemma~\ref{lem:on-the-topology-of-overline-X}
$\overline{X}$ has finite
covering dimension and is locally connected.
It follows therefore from Theorem~\ref{the:topology of FS(X) - X}
that $\FS(X) - \FS(X)^\IR = \FS(X) - X$ is locally connected
and has finite covering dimension.
\\[1mm]
\ref{prop:further-properties-of-FS:proper+cocompact}
This follows from
Theorem~\ref{the:Properness_of_the_geodesic_flow_space}.
\\[1mm]
\ref{prop:further-properties-of-FS:closed-orbits}
For $g \in G$ its translation length on
$(\FS(X),d_{\FS(X),x_0})$ is defined as
\[
l(g) := \lim_{n \to \infty}
             d_{\FS,x_0}(g^n (a,b,t), (a,b,t) ) / n.
\]
By the triangle inequality this definition does not
depend on the choice of $(a,b,t) \in \FS(X)$.
In particular, it depends only on the conjugacy class
of $g$.
Since the isometric $G$-action on $\FS(X)$ is cocompact and proper,
$(G,d_G)$ is quasi-isometric to $(\FS,d_{\FS,x_0})$.
Thus, there are constants $A \geq 1$, $B > 0$
such that
\[
A \cdot l(g) + B \geq \tau(g) := \lim_{n \to \infty}
                 d_G(g^n,1_G).
\]
($\tau(g)$ is the translation length of $g$ on $(G,d_G)$.)
By \cite[Proposition~3.15 on p.465]{Bridson-Haefliger-Buch}
for fixed $C > 0$ the number of conjugacy classes whose
translation length on $(G,d_G)$ is no more than $C$
is finite.
We conclude that the same holds for the translation length
on $(\FS(X),d_{\FS,x_0})$.

Fix $C>0$.
Let $\call$ be the set of all orbits $L$ of
the flow $\phi_\tau$ on $\FS(X)$ with $0 < \per_\phi^G(L) \leq C$,
see Definition~\ref{def:G-period}.
Every $L \in \call$ is a line $(a_L,b_L)_{\FS(X)}$
with $a_L$, $b_L \in \overline{X}$.
For $L \in \call$ there is $g_L \in G$ such that
\[
g_L \cdot (a_L,b_L,t) = (a_L,b_L, t + \per_\phi^G(L))
\]
for $(a_L,b_L,t) \in L$.
In particular $g_L \cdot a_L = a_L$, $g_L \cdot b_L = b_L$.
If $a_L \in X$ or $b_L \in X$, then $g_L$ has finite order because
the action of $G$ on $X$ is proper. 
But this would imply $\per_\phi^G(L) = 0$.
Therefore $a_L$, $b_L \in \dd X$.
By Lemma~\ref{lem:d^times_{x_0}((a,b,t),(a,b,s)right)}
$l(g_L) = \per_\phi^G(L) \leq C$.
Recall that $\dd X \cong \dd G$, because $X$ is quasi-isometric
to $G$, \cite[Theorem~3.9 on p.430]{Bridson-Haefliger-Buch}.
Because every $g \in G$ has at most $2$ fixed points
on $\dd G$
\cite[20.- Corollaire on p.149]{Ghys-Harpe(1990)}
the map $L \mapsto g_L$ is injective.
Because there are only finitely many conjugacy classes
of translation length $\leq C$ this means that
$G \backslash \call$ is finite.
This is what we needed to prove.
\end{proof}


\def\cprime{$'$} \def\polhk#1{\setbox0=\hbox{#1}{\ooalign{\hidewidth
  \lower1.5ex\hbox{`}\hidewidth\crcr\unhbox0}}}

\end{document}